\documentclass[final,leqno,letterpaper]{etna}

\usepackage[T1]{fontenc}
\usepackage{hyphenat} 
\usepackage{microtype}

\usepackage{amsmath} 
\usepackage{amsfonts,amsopn}
\usepackage{nicefrac}
\usepackage{bm} 

\usepackage{siunitx} 
\usepackage{epstopdf}
\usepackage{graphicx}
\usepackage{tikz}
\usepackage{algorithmic}
\usepackage{enumitem}

\usepackage{booktabs,multirow,colortbl,tabularx}
\renewcommand{\arraystretch}{1.05}
\definecolor{tblHighlGray}{RGB}{194,194,194}
\newcolumntype{L}{@{\extracolsep{\fill}}l}
\newcolumntype{R}{@{\extracolsep{\fill}}r}
\newcolumntype{C}{@{\extracolsep{\fill}}c}

\usepackage{listings}
\definecolor{xmlStringColor}{RGB}{169,39,55}
\lstdefinestyle{mystyleXML}{
    commentstyle=\color{black!65!white}\itshape,
    stringstyle=\color{xmlStringColor},
    basicstyle=\ttfamily\footnotesize ,
    breakatwhitespace=true,         
    breaklines=true,                 
    captionpos=t,                    
    keepspaces=true,                 
    numberstyle=\tiny,
    numbers=left,
    frame=lines,                  
    numbersep=5pt,                  
    showspaces=false,                
    showstringspaces=false,
    showtabs=false,                  
    tabsize=2,
    columns=flexible,
    aboveskip=1pt,    
    belowskip=1pt     
}
\usepackage[sort]{cite} 

\usepackage{cleveref}
\crefname{equation}{}{}  
\crefname{lstlisting}{listing}{listings}
\Crefname{lstlisting}{Listing}{Listings}

\usepackage{url}
\newcommand*{\email}[1]{\href{mailto:#1}{\nolinkurl{#1}} } 

\ifpdf
  \DeclareGraphicsExtensions{.eps,.pdf,.png,.jpg}
\else
  \DeclareGraphicsExtensions{.eps}
\fi

\DeclareMathOperator{\Jacobian}{\mathcal{J}}
\DeclareMathOperator{\Residual}{\mathcal{R}}
\DeclareMathOperator{\AdditiveSchwarz}{\mathbf{A_S}}
\newcommand{\Ri}[0]{\mathcal{R}_i}
\newcommand{\TildeRiT}[0]{\widetilde{\mathcal{R}}_i^T}
\newcommand{\expStar}[0]{$^{\star}$}
\newcommand{\CFL}[0]{c_{\text{CFL}}}
\newcommand{\CFLinT}[0]{c_{\text{CFL},T}}
\newcommand{\CFLavg}[0]{c_{\text{CFL},\text{avg}}}
\newcommand{\dx}[0]{\,\mathrm{d}x}
\newcommand{\Reynum}[0]{\textrm{Re}}

\def\fullTitle{Monolithic and Block Overlapping Schwarz Preconditioners\\ for the Incompressible Navier--Stokes Equations}
\def\shortTitle{Monolithic and Block Overlapping Schwarz Preconditioners}
\def\authorListShortA{A. Heinlein, A. Klawonn, J. Knepper, and L. Saßmannshausen} 
\def\authorListShortB{A. Heinlein, A. Klawonn, J. Knepper, L. Saßmannshausen} 

\title{\fullTitle}

\shorttitle{\shortTitle} 
\shortauthor{\authorListShortA}

\author{Alexander Heinlein\footnotemark[1]
\and Axel Klawonn\footnotemark[2]\ \footnotemark[3]
\and Jascha Knepper\footnotemark[2]\ \footnotemark[3]
\and Lea Saßmannshausen\footnotemark[2]\ \footnotemark[3]
}
\renewcommand{\thefootnote}{\fnsymbol{footnote}}

\ifpdf
\hypersetup{
  pdftitle={\fullTitle},
  pdfauthor={\authorListShortB}
}
\fi

\begin{document}

\footnotetext[1]{Delft Institute of Applied Mathematics, Delft University of Technology, Mekelweg 4, Delft, 2628 CD, Netherlands (\email{a.heinlein@tudelft.nl}, \url{https://searhein.github.io})}
\footnotetext[2]{Department of Mathematics and Computer Science, University of Cologne, Weyertal 86--90, Cologne, 50931, Germany (\email{axel.klawonn@uni-koeln.de}, \email{jascha.knepper@uni-koeln.de}, \email{l.sassmannshausen@uni-koeln.de}, \url{https://numerik.uni-koeln.de})}
\footnotetext[3]{Center for Data and Simulation Science, University of Cologne, Albertus-Magnus-Platz, Cologne, 50923, Germany (\url{https://cds.uni-koeln.de})}


\renewcommand{\thefootnote}{\arabic{footnote}} 

\maketitle

\begin{abstract}
Monolithic preconditioners applied to the linear systems arising during the solution of the discretized incompressible Navier--Stokes equations are typically more robust than preconditioners based on incomplete block factorizations. 
Lower number of iterations and a reduced sensitivity to parameters like velocity and viscosity can significantly outweigh the additional cost for their setup. 
Different monolithic preconditioning techniques are introduced and compared to a selection of block preconditioners. 
In particular, two-level additive overlapping Schwarz methods (OSM) are used to set up monolithic preconditioners and to approximate the inverses arising in the block preconditioners. 
GDSW-type (Generalized Dryja--Smith--Widlund) coarse spaces are used for the second level. 
These highly scalable, parallel preconditioners have been implemented in the solver framework \texttt{FROSch} (Fast and Robust Overlapping Schwarz), which is part of the software library \texttt{Trilinos}. 
The new GDSW-type coarse space GDSW\expStar{} is introduced; 
combining it with other techniques results in a robust algorithm.
The block preconditioners PCD (Pressure Convection--Diffusion), SIMPLE (Semi-Implicit Method for Pressure Linked Equations), and LSC (Least-Squares Commutator) are considered to various degrees. 
The OSM for the monolithic as well as the block approach allows the optimized combination of different coarse spaces for the velocity and pressure components, enabling the use of tailored coarse spaces.
The numerical and parallel performance of the different preconditioning methods for finite element discretizations of stationary as well as time-dependent incompressible fluid flow problems is investigated and compared. 
Their robustness is analyzed for a range of Reynolds and Courant--Friedrichs--Lewy (CFL) numbers with respect to a realistic problem setting.
\end{abstract}




\section{Introduction}

Mixed finite element discretizations of the incompressible Navier--Stokes equations lead to systems of saddle-point structure with the unknowns consisting of velocity and pressure variables.
We introduce several monolithic preconditioning techniques and show comparisons with a selection of block preconditioners, where one of the key ingredients is to find a suitable approximation of the Schur complement arising from a block factorization.
Particularly, two-level additive overlapping Schwarz methods are used to set up monolithic preconditioners and to approximate the inverses that appear in the block preconditioners. 
For the first time, this work combines different types of coarse spaces for the velocity and pressure components of the two-level overlapping Schwarz methods, which improves robustness for the monolithic preconditioner. 
By varying Reynolds and CFL numbers in a well-known benchmark problem as well as in a realistic setting, robustness of the monolithic and block overlapping Schwarz approaches is studied and compared. 

We consider the transient and stationary case of the Navier--Stokes equations. 
Newton's method is used for linearization and the finite element method is used to obtain a discrete saddle-point system. 
The discretization with respect to the velocity and pressure space is based on P2--P1, Q2--Q1, Q2--P1-disc., or stabilized P1--P1 or Q1--Q1 elements. 
We solve the linearized saddle-point system using the Generalized Minimal Residual Method (GMRES) Krylov subspace method \cite{saad:1986:gmres}. 
The arising linear system is ill-conditioned, which motivates the use of suitable preconditioners.

To compare the preconditioners, we show weak scalability results, vary the choice of coarse spaces for the second level of the Schwarz method, and offer recycling strategies for coarse problem components (for example, the reuse of symbolic factorizations). 
Furthermore, we show for a range of Reynolds and CFL numbers results for different finite element discretizations and structured and unstructured geometries and domain decompositions. 
Finally, we give strong scalability results for one problem.

Various approaches to build efficient discretizations and solvers for the incompressible Navier--Stokes equations exist, for example, 
\cite{gauthier:2004:fpi,
li:2006:bddc,
ELMAN2008_compParallelBlock,
benzi:2011:rdf,
griffith:2009:projection,
benzi:2013:pim,
cyr_teko_2016,
he:2017:bpi,
hanek:2020:mbddc,
bevilacqua:2024:bddc,
klawonn:2024:nmt}. 
Due to the vast number of different possibilities, we limit ourselves to a selection of discretizations, solution processes, and preconditioners. 
Compared to monolithic preconditioning approaches, there exists extensive research on block preconditioners; see, for example, 
\cite{ELMAN2008_compParallelBlock,
Klawonn_Starke_nonsymSaddlepoints,KLA:1998:BTP,KLA:1998:AOP,
Silvester_PCD_2001,
Klawonn_2000_Mono_Block,
Elman_OG_BFBT,
Elman_fluid_book,
Elman_2001_PrecNS_Mu,
Kay_PCD_2002,
Prec_Twophase_NS_2019,
DEPARIS_unsteady_NS_PCD} 
for a variety of block-preconditioning approaches in the literature.

Since the assembly of the exact Schur complement is usually too expensive, finding a suitable approximation of the Schur complement for the construction of block preconditioners is an integral task; see \Cref{section: Block Preconditioner}.
In this work, we consider different Schur complement approximations. 
For the application of the resulting block preconditioners, the action of the inverse of block matrices is required. 
These inverses are further approximated with a two-level additive overlapping Schwarz method \cite{smith:1996:dpm,toselli:2005:ddm}.

The definition of a suitable approximation is highly problem dependent. 
For a Stokes problem, a different preconditioner can be used than for a Navier--Stokes problem, particularly the nonlinearity -- influenced by the Reynolds number -- plays a vital role. 
The preconditioning approach can vary if the problem is stationary or time-dependent, if Newton's method is used for linearization or a Picard iteration, if it is compressible flow or incompressible, Newtonian or non-Newtonian, and also the choice of the time and space discretization matters.

When dealing with a steady-state Stokes problem and a stable discretization, the Schur complement can be replaced with a spectrally equivalent scaled pressure mass matrix; see \cite{Elman_2001_PrecNS_Mu,Klawonn_Starke_nonsymSaddlepoints}. 
Consequently, for low Reynolds numbers, we may expect that the pressure mass matrix achieves sufficiently good results. 
For higher Reynolds numbers, the pressure mass matrix cannot account for the dominating advective forces.
Therefore, for larger advective forces, a different approximation of the Schur complement should be used. 

The considered block methods are the Pressure Convection--Diffusion (PCD) \cite{Silvester_PCD_2001, Kay_PCD_2002, Elman_fluid_book, prec_nav_stokes, Elman_BC_Navier_Stokes, PCD_bc_Tuminaro, DEPARIS_unsteady_NS_PCD}
and Least-Squares Commutator (LSC) \cite{Elman_OG_BFBT, Elman_fluid_book, prec_nav_stokes, Elman_LSC_2006,Elman_BC_Navier_Stokes} triangular block preconditioners and the SIMPLE (Semi-Implicit Method for Pressure Linked Equations) \cite{Patankar_SIMPLE_1972,  SIMPLE_pernice_multigrid-preconditioned_2001, ELMAN2008_compParallelBlock, prec_nav_stokes, DEPARIS_unsteady_NS_PCD} preconditioner. 
In this work, we restrict ourselves to the SIMPLE variant SIMPLEC and note that variations like SIMPLER exist; 
see, for example, \cite{SIMPLER_OG_Doorman,prec_nav_stokes, simple-like_2016}.
These methods each pursue a different approach in finding a suitable approximation of the Schur complement. 
In the \texttt{Trilinos} package \texttt{Teko}, the desired block preconditioners are predefined and can also be combined with other techniques (e.g., algebraic multigrid) for individual blocks \cite{Cyr_teko_stabilized, cyr_teko_2016}. 

In a monolithic approach, the preconditioner is built for the entire saddle-point system \cite{Klawonn_1998_MixedFE, Klawonn_2000_Mono_Block, Heinlein_2019_MonoFluidFlow, heinlein_reduced_2019, diss_hochmuth_2020}. 
Consequently, a Schur complement approximation as for block preconditioners is not required. 
Both preconditioning techniques could be extended to handle related problems, such as fluid--structure interaction.

An area of application is the simulation of realistic hemodynamics. 
This can entail varying viscosities and high velocities leading to Reynolds numbers up to a few thousand; 
see also \cite{DEPARIS_unsteady_NS_PCD} for parallel preconditioners in the context of hemodynamics. 
A preconditioner that is robust to changes in the velocity and Reynolds number -- for example due to heartbeats -- is advantageous. 
An additional emphasis is placed on achieving parallel scalability.

The remainder of the paper is organized as follows. 
In \Cref{Sec: Navier-Stokes Equations} we introduce the Navier--Stokes equations along with the boundary conditions used in this work, describe the discretization in space and time and the Newton--Krylov method used to solve the nonlinear problem. 
In \Cref{section: Block Preconditioner} we derive the block preconditioners that are used later to obtain numerical results. 
The inverse of the monolithic system and inverses that arise during the construction of block preconditioners are approximated with additive overlapping Schwarz preconditioners; 
a detailed description is given in \Cref{sec:schwarz}. 
Different coarse spaces that can be used for the coarse level of the overlapping Schwarz method are specified in \Cref{subsec: Coarse Spaces}. 
Specifically, we make use of the coarse spaces GDSW \cite{dohrmann_domain_2008,dohrmann_family_2008}, RGDSW \cite{dohrmann_design_2017,heinlein_reduced_2019}, and introduce the new variant GDSW\expStar{}. 
A description of the parallel implementation is given in \Cref{Sec: Implementation}. 
It consists of the \texttt{FEDDLib} \cite{feddlib}, which is used for the assembly and as a simulation tool that interfaces \texttt{Trilinos}, and the \texttt{Trilinos} packages for block and monolithic preconditioning: \texttt{Teko} is used for block preconditioners and \texttt{FROSch} contains overlapping Schwarz preconditioners used to approximate the inverse of the monolithic system and the inverses arising in the block preconditioners. 
\Cref{Sec: Results} contains results for different problem settings. 
We consider a backward-facing step geometry for stationary and transient fluid flow and a realistic artery. 
In both test cases, we consider Reynolds numbers between 200 and 3\,200 and a range of CFL numbers to test the robustness of the solvers. 
The results indicate that, using a combination of the new coarse space GDSW\expStar{} for the velocity with RGDSW for the pressure, our monolithic preconditioner is more robust than the considered block preconditioners PCD, LSC, and SIMPLEC. 
In a realistic setting, varying the hemodynamic properties via the viscosity or Reynolds number, the monolithic approach performs best. 
Specifically, the total number of linear iterations is lower. 
Despite a higher cost to set up the monolithic preconditioner, the total compute time is usually lower; 
particularly for high CFL numbers, the difference is significant.


\section{Navier--Stokes equations}\label{Sec: Navier-Stokes Equations}

We consider the transient, incompressible Navier--Stokes equations to model a Newtonian fluid given by
\begin{equation}\label{eq: transient Navier-Stokes}
\begin{split}
	\rho\Big (\frac{\partial \mathbf{u}}{\partial t} + (\mathbf{u}\cdot \nabla) \mathbf{u} \Big )  - \mu \Delta \mathbf{u} + \nabla p &= \rho \mathbf{f}, \\
	\operatorname{div}(\mathbf{u} ) &= 0, 
	\end{split}
\end{equation}
with velocity $\mathbf{u}$, pressure $p$, density $\rho$, dynamic viscosity $\mu$, and volume force $\mathbf{f}$. 
The kinematic viscosity is given by $\nu=\nicefrac{\mu}{\rho}$. 
If the set of equations \cref{eq: transient Navier-Stokes} reaches a steady state for $t\to\infty$, we have $\partial_t u = 0$ and obtain the stationary Navier--Stokes equations
\begin{equation}\label{eq: stationary Navier-Stokes}
	\begin{split}
	\rho (\mathbf{u}\cdot \nabla) \mathbf{u}  - \mu \Delta \mathbf{u} + \nabla p &= \rho \mathbf{f}, \\
	\operatorname{div}(\mathbf{u} ) &= 0.
	\end{split}
\end{equation}
We denote the fluid domain by $\Omega\subset\mathbb{R}^{3}$, and its boundary by $\partial \Omega$. We assume that there are subsets $\partial \Omega_{\text{in}}$, $\partial \Omega_{\text{out}} \subset \partial \Omega$, which denote the inflow and outflow parts of the boundary. 
We assume that $\partial \Omega_{\text{in}}$ and $\partial \Omega_{\text{out}}$ have a positive surface measure. 
On $\partial \Omega_{\text{in}}$ a velocity field is prescribed for $\mathbf{u}$; 
on $\partial \Omega_{\text{out}}$ the Neumann boundary condition
\begin{equation*}
\nu \nabla\mathbf{u}\cdot\mathbf{n} - p\mathbf{n} = 0 \quad \text{ on }  \partial \Omega_{\text{ out }}
\end{equation*}
is prescribed. 
The other boundaries of the fluid domain are walls, $\partial \Omega_{\text{wall}}\subset \partial \Omega$, which prevent fluid from leaving the computational domain. 
We use a no-slip condition, such that $\mathbf{u} = 0$ on $\partial \Omega_{\text{wall}}$.
We have 
\begin{equation*}
\partial \Omega = \partial \Omega_{\text{in}} \cup \partial \Omega_{\text{out}} \cup \partial \Omega_{\text{wall}}
\end{equation*}
and 
\begin{equation*}
\partial \Omega_{\text{in}} \cap \partial \Omega_{\text{wall}} = \emptyset,
\quad\partial \Omega_{\text{out}} \cap \partial \Omega_{\text{wall}} = \emptyset,
\quad\partial \Omega_{\text{in}} \cap\partial  \Omega_{\text{out}} = \emptyset.
\end{equation*}


\subsection{Finite element and time discretization} \label{sec:disc}

For the spatial discretization of the incompressible fluid flow problem, we use different pairs of mixed finite elements. 
We partition the domain with characteristic mesh size $h$ into tetrahedra or hexahedra. 
Let $\tau_h:=\tau_h(\Omega)$ denote the triangulation, $\mathcal{P}_k(T)$ the polynomial space of order $k$, and $\mathcal{Q}_k(T)$ the polynomial space with order up to $k$ in each coordinate direction. 
For the partition into tetrahedral elements, we use continuous and either piecewise linear or quadratic functions for the velocity space
\begin{equation*}
	\mathbf{V}^{h}_{Pk} (\Omega) 
	=
	\Big  \{ 
		\mathbf{v}_h \in  \big(C(\Omega)\big)^{3} \cap \big(H^{1}(\Omega)\big)^{3} 
		: 
		\mathbf{v}_h|_{T} \in \big(\mathcal{P}_k(T)\big)^3\ \forall T \in \tau_h
	\Big\},\ k\in\{1,2\}.
\end{equation*}
Similarly, we can define the velocity space for hexahedral elements,
\begin{equation*}
	\mathbf{V}^{h}_{Qk} (\Omega) 
	=\Big \{ 
		\mathbf{v}_h \in \big(C(\Omega)\big)^{3} \cap \big(H^{1}(\Omega)\big)^{3} 
		: 
		\mathbf{v}_h|_{T} \in \big(\mathcal{Q}_k(T)\big)^3\ \forall T \in \tau_h 
	\Big\},\ k\in\{1,2\}. 
\end{equation*}
The continuous and piecewise linear pressure spaces for the different partitions are defined as
\begin{align*} 
Q^h_{P1}  (\Omega) &=\Big \{ q_h \in C(\Omega) \cap L^2(\Omega): q_h|_{T} \in \mathcal{P}_1(T)\ \forall T \in \tau_h\Big \},\\ 
Q^h_{Q1}  (\Omega) &=\Big \{ q_h \in C(\Omega) \cap L^2(\Omega): q_h|_{T} \in \mathcal{Q}_1(T)\ \forall T \in \tau_h\Big \}.
\end{align*}
Additionally, we consider the combination of Q2--P1-disc., where the pressure space consists of discontinuous piecewise linear functions:
\begin{equation*}
Q^h_{\text{P1-disc.}} (\Omega)  =\Big \{ q_h \in  L^2(\Omega): q_h|_{T} \in \mathcal{P}_1(T)\ \forall	T \in \tau_h\Big \}.
\end{equation*}
If the type of element is unspecified, we write $\mathbf{V}^h$ and $Q^h$ for the velocity and pressure space, respectively. 
Let $\mathbf{u}_h \in \mathbf{V}^{h}$ and $p_h \in Q^h$ be the approximations of the solutions $\mathbf{u}$ and $p$. 
Then the resulting variational formulation for the stationary case \cref{eq: stationary Navier-Stokes} is as follows: Find $\mathbf{u}_h \in \mathbf{V}^{h}$ and $p_h \in Q^h$ such that 
\begin{align*}
\nu \!\int_\Omega \!\!\nabla \mathbf{u}_h \!:\!\nabla \mathbf{v}_h\dx 
+\!\! \int_\Omega \!((\mathbf{u}_h \cdot \nabla) \mathbf{u}_h) \cdot \mathbf{v}_h \dx 
-\!\! \int_\Omega \! p_h (\nabla \cdot \mathbf{v}_h)\dx 
&= \!\! \int_\Omega \!\mathbf{f} \cdot \mathbf{v}_h \dx 
&\forall \mathbf{v}_h \in \mathbf{V}^{h}, \\
\int_{\Omega} \! q_h (\nabla \cdot \mathbf{u}_h)\dx 
&= 0  
& \forall q_h \in Q^h.
\end{align*}
We denote by $\bm{\varphi}_1,\dots,\bm{\varphi}_{N_u}$ the finite element basis of $\mathbf{V}^{h}$ and by $\psi_1,\dots,\psi_{N_p}$ the basis of $Q^h$. Then, we can write
\begin{equation*}
\mathbf{u}_h(x) = \sum_{i=1}^{N_u} u_i \bm{\varphi}_i(x)\ \text{ and }\ p_h(x) = \sum_{i=1}^{N_p} p_i \psi_i(x).
\end{equation*}
We will also refer to the corresponding coefficient vectors of $\mathbf{u}_h(x)$ and $p_h(x)$ by $\mathbf{u}_h$ and $p_h$, respectively. 
We define the following matrix operators corresponding to the weak formulation,
\begin{align*}
A &= [a_{ij}], 
& \! B&=[b_{ij}], 
& \! N(\mathbf{u}_h) &= [n_{ij}(\mathbf{u}_h)], \\
a_{ij} &= \int_\Omega  \nabla \bm{\varphi}_i : \nabla \bm{\varphi}_j \dx, 
& \! b_{ij} &= -\!\int_\Omega \psi_i (\nabla \cdot \bm{\varphi}_j) \dx, 
& \! n_{ij}(\mathbf{u}_h) &= \int_\Omega ((\mathbf{u}_h \cdot \nabla) \bm{\varphi}_j)\cdot \bm{\varphi}_i\dx,
\end{align*}
and obtain with the discretization $\mathbf{f}_u$ of $\mathbf{f}$ the system of linear equations for the stationary Navier--Stokes equations,
\begin{equation}\label{Eq: NS Steady nonlinear}
\begin{bmatrix}
F_O(\mathbf{u}_h) & B^{T} \\
  B        & 0
\end{bmatrix}
\begin{bmatrix}
\mathbf{u}_h \\
p_h
\end{bmatrix}
=
\begin{bmatrix}
\mathbf{f}_u \\
0
\end{bmatrix},
\end{equation}
with $F_O(\mathbf{u}_h) = \nu  A +  N(\mathbf{u}_h)$. 
The resulting Picard iteration with the iterates $\mathbf{u}_h^k$ and $p_h^k$ reads \cite[p.~346]{Elman_fluid_book}
\begin{equation*}
\begin{bmatrix}
F_O(\mathbf{u}_h^k) & B^{T} \\
  B        & 0
\end{bmatrix}
\begin{bmatrix}
\mathbf{u}_h^{k+1} \\
p_h^{k+1}
\end{bmatrix}
=
\begin{bmatrix}
\mathbf{f}_u \\
0
\end{bmatrix},
\end{equation*}
defining the discretized Oseen problem. 
We apply the Newton method to solve the nonlinear problem:
\begin{equation*}
\begin{aligned}
&\text{Solve: } & \Jacobian(X^k) \delta X^{k} &= -\Residual(X^k), \\
&\text{Update: } & X^{k+1} &= X^k + \delta X^k,
\end{aligned}
\end{equation*}
where the $\delta$ notation is defined as $\delta X^k := X^{k+1}- X^k$. 
The Newton residual $\Residual$ for the stationary case is defined as
\begin{equation*}
\Residual_{\text{stat}}
\left(
\mathbf{u}_h^k,
p_h^k
\right)
= 
\begin{bmatrix}
F_O(\mathbf{u}_h^k) & B^{T} \\
  B        & 0
\end{bmatrix}
\begin{bmatrix}
\mathbf{u}_h^{k} \\
p_h^{k}
\end{bmatrix}
-
\begin{bmatrix}
\mathbf{f}_u \\
0
\end{bmatrix}.
\end{equation*}
We need to linearize \cref{Eq: NS Steady nonlinear} by adding an additional term $W(\mathbf{u}_h^k)$ to $F_O(\mathbf{u}_h^k)$ from the derivative of the convection term:
\begin{equation*}
W(\mathbf{u}_h^k)
= [w_{ij}(\mathbf{u}_h^k)],
\quad
w_{ij}(\mathbf{u}_h^k)
= \int_\Omega ((\bm{\varphi}_j \cdot \nabla) \mathbf{u}_h^{k})\cdot \bm{\varphi}_i \dx,
\end{equation*}
with the velocity solution from the previous $k$th Newton step $\mathbf{u}_h^{k}$. 
This leads to 
\begin{equation*}
\Jacobian
\left(
\mathbf{u}_h^k
\right)
\begin{bmatrix}
\delta \mathbf{u}_h^{k} \\
\delta p_h^{k}
\end{bmatrix} = 
\begin{bmatrix}
F_N(\mathbf{u}_h^{k})  & B^{T} \\
  B & 0 
\end{bmatrix}
\begin{bmatrix}
\delta \mathbf{u}_h^{k} \\
\delta p_h^{k}
\end{bmatrix}
=
-
\Residual_{\text{stat}}
\left(
\mathbf{u}_h^k,
p_h^k
\right),
\end{equation*}
with $F_N(\mathbf{u}_h^k) = F_O(\mathbf{u}^k_h) + W(\mathbf{u}_h^{k})$.
For transient problems we include the BDF-2 time discretization \cite{BDF2} and expand the system to 
\begin{equation}\label{Eq: Sys with BDF2}
    \Jacobian\!
    \left(
    \mathbf{u}_h^{k,n+1}
    \right)
    \!
    \begin{bmatrix}
    \delta \mathbf{u}_h^{k,n+1} \\
    \delta p _h^{k,n+1}
    \end{bmatrix}
=
    -\Residual_{\text{stat}}\!
    \left(
    \mathbf{u}_h^{k,n+1},
    p_h^{k,n+1}
    \right)
    + 
    \begin{bmatrix}
    \frac{2}{\Delta t} M_u \mathbf{u}_h^{n} - \frac{1}{2\Delta t}M_u \mathbf{u}_h^{n-1}\\
    0
    \end{bmatrix}\!,
\end{equation}
where $M_u$ is the velocity mass matrix, $\mathbf{u}_h^n$ and $\mathbf{u}_h^{n-1}$ are the solutions to the last two time steps, and 
\begin{equation*}
    \Jacobian\!
    \left(
    \mathbf{u}_h^{k,n+1}
    \right)
    =
    \begin{bmatrix}
    \frac{3}{2\Delta t_n} M_u +F_N(\mathbf{u}_h^{k,n+1}) & B^{T} \\
    B & 0 
    \end{bmatrix}.
\end{equation*}
The resulting saddle-point problem for the different discretizations has the generic form
\begin{equation}\label{Eq: Saddle point problem Navier Stokes}
	\mathcal{F} := \begin{bmatrix}
	F & B^{T} \\ B & -C 
	\end{bmatrix},
\end{equation}
where the fluid component $F$ depends on the underlying problem setting, and $C\neq 0$ in case of a discretization that is not inf-sup stable.

We focus on three-dimensional problems for which the P2--P1, Q2--Q1, and Q2--P1-disc. finite element pairs satisfy the inf--sup condition; cf. \cite{Boffi_MixedFE_2013}.
In case of the non inf-sup stable finite elements P1--P1 and Q1--Q1, the Bochev--Dohrmann stabilization \cite{Dohrmann_stabilization_2004} is used, which penalizes nonphysical pressure variations by inserting a pressure term into the conservation of mass equation (see \cite[Sect.~2]{Dohrmann_stabilization_2004} and \cite[Sect. 5.6.1]{diss_hochmuth_2020}). 
Its bilinear form reads
\begin{equation*}
C(p,q) = \frac{1}{\nu} \int_{\Omega} (p-\rho_0 p) ( q-\rho_0 q)\dx \quad \forall q \in L^2(\Omega),
\end{equation*}
where $\rho_0$ is an $L^2$-orthogonal projection onto the space $\mathcal{P}_0(\tau_h)$ of discontinuous, piecewise constant functions;
the projection is given by 
\begin{equation*}
(\rho_0 q)|_T = \left(\int_T 1\dx\right)^{-1} \int_T q \dx\quad\forall T\in\tau_h(\Omega),
\end{equation*}
which can be computed efficiently due to its independent definition on each element $T\in\tau_h(\Omega)$.

As a measure that indicates how much information travels through a computational grid cell in a unit of time, that is, which relates the velocity of transport, time step, and mesh resolution, we define a Courant(--Friedrichs--Lewy)-type number (CFL) for structured meshes as 
\begin{equation*}
\CFL(t) 
= \frac{u^h_\infty(t) \Delta t}{h}, 
\quad u^h_\infty(t)
=\max_{x^h\in\tau_h(\Omega)} \|\mathbf{u}_h(x^h,t)\|_2,
\end{equation*}
where $x^h$ is a finite element node. 
For an unstructured mesh, we take varying element sizes and shapes into account and use 
\begin{align*}
\CFL(t)
&= \max_{T\in\tau_h(\Omega)} \CFLinT(t), \\
	 \CFLinT(t) &= \frac{u^h_\infty(T,t) \Delta t}{d_I(T)}\quad\qquad \text{ with }
	 u^h_\infty(T,t) = \max_{x^h\in T} \|\mathbf{u}_h(x^h,t)\|_2,
\end{align*}
where $T\in\tau_h(\Omega)$ is a finite element, and $d_I(T)$ is the diameter of the largest incircle of~$T$. 

We will also use the average of $\CFLinT(t)$ over all elements $T\in\tau_h(\Omega)$, denoted by $\CFLavg(t)$. 
By $\CFL$ and $\CFLavg$, we refer to the corresponding maximum CFL number over time.
In case the CFL number is too large, instability or unphysical behavior may be encountered. 
Generally, implicit time stepping schemes are more stable with respect to the CFL number. 
High CFL numbers tend to be more demanding for preconditioners; see, for example, \cite[Fig.~9]{Cyr_teko_stabilized} and \Cref{figure: Iterations CFL comparison,figure: Iterations RE transient}.


\subsection{Inexact Newton method: Newton--Krylov}\label{Sec: Inexact Newton + Forcing Term}

We apply an inexact Newton method \cite{InexactNewton} by solving the linearized system with a Krylov subspace method up to a certain tolerance; 
specifically, we employ the generalized minimal residual (GMRES) method \cite{saad:1986:gmres}. 
With an initial guess for $\delta X^k$, we use GMRES to solve
\begin{equation*}
\Vert \Jacobian(X^k) \delta X^k + \Residual(X^k) \Vert_2 \leq \eta_k \Vert \Residual(X^k)\Vert_2,
\quad \eta_k \in [0,1),
\end{equation*}
and update the Newton iterate via $X^{k+1} = X^{k} + \delta X^{k}$. 
The linear solver tolerance $\eta_k$ can be constant throughout the Newton iterations or change between consecutive iterations. 
In the second case, $\eta_k$ is referred to as forcing term~\cite{ForcingTermNewton}, which can change depending on the underlying strategy. 
We use the \textit{Choice 2} forcing term from~\cite{ForcingTermNoxImplementation},
\begin{equation*}
\eta_k = \gamma \Bigg( \frac{\Vert \Residual(X^{k}) \Vert_2}{\Vert \Residual(X^{k-1}) \Vert_2} \Bigg)^{\alpha},
\end{equation*}
with $\gamma \in [0,1]$ and $\alpha \in (1,2]$. To make use of this strategy, we select an $\eta_0 \in [0,1)$ for the first iteration. We impose a lower and upper bound for the forcing term by
\begin{equation}\label{eq: forcing term bounds}
\max\Big \{\gamma \eta_{k-1}^{\alpha},\eta_{\min} \Big \} \leq \eta_k \leq \eta_{\max},
\end{equation}
with prescribed bounds $\eta_{\min}$ and $\eta_{\max}$.


\section{Block preconditioners}\label{section: Block Preconditioner}

We introduce a block LDU decomposition of the Navier--Stokes saddle-point problem~\cref{Eq: Saddle point problem Navier Stokes},
\begin{equation}\label{eq:Navier Stokes LDU}
\mathcal{F} 
= 
	\begin{bmatrix}
		F & B^{T} \\ B & -C 
	\end{bmatrix}
=
	\begin{bmatrix}
	 I  &  0 \\
	 B F^{-1} &  I
	\end{bmatrix}
	\begin{bmatrix}
	 F  &  0 \\
	 0 & S
	\end{bmatrix}
	\begin{bmatrix}
	 I  &   F^{-1} B^{T} \\
	0 &  I
	\end{bmatrix},
\end{equation}
where $S= -C - B F^{-1}B^{T}$ is the Schur complement. 
Generally, we seek to replace $F$ and the Schur complement $S$ with suitable approximations $\hat{F}$ and $\hat{S}$ that facilitate the computation of $\hat{F}^{-1}$ and $\hat{S}^{-1}$. 
Suitable approximations depend on the problem setting. 
The approximation of $F^{-1}$ will be based on a Schwarz preconditioner (see \Cref{sec:schwarz}), and the approximation of $S^{-1}$ on the approaches LSC, PCD, SIMPLE, and SIMPLEC. 
These preconditioners contain inverses of matrices that will be further approximated with a Schwarz method. 
A resulting upper block-triangular preconditioner reads 
\begin{equation}
\label{eq:block-tri precond}
\mathcal{B}_{\text{Tri}}^{-1}
= 
	\begin{bmatrix}
	\hat{F} & B^{T} \\
	0 & \hat{S} 
	\end{bmatrix} ^{-1}
= 
	\begin{bmatrix}
	\hat{F}^{-1} & - \hat{F}^{-1} B^{T} \hat{S}^{-1} \\
	0 & \hat{S}^{-1} 
	\end{bmatrix},
\end{equation}
which encapsulates the process of a block-Gaussian elimination.

For a Navier--Stokes problem with kinematic viscosity $\nu$, velocity $v$, and a characteristic length $L$ (for example, the diameter of a vessel), the Reynolds number is defined as 
\begin{equation}\label{eq:Reynolds number:general}
\Reynum = \frac{v L }{\nu}.
\end{equation}
For steady-state Stokes problems and a stable discretization, the Schur complement can be replaced with the spectrally equivalent, scaled pressure mass matrix $-\frac{1}{\nu} M_p$ \cite{Elman_2001_PrecNS_Mu,Klawonn_Starke_nonsymSaddlepoints}. 
We may, thus, expect that this approach also works well for small Reynolds numbers. 
For higher Reynolds numbers, the pressure mass matrix cannot account for the dominating advective forces in $F$.
Therefore, for larger advective forces, a different approximation of the Schur complement should be used. 
We will present the PCD (Pressure Convection--Diffusion) and LSC (Least-Squares Commutator) block-triangular preconditioners and the SIMPLE method (Semi-Implicit Method for Pressure Linked Equations) in the following sections.

The PCD and LSC Schur complement approximations are motivated by finding a suitable approximation of $F^{-1}$, which is used to construct the Schur complement. 
Here and in the following, for the construction of approximations of $F^{-1}$, we assume that $F=F_O$, the discrete Oseen operator. 
In an implementation, however, we will use the general $F$ in \cref{Eq: Saddle point problem Navier Stokes} that includes an additional term from the Newton linearization. 
For the consideration of the time discretization, see \cref{Eq: Commutator_time}.

The discrete Oseen operator is the discretization of a convection--diffusion differential operator for the velocity, 
\begin{equation*}
\mathcal{L} := - \nu \Delta + \mathbf{w}_h \cdot \nabla,
\end{equation*}
where $\mathbf{w}_h$ is the approximation to the velocity solution from the previous nonlinear iteration. 
The discretization of $\mathcal{L}$ corresponds to $F(\mathbf{w}_h)=\nu A + N(\mathbf{w}_h)$. 
We denote the corresponding differential operator on the pressure space by
\begin{equation*}
\mathcal{L}_p := (-\nu \Delta + \mathbf{w}_h \cdot \nabla)_p.
\end{equation*}
Its discretization is given by $F_p(\mathbf{w}_h) = \nu A_p + N_p(\mathbf{w}_h) = [f_{p,ij}(\mathbf{w}_h)]$ with 
\begin{equation*}
f_{p,ij}(\mathbf{w}_h)
= \nu \int_\Omega \nabla \psi_j \cdot \nabla \psi_i\dx
+\int_{\Omega} ( \mathbf{w}_h \cdot \nabla \psi_j) \psi_i \dx.
\end{equation*}
Using the convection--diffusion operators, we can define a commutator \cite[(9.13)]{Elman_fluid_book}
\begin{equation}\label{Eq: Commutator}
\mathcal{E} := \nabla \cdot (-\nu \Delta + \mathbf{w}_h \cdot \nabla ) - (-\nu \Delta + \mathbf{w}_h \cdot \nabla)_p \nabla \cdot 
\end{equation}
to find an approximation~$\hat{S}$ of the Schur complement~$S$. 
Note that it is assumed that the differential operators are applied to functions that are sufficiently regular; 
they are not well defined for the used finite element spaces. 

We suppose that the commutator \cref{Eq: Commutator} is small in some sense in order to define the PCD and LSC preconditioner.
The PCD Schur complement approximation is derived from a discretization of the commutator. 
In turn, to define LSC, the matrix $F_p$ is replaced with a new matrix that makes \cref{Eq: Commutator} small in a least-squares sense.

\Cref{Table: Teko and Monolithic} lists the matrices that are required to set up the different preconditioners that are introduced in the following sections; 
these matrices are usually further approximated via, for example, lumping or a Schwarz method.

\begin{table}[!tb]
\caption{\textbf{Operators appearing in the different preconditioning strategies.}
These operators are usually further approximated; 
for example, $F^{-1}$ can be approximated with a Schwarz method. 
Often the inverse of a mass matrix $M$ is replaced by $H_M^D = \operatorname{diag}(M)^{-1}$, the inverse of the diagonal entries of $M$. 
Then, an elaborate strategy (e.g., a Schwarz method) for the approximation of $M^{-1}$ is not necessary.}\label{Table: Teko and Monolithic}
\begin{tabular*}{\linewidth}{@{} LCCCCCC@{} }
\toprule
Operator                     & Monolithic & PCD        & SIMPLE/SIMPLEC    & LSC        & LSC$_{A_p}$ & LSC$_{\text{stab},A_p}$ \\ \midrule
$\mathcal{F}^{-1}$           & \checkmark &            &                   &            &             &            \\ \midrule
$F^{-1}$                     &            & \checkmark & \checkmark        & \checkmark & \checkmark  & \checkmark \\ \midrule
$A_p^{-1}$                   &            & \checkmark &                   &            & \checkmark  & \checkmark \\ \midrule
$M_p^{-1}$                   &            & \checkmark &                   &            &             &            \\ \midrule
$M_u^{-1}$                   &            &            &                   & \checkmark & \checkmark  & \checkmark \\ \midrule
$\rho_\text{eig}^{-1}D^{-1}$ &            &            &                   &            &             & \checkmark \\ \midrule
$(B M_u^{-1} B^{T})^{-1}$    &            &            &                   & \checkmark &             &            \\ \midrule
$(-C-B H_F B^{T})^{-1}$      &            &            & \checkmark        &            &             &            \\
\bottomrule
\end{tabular*}
\end{table}


\subsection{PCD (Pressure Convection--Diffusion)} \label{subsection: PCD}

We present a derivation following \cite{prec_nav_stokes} and \cite[chapter~9]{Elman_fluid_book} of the pressure convection--diffusion preconditioner, which was introduced in \cite{Silvester_PCD_2001,Kay_PCD_2002}.

In the velocity space with zero boundary conditions, the matrix representation of the discrete divergence operator is $-M_p^{-1} B$, and $M_u^{-1} B^{T}$ is the representation of the discrete gradient operator \cite[p.~409]{Elman_fluid_book}. 
$M_u$ is the velocity mass matrix and $M_p$ the pressure mass matrix:
\begin{equation*}
M_u = [m^u_{ij}],\quad 
m^u_{ij} = \int_{\Omega} \bm{\varphi}_i \cdot \bm{\varphi}_j \dx,\quad
M_p = [m^p_{ij}],\quad 
m^p_{ij} = \int_{\Omega} \psi_i \psi_j \dx.
\end{equation*}
With these representations, we can define a discrete version of the commutator~\cref{Eq: Commutator}:
\begin{equation}\label{eps_h}
\mathcal{E}_h = (M_p^{-1} B ) (M_u^{-1} F ) -  (M_p^{-1} F_p)(M_p^{-1} B).
\end{equation}
Suppose that the commutators $\mathcal{E}$ and $\mathcal{E}_h$ are small in some operator norm (cf. \cite[Remark~9.5]{Elman_fluid_book} and \Cref{subsection: LSC}).
Pre-multiplication of \cref{eps_h} with $M_p F_p^{-1}M_p$, post-multiplication with $F^{-1}B^{T}$, and assuming the commutator is small leads to the Schur complement approximation
\begin{equation}\label{eq: BFB approx}
-B F^{-1} B^{T} \approx   -M_p F_p^{-1}  B M_u^{-1} B^{T}.
\end{equation}
Note that $F$ is not symmetric and that, in this work, the approximation sign is not meant to denote spectral equivalence, but that the difference is approximately zero in some matrix norm.
We require the application of the inverse of the Schur complement, that is, also of its approximation $B M_u^{-1}B^{T}$. 
This expensive part is replaced with the pressure-Laplace operator $A_p$, which is spectrally equivalent for an enclosed-flow problem (cf. \cite[p.~366]{Elman_fluid_book} and \cite{prec_nav_stokes}): 
\begin{align*}
q_h^{T} A_p p_h
= \langle \nabla p_h,\nabla q_h \rangle_{L_2} 
&\approx \langle M_u^{-1} B^{T} p_h , M_u^{-1} B^{T} q_h \rangle_{M_u} \\
&= \langle B^{T} p_h , M_u^{-1} B^{T} q_h \rangle_{l_2} \\
&= q_h^{T} B M_u^{-1} B^{T} p_h.
\end{align*}
This gives us
\begin{equation}
\label{eq:PCD:Schur complement approximation}
S = -B F^{-1} B^{T} \approx - M_p F_p^{-1} A_p =: S_{\text{PCD}}.
\end{equation}
Then, we can define the PCD block-triangular preconditioner
\begin{equation*}
\mathcal{B}_{\text{PCD}}^{-1} = \begin{bmatrix}
F & B^{T} \\
0 & S_{\text{PCD}} 
\end{bmatrix}^{-1}.
\end{equation*}

As in \cite{Silvester_PCD_2001} we can extend the operator $F_p$ to also reflect the contribution of the BDF-2 time discretization for a transient Navier--Stokes problem:
\begin{equation}\label{F_p}
F_p(\mathbf{w}_h) = \frac{3}{2\Delta t} M_p + N_p(\mathbf{w}_h) + \nu A_p.
\end{equation}
The motivation is similar to before; 
the commutator is extended to
\begin{equation}\label{Eq: Commutator_time}
\mathcal{E} = \nabla \cdot \Big (\frac{\partial}{\partial t}-\nu \Delta + \mathbf{w}_h \cdot \nabla \Big) - \Big (\frac{\partial}{\partial t}-\nu \Delta + \mathbf{w}_h \cdot \nabla\Big )_p \nabla \cdot 
\end{equation}
to include a time derivative.

Note that the definition of $S_{\text{PCD}}$ is based on the assumption of a small error between the operators $\mathcal{L}$ and $\mathcal{L}_p$. 
This is only directly applicable to the fixed-point method, as Newton's method introduces an additional term for the linearization. 
Consequently for Newton's method, we cannot necessarily assume the error as defined above to be small. 
Nonetheless, it is observed that Newton's method combined with PCD preconditioning yields good results; see for example \cite[Chapters 9.3, 9.4]{Elman_fluid_book}.


\subsubsection{Stabilized discretizations}
\label{sssection:PCD:stabilized}

For elements which are not inf-sup stable, the PCD preconditioner also extends to stabilized approximations (e.g., P1--P1); see \cite[Chapter~9]{Elman_fluid_book}. 
A small commutator $\mathcal{E}$ (see \cref{Eq: Commutator}) implies
\begin{equation}\label{eq:stabilization pcd 1}
\nabla \cdot (-\nu \Delta + \mathbf{w}_h \cdot \nabla )^{-1} \nabla \approx (-\nu \Delta + \mathbf{w}_h \cdot \nabla)_p^{-1} \nabla \cdot \nabla .
\end{equation}
Similarly to before, the left-hand side of \cref{eq:stabilization pcd 1} can be expressed discretely with 
\begin{equation*}
(-M_p^{-1} B ) (M_u^{-1} F)^{-1} (M_u^{-1} B^{T}) = M_p^{-1} (-B F^{-1} B^T).
\end{equation*}
This is the scaled Schur complement matrix in the case of unstabilized mixed finite element approximations. 
If our discretization is stabilized, the Schur complement expands to $S= - B F^{-1} B^T - C$, such that we obtain
\begin{equation*}
M_p^{-1} (-B F^{-1} B^T - C)
\end{equation*}
for the left-hand side of \cref{eq:stabilization pcd 1}. 
For the right-hand side we get
\begin{equation*}
(M_p^{-1} F_p)^{-1}(-M_p^{-1} B)(M_u^{-1} B^{T}) = F_p^{-1} (-B M_u^{-1} B^{T}).
\end{equation*}
Thus, for stabilized elements, we have
\begin{equation*}
M_p^{-1} (-B F^{-1} B^T - C) \approx F_p^{-1} (-B M_u^{-1} B^{T}).
\end{equation*}
The matrix $-BM_u^{-1}B^T$ stems from the discretization of $\nabla\cdot\nabla$, but it is not a good representation of a discrete Laplacian in cases where the underlying mixed approximation is not uniformly stable \cite[Chapter~9, p.~369]{Elman_fluid_book}. 
Furthermore, for inherently unstable discretizations, highly oscillating components are in the null space of $B^T$, and, consequently, $B M_u^{-1} B^{T}$ is singular. 
However, we can fix this by using $-A_p$ as a more direct discretization of $\Delta$, as described in the previous section. 
As a result, in the stabilized case, we obtain the Schur complement approximation $S\approx -M_p F_p^{-1} A_p=S_{\text{PCD}}$, which is the same as in the stable case; cf. \cite[p.~351]{Cyr_teko_stabilized}.


\subsubsection{Boundary conditions}\label{Sec: PCD Boundary}

The assembly of the original Schur complement $S$ contains information about the boundary conditions. 
Similarly, we need to account for this information in the components of $S_{\text{PCD}}$. 
For a detailed discussion on this topic, we refer to \cite{Elman_fluid_book,Elman_BC_Navier_Stokes}, where different boundary conditions are considered. 
Originally in \cite{Kay_PCD_2002}, $F_p$ and $A_p$ were constructed with Neumann, i.e., do-nothing boundary conditions. 
Then, based on the analysis of a one- and two-dimensional problem, other boundary conditions were derived; see \cite{Elman_fluid_book, Elman_BC_Navier_Stokes}.

The derivation of the preconditioner using homogeneous Neumann boundary conditions is based on enclosed-flow problems (e.g., lid-driven cavity). 
Consequently, we should expect it to perform better in these types of problems compared to nonenclosed-flow problems. 
In in- and outflow problems, $A_p$ and $F_p$ can be treated using different boundary conditions to improve performance. 
Note that $A_p$ is not only used in the PCD Schur complement approximation $-M_p F_p^{-1} A_p$ but also for the assembly of $F_p=\tfrac{3}{2\Delta t}M_p+N_p+\nu A_p$.
The setup of the two instances is based on the same Neumann matrix corresponding to $A_p$ but, subsequently, boundary conditions are set separately for $A_p$ in $S_{\text{PCD}}$ and for $F_p$. 
We will revise a few results from \cite{Elman_fluid_book, Elman_BC_Navier_Stokes,PCD_bc_Tuminaro}. 
The main idea is to construct $F_p$ \cref{F_p} and $A_p$ as though they come from the corresponding boundary value problem. 
For the discrete Laplace operator~$A_p$, forcing the discrete pressure $p_h$ to satisfy a homogeneous Dirichlet boundary condition along the outflow boundary of $\partial \Omega_N$ is recommended; see \cite{Elman_fluid_book}. 
For the $\partial \Omega_D=\partial\Omega\setminus\partial\Omega_N$ part of the boundary of the domain, a Neumann condition $\frac{\partial p_h}{\partial \mathbf{n}} = 0$ should be applied in $A_p$. 
For the discrete convection--diffusion operator $F_p$ it is suggested that the Robin condition
\begin{equation*}
- \nu \frac{\partial p_h}{\partial \mathbf{n}} + (\mathbf{w}_h \cdot \mathbf{n}) p_h = 0
\end{equation*}
for all boundary edges is appropriate; see, for example, \cite{Elman_fluid_book,Elman_BC_Navier_Stokes}. 
Note that this includes edges where Dirichlet velocity boundary conditions are applied, such as the inflow boundary or the walls. 
The use of Dirichlet conditions for $F_p$ along the outlet $\partial \Omega_N$ can be beneficial. 
In \cite[Section 7]{Elman_BC_Navier_Stokes} different boundary conditions for $F_p$ were discussed, and the use of a Robin boundary condition along the inlet and Dirichlet boundary conditions along the outlet showed good results for Reynolds numbers between 10 and 100. For higher Reynolds number in the range 200 to 400, results for Robin on the inlet and Neuman boundary conditions on the outlet proved best. 
In \cite{cyr_teko_2016}, using a similar implementation to ours, the use of scaled Dirichlet conditions was motivated. Other scaling options are mentioned in \cite{Elman_BC_Navier_Stokes}. 

Since it is widely emphasized \cite{Elman_fluid_book, Elman_BC_Navier_Stokes, PCD_bc_Tuminaro,cyr_teko_2016} that the boundary conditions applied to PCD can have a substantial effect on the convergence (see also \Cref{table: Block PCD BC stationary}), we will consider the options in \Cref{Table: PCD BC-Table} for simulations but keep the recommended one (BC--3) as a default, where $A_p$ has a Dirichlet condition on the outlet and $F_p$ a Robin condition on the inlet.

\begin{table}[!tb]
\setlength{\tabcolsep}{0pt}
\centering
\caption{\textbf{Boundary conditions used in the PCD preconditioner.} Different strategies for setting homogeneous boundary conditions in $A_p$ and $F_p$ denoted as D: Dirichlet, N: Neumann, and R: Robin. The default strategy (BC--3) is highlighted.} \label{Table: PCD BC-Table}
\begin{tabularx}{\linewidth}{@{} p{2cm} >{\centering\arraybackslash}X>{\centering\arraybackslash}X>{\centering\arraybackslash}X p{1cm} >{\centering\arraybackslash}X>{\centering\arraybackslash}X>{\centering\arraybackslash}X @{} }
\toprule
       & \multicolumn{3}{c}{$A_p$} && \multicolumn{3}{c}{$F_p$} \\
\cmidrule{2-4}\cmidrule{6-8}
       & $\partial\Omega_{\text{out}}$ & $\partial\Omega_{\text{in}}$ & $\partial\Omega_{\text{wall}}$ &&
			   $\partial\Omega_{\text{out}}$ & $\partial\Omega_{\text{in}}$ & $\partial\Omega_{\text{wall}}$
				\\
\,(BC--1) & D & N & N &&
         D & N & N
				\\
				
\,(BC--2) & D & N & N && 
         D & R & N
				\\
\rowcolor{tblHighlGray} \,(BC--3) & D & N & N &&
         N & R & N \\
\bottomrule
\end{tabularx}
\end{table}


\subsection{LSC (Least-Squares Commutator)}\label{subsection: LSC}

The construction of the LSC preconditioner follows \cite{Elman_LSC_2006,Elman_OG_BFBT} and \cite[Section~9.2.3]{Elman_fluid_book}. 
Similarly to before, we try to define a suitable approximation of the Schur complement by reducing the error of a commutator. 
Unlike in \Cref{section: Block Preconditioner}, however, we do not use the corresponding discrete commutator directly to define a Schur complement approximation with the help of a newly assembled matrix $F_p$. 
Instead, $F_p$ is replaced by a matrix that does not require the assembly of an additional matrix (we want to avoid having to assemble the discrete convection--diffusion operator on the pressure) and that results in a small commutator with respect to some norm.

The commutator $\varepsilon=\nabla \cdot (-\nu \Delta + \mathbf{w}_h \cdot \nabla ) - (-\nu \Delta + \mathbf{w}_h \cdot \nabla)_p \nabla \cdot$ from \Cref{section: Block Preconditioner} can be regarded as the adjoint commutator to (cf. \cite[Remark~9.3]{Elman_fluid_book})
\begin{equation*}
\varepsilon^*:=(-\nu \Delta - \mathbf{w}_h \cdot\nabla ) \nabla - \nabla (-\nu \Delta - \mathbf{w}_h \cdot \nabla )_p.
\end{equation*}
This time we base the approximation on the minimization of the commutator in the least-squares sense. 
Specifically, we seek to minimize the operator norm 
\begin{equation*}
\sup_{p_h \neq 0} \frac{\left\Vert\varepsilon^* p_h \right\Vert_{L^2(\Omega)}}{\Vert p_h \Vert_{L^2(\Omega)}}.
\end{equation*}
Similar to before, we can formulate this in a discrete sense and minimize
\begin{equation}\label{eq:LSC:norm discrete commutator}
	\sup_{p_h \neq 0} \frac{
		\left\Vert\varepsilon_h^* p_h \right\Vert_{M_u}
		}{
		\Vert p_h \Vert_{M_p}},
		\quad
		\varepsilon_h^* = (M_u^{-1} F) (M_u^{-1} B^{T} ) - (M_u^{-1} B^{T})(M_p^{-1} F_p),
\end{equation}
where in this context $p_h$ corresponds to its coefficient vector, $\Vert \mathbf{v} \Vert_{M_u} = \langle M_u \mathbf{v},\mathbf{v} \rangle ^{\frac{1}{2}}$, and $\Vert p_h \Vert_{M_p} = \langle M_p p_h, p_h \rangle^{\frac{1}{2}}$. 
Instead of using $F_p$ as the discrete convection--diffusion operator on the pressure, we replace it with a matrix that yields a small (not necessarily minimal) value for \cref{eq:LSC:norm discrete commutator}. 
One way to achieve this is to minimize the individual vector norms of the columns of the discrete commutator, that is, by deﬁning the $j$th column $[F_p]_j$ of $F_p$ to solve the weighted least-squares problem
\begin{equation*}
\min_{[F_p]_j} \left\Vert 
	\left[ M_u^{-1} F M_u^{-1} B^{T}\right]_j - M_u^{-1} B^{T} M_p^{-1} \left[F_p\right]_j 
	\right\Vert_{M_u}.
\end{equation*}
The associated normal equations are
\begin{equation*}
M_p^{-1} B M_u^{-1} B^{T} M_p^{-1} \left[F_p \right]_j 
= 
\left[ M_p^{-1}B M_u^{-1} F M_u^{-1} B^{T} \right]_j,
\end{equation*}
which leads to the following definition of $F_p$:
\begin{equation}
\label{eq:LSC:Fp}
F_p = M_p (B M_u^{-1} B^{T})^{-1} (B M_u^{-1} F M_u^{-1} B^{T} ).
\end{equation}
Similary to how we have derived the Schur complement approximation \cref{eq: BFB approx}, we can pre-multiply $\varepsilon_h^*$ with $B F^{-1} M_u$ and post-multiply it with $F_p^{-1} M_p$. 
Under the assumption that the commutator is small, we obtain 
\begin{equation*}
-B F^{-1} B^T \approx -B M_u^{-1} B^T F_p^{-1} M_p.
\end{equation*}
Substituting $F_p$ from \cref{eq:LSC:Fp} into this Schur complement approximation gives 
\begin{equation*}
-B F^{-1} B^T
\approx 
	-\left( B M_u^{-1} B^T \right) 
	\left(B M_u^{-1} F M_u^{-1} B^{T} \right)^{-1} 
	\left(B M_u^{-1} B^{T}\right).
\end{equation*}
The approximate inverse of the Schur complement for the LSC block preconditioner for a stable discretization reads
\begin{equation}\label{Eq: S_LSC}
S_{\text{LSC}}^{-1} := -(B M_u^{-1} B^{T})^{-1} (B M_u^{-1}F M_u^{-1}B^{T}) (B M_u^{-1} B^{T})^{-1}.
\end{equation} 
We define the LSC block-triangular preconditioner
\begin{equation*}
\mathcal{B}_{\text{LSC}}^{-1} = \begin{bmatrix}
F & B^{T} \\
0 & S_{\text{LSC}}
\end{bmatrix}^{-1}.
\end{equation*}
If the mass matrix $M_u$ is replaced with the identity matrix, the resulting Schur complement approximation, 
\begin{equation}\label{Eq: S_BFBT}
S_{\text{BFBt}}^{-1} = -(B B^{T})^{-1} (B F B^{T}) (B B^{T})^{-1},
\end{equation}
is the classical BFBt preconditioner \cite{Elman_OG_BFBT}, opposed to the LSC preconditioner being a scaled BFBt method.

As for the construction of the PCD preconditioner, $B M_u^{-1} B^{T}$ can be replaced with the pressure-Laplace matrix:
\begin{equation}\label{Eq: S_LSC_Ap}
S_{\text{LSC}_{A_p}}^{-1}
=
 -A_p^{-1} (B M_u^{-1}F M_u^{-1}B^{T}) A_p^{-1}.
\end{equation}
We have chosen to use this replacement strategy in our simulations; cf. \Cref{sect:teko:lsc}. 
As before, $A_p$ can be assembled using different boundary conditions \cite[Sect.~9.2.2]{Elman_fluid_book}; see also \cite{Elman_BC_Navier_Stokes}. 
In this work, we will only use (BC--2) for LSC variants, that is, a Dirichlet condition at the outflow of the domain.


\paragraph{Stabilization}

In non inf-sup stable discretizations, the LSC approach requires some modification. 
As mentioned in \Cref{sssection:PCD:stabilized}, the approximation $B M_u^{-1} B^{T}$ of the discrete Laplacian is singular due to the nontrivial null space of~ $B^{T}$. 
In \cite[Sect.~4.2]{elman:2008:lsp} (see also \cite[Sect. 3.3.3]{Cyr_teko_stabilized} and \cite[Sect. 3.2.2]{ELMAN2008_compParallelBlock}) a modification was proposed by adding scaled stabilization operators:
\begin{equation*}
S_{\text{LSC}_{\text{stab}}}^{-1}
= 
-(B M_u^{-1} B^{T} + \gamma C)^{-1} (B M_u^{-1}F M_u^{-1}B^{T}) (B M_u^{-1} B^{T} + \gamma C)^{-1} - \alpha D^{-1},
\end{equation*}
where, using $\operatorname{diag}(\cdot)$ for the diagonal part and $\rho_{\text{eig}}(\cdot)$ for the spectral radius of a matrix,
\begin{align*}
\gamma &:= \frac{\rho_{\text{eig}}\left(M_u^{-1} F\right)}{3}, && \\
D &:= \operatorname{diag}\left( \hat{S}_{\text{diag}}\right), 
&\hat{S}_{\text{diag}} &:= -B \operatorname{diag}(F)^{-1} B^{T} - C, \\
\alpha &:= \frac{-1}{\rho_{\text{eig}}\left(\hat{S}_{\text{diag},C=0} D^{-1} \right) },
&\hat{S}_{\text{diag},C=0} &:= -B \operatorname{diag}(F)^{-1}B^{T}.
\end{align*}
As in \cite[Sect. 3.3.3]{Cyr_teko_stabilized}, we replace $BM_u^{-1} B^{T} + \gamma C$ with $A_p$ with no special consideration to the stabilization matrix~$C$:
\begin{equation*}
S_{\text{LSC}_{\text{stab},A_p}}^{-1}
= 
-A_p^{-1} (B M_u^{-1}F M_u^{-1}B^{T}) A_p^{-1} - \alpha D^{-1}.
\end{equation*}


\paragraph{Inverse of velocity mass matrix}

Similar to a strategy applied in the construction of the SIMPLE preconditioner (see \Cref{subsection: SIMPLE}), we will approximate the inverse of $M_u$ with one of the following two matrices: 
Let $\delta^{ij}$ denote the Kronecker delta; then 
\begin{equation*}
H_{M_u}^D := \operatorname{diag}(M_u)^{-1}, \qquad 
H_{M_u}^{\Sigma} := \delta^{ij}\Big(\sum_{k=1}^{N_u} |(M_u)_{i,k}|\Big )^{-1}.
\end{equation*}


\subsection{SIMPLE}\label{subsection: SIMPLE}

The SIMPLE (Semi-Implicit Method for Pressure-Linked Equations) preconditioner was originally introduced by Patankar and Spalding in \cite{Patankar_SIMPLE_1972} to solve the Navier--Stokes equations, and the SIMPLE block preconditioner is based on it; see also \cite{Quarteroni_numApprox_NavSto,ELMAN2008_compParallelBlock}.

First, the lower triangular and the block-diagonal component of the LDU decomposition \cref{eq:Navier Stokes LDU} are grouped together in an $LU$ decomposition:
\begin{equation*}
\begin{bmatrix}
F & 0 \\
B & S
\end{bmatrix}
\begin{bmatrix}
I & F^{-1} B^{T} \\
0 & I 
\end{bmatrix}.
\end{equation*} 
Then, the inverse of $F$ in the upper right block is replaced by an approximation $H_F$, which differs whether SIMPLE or SIMPLEC is used; 
a SIMPLE variant that is not covered in this work is SIMPLER \cite{simple-like_2016,SIMPLER_OG_Doorman}. 
The Schur complement $S$ is approximated with $S_{\text{SIMPLE}}=-C-B H_F B^{T}$. 
The introduction of an under-relaxation parameter $\alpha$ completes the construction of the preconditioner $\mathcal{B}_{\text{SIMPLE}}^{-1}$:
\begin{equation*}
\mathcal{B}_{\text{SIMPLE}} = 
\begin{bmatrix}
F & 0 \\
B & S_{\text{SIMPLE}} 
\end{bmatrix}
\begin{bmatrix}
I & \frac{1}{\alpha} H_F B^{T} \\
0 & \frac{1}{\alpha} I 
\end{bmatrix}.
\end{equation*} 

$H_F$ is a diagonal matrix and depends on the specific SIMPLE variant. 
We denote the matrix of the default variant SIMPLE as $H_F^D$ and the one of SIMPLEC as $H_F^{\Sigma}$:
\begin{align*}
H_F^D &= \operatorname{diag}(F)^{-1}, \tag{SIMPLE}\\
H_F^{\Sigma} &= \delta^{ij}\Big(\sum_{k=1}^{N_u} |F_{i,k}|\Big )^{-1}. \tag{SIMPLEC}
\end{align*}
If we expand $F$ into its parts 
\begin{equation*}
F = \frac{3}{2\Delta t}  M_u + \nu  A +  N + W,
\end{equation*}
with velocity mass matrix $M_u$, discrete Laplacian $-A$, and advective terms $N+W$, the well-conditioned mass matrix dominates $F$ for small time steps $\Delta t$. 
Since $M_u$ is spectrally equivalent to the lumped mass matrix $\operatorname{diag}(M_u)$, approximating $F$ via its diagonal is sufficient for small $\Delta t$. 
The performance of SIMPLE and SIMPLEC can deteriorate for stationary flow problems if $\Delta t$ is large or if advective forces dominate. 


\section{Additive overlapping Schwarz preconditioners} \label{sec:schwarz}

As discussed in~\Cref{sec:disc,Sec: Inexact Newton + Forcing Term}, we solve the Navier--Stokes equations using a Newton--Krylov approach, that is, by solving the linearized systems arising from Newton's method using a Krylov method. 
Since the tangent matrix~\cref{Eq: Saddle point problem Navier Stokes} is not positive definite and generally unsymmetric, we employ the preconditioned GMRES method~\cite{saad:1986:gmres} as the Krylov method. 
In~\Cref{section: Block Preconditioner}, we have discussed several options to block-precondition the system matrix~\cref{Eq: Saddle point problem Navier Stokes}, where each application of a block-preconditioner requires the separate solution of velocity and pressure systems. 
For large-scale problems, it is often unfeasible to use direct solvers for the solution of these subsystems. 
In this paper, we will approximate the involved inverses with overlapping Schwarz domain decomposition methods. 
Specifically, we use block Schwarz preconditioners and monolithic Schwarz preconditioners, the latter of which are constructed by directly employing Schwarz preconditioners on the system matrix; see \Cref{sec: BlockSchwarz,Section: OSP: Monolithic} for details on the block and monolithic approaches, respectively. 
For the implementation of both techniques, we employ the \texttt{Trilinos} package \texttt{FROSch}; cf.~\Cref{sec:frosch}. 
We introduce Schwarz preconditioners under simplified assumptions, that is, for a linear system
\begin{equation*}
	K x = b,
\end{equation*}
resulting from the discretization of a Laplacian model problem on the computational domain $\Omega$ using, for example, piecewise linear finite elements with the finite element space $V^h = V^h(\Omega)$. 
In this case, the matrix~$K$ is symmetric positive definite, and theoretical results for the convergence are available. We refer to~\cite{toselli:2005:ddm} for a more detailed introduction. 
The following construction is essentially the same for general problems but is then not based on a robust theoretical background.

Let $\Omega$ be decomposed into nonoverlapping subdomains $\{ \Omega_i \}_{i=1}^{N}$ that are conforming with the finite element discretization; that is, the closure of every subdomain is the union of its finite elements. 
We extend the subdomains by $k$ layers of finite elements, resulting in an overlapping domain decomposition $\{ \Omega_i' \}_{i=1}^{N}$ with overlap $\delta = kh$. 
Note that extending the subdomains by layers of elements is a geometric operation based on the mesh. 
To only use information available via the matrix~$K$, the subdomains can be algebraically extended by iteratively adding those degrees of freedom that are in the neighborhood of the (overlapping) subdomain, defining neighborhood via the sparsity pattern of~$K$; 
cf. \cite{Heinlein:2020:FRO,heinlein_fully_2021}.  
For example, if $\Omega_i$ is extended by one layer, for all degrees of freedom $l$ in $\Omega_i$, degrees of freedom $j$ for which $K_{l,j}\neq 0$ are added to $\Omega_i$ to obtain $\Omega_i'$.
This algebraic overlap construction is used in our implementation. 
We denote the size of the algebraically determined overlap by $\hat{\delta} = k$ for $k$~iterations of the overlap-constructing algorithm. 
Note that in the context of conforming Lagrange finite elements, $\hat{\delta}=k$ amounts to a geometric overlap of $\delta=(k+1)h$.

Based on the overlapping domain decomposition, we define a restriction operator $R_i: V^h \rightarrow V_i^h$, $i=1,\dots,N$, to map from the global finite element space $V^h$ to the local finite element space $V_i^h := V^h(\Omega_i')$ on the overlapping subdomains $\Omega_i'$. 
Moreover, we define the corresponding prolongation operator $R_i^T: V_i^h \rightarrow V^h$, which is the transposed of the restriction operator and extends a local function on $\Omega_i'$ by zero outside of $\Omega_i'$. 
We define a second prolongation $\widetilde{R}_i^T$: 
Standard additive Schwarz is given by $\widetilde{R}_i^T=R_i^T$. 
Restricted additive Schwarz \cite{SAS_RAS_Cai} uses $R_i^T$ but sets values to zero such that $\sum_{i=1}^N\widetilde{R}_i^T R_i$ is the identity operator. 
In this work, we define $\widetilde{R}_i^T$ to satisfy the identity-operator property via an inverse multiplicity: 
A row of $R_i^T$ corresponding to the finite element node $x^h$ is scaled with the inverse of $|\{j \in \{1,\dots,N\}: x^h\in\Omega_j'\setminus\partial\Omega_j'\}|$, where $| \cdot |$ is the cardinality of the set; cf. \cite[Remark~2.7]{SAS_RAS_Cai} and \cite[Sect.~4.1]{heinlein_reduced_2019}.

Using the aforementioned operators $R_i$, $R_i^{T}$, and $\widetilde{R}_i^T$, and defining local overlapping stiffness matrices $K_i := R_i K R_i^{T}$, $i=1,\ldots,N$, the (scaled) additive one-level overlapping Schwarz preconditioner is given by
\begin{equation} \label{eq:one_level schwarz}
	M_{\text{OS-1}}^{-1}
	=
    \sum_{i=1}^{N}  \widetilde{R}_i^T K_i^{-1} R_i.
\end{equation}

Even for a scalar elliptic Laplacian model problem, the one-level Schwarz preconditioner~\cref{eq:one_level schwarz} is generally not numerically scalable: The number of Krylov iterations will increase with an increasing number of subdomains. 
In order to define a numerically scalable overlapping Schwarz preconditioner, which yields a convergence rate independent of the number of subdomains, a coarse level can be introduced. 
With a coarse interpolation operator $\Phi:V_0\to V^h$, which is still to be specified, the additive two-level overlapping Schwarz preconditioner reads
\begin{equation} \label{eq:two_level schwarz}
	M_{\text{OS-2}}^{-1}
	=
	\Phi K_0^{-1} \Phi^{T}
	+
	\underbrace{\sum_{i=1}^{N} \widetilde{R}_i^T K_i^{-1} R_i}_{= M_{\text{OS-1}}^{-1}},
\end{equation}
where $K_0 = \Phi^{T} K \Phi$. 
The columns of $\Phi$ form the basis of the coarse space, $V_0$, and $K_0$ is a Galerkin projection of~$K$ into the coarse space. 
A natural choice for coarse basis functions are Lagrangian basis functions defined on a coarse triangulation, which yields a scalable preconditioner. 
However, the definition of a coarse triangulation heavily depends on geometric information and cannot be easily done in a conforming way for complex domain geometries. 
In the next section, we will discuss alternatives.


\subsection{GDSW-type coarse spaces}\label{subsec: Coarse Spaces}

For this work, we choose a very flexible approach for the construction of a coarse space that was introduced with the Generalized Dryja--Smith--Widlund (GDSW) preconditioner in~\cite{dohrmann_domain_2008,dohrmann_family_2008}. 
It is inspired by nonoverlapping domain decomposition methods, such as those from FETI--DP~\cite{farhat_scalable_2000,farhat_feti-dp_2001} and BDDC~\cite{cros_preconditioner_2003,dohrmann_preconditioner_2003}. 
The approach allows to define many different variants, three of which (GDSW, RGDSW, GDSW$^\star$) we will introduce and use for simulations. 
They require little to no additional information to the system matrix. 
In this approach, the nonoverlapping subdomains $\{ \Omega_i \}_{i=1}^{N}$ are chosen as the elements of the coarse triangulation, and the coarse basis functions are defined via extension of trace functions defined on the interface
\begin{equation*}
	\Gamma 
	= 
	\big \{ 
	x \in \big( \Omega_i \cap \Omega_j \big) \setminus \partial \Omega_D 
	: 
	i \neq j, 1 \leq i, j \leq N 
	\big \},
\end{equation*}
where $\partial \Omega_D$ is the global Dirichlet boundary; note the comment in \Cref{subsubsec: disc} for discontinuous finite elements, in which case the interface does not contain any degrees of freedom. 
More precisely, first, a partition of unity is defined on $\Gamma$. 
The different variants of the coarse spaces differ only in the choice of the partition of unity on the interface. 
Then the functions are extended into the interior of the subdomains. 
The extension is defined algebraically using submatrices of~$K$ and corresponds to a discrete harmonic extension in the case of a Laplace problem. 
In this specific case, the extension of a partition of unity on the interface gives a partition of unity on the whole domain $\Omega$. 
Given interface values $\Phi_\Gamma$, the full matrix $\Phi$ is obtained by solving for all degrees of freedom $I$ that do not correspond to the interface~$\Gamma$:
\begin{equation}\label{eq: energy minimizin extension}
	\Phi = \begin{bmatrix}
		\Phi_I \\
		\Phi_\Gamma 
	\end{bmatrix}
	= 
	\begin{bmatrix}
		-K_{II}^{-1} K_{I\Gamma} \Phi_\Gamma \\
		\Phi_\Gamma
	\end{bmatrix}.
\end{equation}
The blocks $K_{II}$ and $K_{I\Gamma}$ are obtained by reordering and partitioning $K$ as 
\begin{equation*}
	K =
	\begin{bmatrix}
		K_{II}        & K_{I \Gamma} \\
		K_{\Gamma I } & K_{\Gamma \Gamma} 
	\end{bmatrix}.
\end{equation*}
Note that $K_{II}$ is block-diagonal, and hence the solver required for $K_{II}^{-1}$ can be dealt with independently for each subdomain. 
Moreover, as only submatrices of $K$ are required, the construction of the coarse space is algebraic. 
For other problems, like elasticity, additional geometric information may be required; cf., e.g., \cite[sect.~4]{dohrmann:2009:osa}. 
In this work, for the fluid problems, we do not make use of geometric information and obtain an algebraically constructed coarse space. 
The description above outlines the construction for a scalar problem; the extension to vector-valued problems is given in \Cref{subsubsec: vector}. 
The coarse spaces used in this work are all based on GDSW-type preconditioners and on extending the interface values into the interior via~\cref{eq: energy minimizin extension}; they only differ in the definition of the interface partition of unity $\Phi_\Gamma$. 
We review the different approaches implemented in \texttt{FROSch} in~\Cref{Sec: Results}. 
We consider three coarse space variants: GDSW~\cite{dohrmann_domain_2008,dohrmann_family_2008} and reduced-dimension GDSW (RGDSW)~\cite{dohrmann_design_2017} (``Option~1'') -- due to its significantly reduced dimension, the RGDSW coarse space improves parallel scalability; cf.~\cite{heinlein_improving_2018} -- and we will introduce an intermediate approach that we denote as GDSW$^\star$. 
It yields smaller coarse spaces than the classical GDSW approach but larger coarse spaces than the RGDSW approach; in two dimensions, GDSW$^\star$ and RGDSW are identical. 
Preliminary results for GDSW$^\star$ based on its implementation in \texttt{FROSch} can be found in~\cite[p.~160]{diss_hochmuth_2020}, there under the working title \textit{GDSW-Star}. 
A related approach was taken in the context of adaptive coarse spaces in \cite{knepper:2022:dissertation}, there under the name R--WB--GDSW.


\subsubsection[Interface partition of unity of GDSW, GDSW*, and RGDSW coarse spaces]{Interface partition of unity of GDSW, GDSW$^\star$, and RGDSW coarse spaces}\label{subsubsec: IPOU}

\begin{figure}[!tb]
	\hspace*{\fill}
	\includegraphics[width=0.23\textwidth]{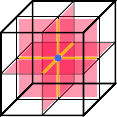}
	\hspace*{\fill}
	\includegraphics[width=0.23\textwidth]{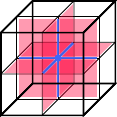}
	\hspace*{\fill}
	\includegraphics[width=0.23\textwidth]{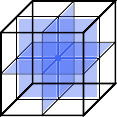}
	\hspace*{\fill}%
	\caption{%
		\textbf{Visualization of three types of (overlapping) interface decompositions in a $2\times 2\times 2$ domain decomposition.}
		\textbf{Left:} GDSW interface components: 6~edges, 1~vertex, and 12~faces. 
		\textbf{Center:} GDSW$^\star$ interface components: 1~vertex-based component (union of vertex and 6~edges), 12~faces. 
		\textbf{Right:} RGDSW interface component: 1~vertex-based component (union of vertex, 6~edges, 12~faces).}
	\label{figure: IPOU}
\end{figure}

The classical GDSW coarse space defines the interface partition of unity based on a nonoverlapping decomposition of the interface degrees of freedom into interface components of faces, edges, and vertices; cf. \Cref{figure: IPOU}. 
Contrary, the RGDSW coarse spaces from \cite{dohrmann_design_2017} employ overlapping interface components, combining each vertex with all its adjacent faces and edges; 
the resulting decomposition overlaps in the faces and edges. 
Combining vertices with edges and faces reduces the number of interface components from the number of faces, edges, and vertices to only the number of vertices. 
In~\cite{heinlein_adaptive_2022}, a variant of the RGDSW coarse space is introduced in which the RGDSW interface decomposition is made nonoverlapping without increasing the number of interface components. 
The new GDSW$^\star$ coarse space employs a partition of unity of intermediate size, where each vertex is only combined with the adjacent edges, and the faces are kept separate. 

The interface partition of unity is obtained by defining one function for each interface component. 
As a result, the number of coarse basis functions, which is equal to the coarse space dimension, is reduced for GDSW$^\star$ and even further reduced for RGDSW; cf.~\Cref{Table: coarse space dim}. 
See also \cite[Fig.~5.4]{knepper:2022:dissertation} for a comparison of the number of interface components for four different types of meshes and domain decompositions; note that, although the precise definition of the interface components differs in \cite{knepper:2022:dissertation}, their quantity is the same. 
There, the analogon of GDSW\expStar{} is denoted as R--WB--GDSW.

\begin{table}[!tb]
\centering
\caption{\textbf{Coarse space dimensions for scalar Laplace problem and the different GDSW-type coarse spaces used in this work.}
Vertices, edges, and faces refer to the respective interface components of the domain decomposition.} \label{Table: coarse space dim}
\begin{tabularx}{\linewidth}{@{} >{\raggedright\arraybackslash}X>{\raggedright\arraybackslash}X @{}}
\toprule
Coarse Space & Coarse Space Dimension \\
\midrule
GDSW         & number of vertices + faces + edges \\
GDSW$^\star$ & number of vertices + faces \\
RGDSW        & number of vertices \\
\bottomrule
\end{tabularx}
\end{table}

To define the interface partition of unity of the GDSW, GDSW$^\star$, and RGDSW coarse spaces, we decompose the interface $\Gamma$ into components $\gamma_k \subset \Gamma$ such that
\begin{equation*}
	\Gamma = \bigcup_k \gamma_k. 
\end{equation*}
This union can be disjoint or overlapping. 
Then, we define functions $\Phi_{\Gamma,k}$ on $\Gamma$, with $\operatorname{supp}(\Phi_{\Gamma,k}) \subset \gamma_k$, such that
\begin{equation*}
	\sum_k \Phi_{\Gamma,k} \equiv 1 \text{ on } \Gamma.
\end{equation*}
Then, we gather the functions $\Phi_{\Gamma,k}$ (that is, their respective coefficient vectors) as columns of a matrix 
\begin{equation*}
	\Phi_\Gamma = \begin{bmatrix}
		\Phi_{\Gamma,1} & \Phi_{\Gamma,2} & \hdots &
	\end{bmatrix},
\end{equation*}
which is subsequently extended using~\cref{eq: energy minimizin extension}.

It remains to define the partition-of-unity functions $\Phi_{\Gamma,k}$ for GDSW, GDSW$^\star$, and RGDSW. 
We introduce the following formal definitions of faces, edges, and vertices. 
In three dimensions:
\begin{itemize}[leftmargin=2.0em,topsep=4pt,itemsep=2pt]
	\item A face is a set of nodes that belongs to the same two subdomains. Let $\mathcal{F}$ be the set of all faces.
	\item An edge is a set of at least two nodes that belongs to the same (more than two) subdomains. Let $\mathcal{E}$ be the set of all edges.
	\item A vertex is a single node that belongs to the same (more than two) subdomains. Let $\mathcal{V}$ be the set of all vertices.
\end{itemize}

These interface components are exemplarily shown in \Cref{figure: IPOU} (left); 
all of them are used to construct GDSW coarse spaces. 
In two dimensions, the definition of an edge is the same as that of a face in 3D, and the definition of the vertex remains the same. 
In the following, we will introduce the different partitions of unity considered in this work, focusing on the three-dimensional case.


\paragraph{GDSW}

For the construction of the coarse space, the $\gamma_k$ are given by faces, edges, and vertices. 
We define the discrete partition-of-unity functions via 
	\begin{equation*}
		\chi_f, f \in \mathcal{F}, 
		\quad
		\chi_e, e \in \mathcal{E},
		\quad
		\chi_v, v \in \mathcal{V};
	\end{equation*}
$\chi_s$ is the discrete characteristic function that is one in the finite elements nodes contained in~$s$ and zero elsewhere. 
Then, we have that 
\begin{equation*}
	\sum_{f \in \mathcal{F}} \chi_f
	+
	\sum_{e \in \mathcal{E}} \chi_e
	+
	\sum_{v \in \mathcal{V}} \chi_v
	\equiv
	1
	\text{ on } \Gamma,
\end{equation*}
since the vertices, edges, and faces form a nonoverlapping decomposition of $\Gamma$.


\paragraph{RGDSW}

Let $\mathcal{V} = \lbrace v_1,\ldots,v_{N_\mathcal{V}} \rbrace$. 
An interface component $\gamma_k$ is the union of $v_k \in \mathcal{V}$ and its adjacent faces and edges. 
Note that, as discussed in~\cite{dohrmann_design_2017}, the definition is slightly more involved for general domain decompositions, where not every edge and face may have an adjacent vertex; for simplicity, we do not discuss these cases. 

In \cite{dohrmann_design_2017,heinlein_adaptive_2022} various options for defining the corresponding partition-of-unity functions are given. 
Here, we employ and describe the algebraic variant ``Option~1'' from~\cite{dohrmann_design_2017}, which is based on a simple multiplicity scaling (cf.~\cite[eq.~(1)]{dohrmann_design_2017}). 
In particular, we first define a component multiplicity $\mathcal{M}(n)$ of an interface node~$n$ as the number of (RGDSW) interface components $\gamma_k$ that contain~$n$: 
\begin{equation*}
	\mathcal{M}(n) := |\{ k : n\in\gamma_k \}|.
\end{equation*}
For example, for an edge with two adjacent vertices, we have $\mathcal{M}(n)=2$ for all nodes $n$ on the edge. 
We can then define a function $\Phi_{\Gamma,k}$ associated with $\gamma_k$ by using the inverse of the multiplicity:
\begin{equation} \label{eq:mult_scaling}
	\Phi_{\Gamma,k}(n) := 
	\begin{cases}
		\nicefrac{1}{\mathcal{M}(n)}, & n \in \gamma_k,\\
		0, & \text{else.}
	\end{cases}
\end{equation}
The partition-of-unity property is satisfied as in every interface node $n \in \Gamma$, we have 
\begin{equation*}
	\sum_k \Phi_{\Gamma,k} (n) = 1.
\end{equation*}


\paragraph{GDSW$^\star$}

Finally, we introduce the partition of unity of the new GDSW$^\star$ coarse space. 
The definition is based directly on those of GDSW and RGDSW. 
In particular, we define two types of interface components:
\begin{enumerate}[leftmargin=2.5em,topsep=4pt,itemsep=2pt]
\item A vertex $v_k \in \mathcal{V}$ and its adjacent edges. 
    These components overlap.
\item A face $f_k \in \mathcal{F}$. 
    Faces do not overlap with the vertex-based interface components or other faces.
\end{enumerate}
We define the partition of unity using the inverse multiplicity scaling~\cref{eq:mult_scaling}, this time defined via the set of GDSW$^\star$ interface components. 
Note that, because the faces do not overlap with other interface components, the resulting functions coincide with the discrete characteristic functions $\chi_f$ of the faces $v \in \mathcal{F}$.

As faces do not exist in two-dimensions, the partition of unity of GDSW$^\star$ and RGDSW are the same in this case.


\subsubsection{Extension to vector-valued problems}\label{subsubsec: vector}

The approaches discussed before are directly applicable to scalar problems, such as problems related to the pressure block. However, the application to vector-valued problems, such as problems related to the velocity block, requires a small modification. 

In case of, for example, a multi-dimensional Laplace operator $\Delta \mathbf{u}$ or a linearized elasticity problem, the partition-of-unity functions $\gamma_k$ are multiplied with a basis of the null space of the corresponding partial differential equation with homogeneous Neumann conditions on the boundary. 
This ensures that the null space can be represented by the coarse space, which is required for numerical scalability in case of symmetric positive semidefinite problems; cf.~\cite[p.~132]{smith:1996:dpm}. 
In our case, we use translations for the velocity block, which is motivated by $\Delta \mathbf{u}$ appearing in the diagonal velocity block of the discrete Stokes matrix. 
For the pressure, we use constant functions; note that these are in the null space of $B^T$ for an enclosed-flow problem.

Specifically, for the two-dimensional case of a velocity coarse function, the partition-of-unity functions are multiplied with 
\begin{equation*}
	r_1 := \begin{bmatrix}
		1\\
		0
	\end{bmatrix},~
	r_2:= \begin{bmatrix}
		0 \\
		1
	\end{bmatrix},
\end{equation*}
and in three dimensions with
\begin{equation*}
	r_1 := \begin{bmatrix}
		1 \\
		0 \\
		0
	\end{bmatrix},~
	r_2 := \begin{bmatrix}
		0 \\
		1 \\
		0
	\end{bmatrix},~
	r_3 := \begin{bmatrix}
		0 \\
		0 \\
		1
	\end{bmatrix},
\end{equation*}
for each interface node. 
As a result, the number of coarse space functions associated with the velocity or pressure is the space dimension multiplied with the number of partition-of-unity functions. 
This discussion is compatible with the scalar elliptic case. In particular, for the construction of pressure functions, the node-wise multiplication of the partition-of-unity functions with the scalar constant function, i.e., 
\begin{equation*}
    r := \begin{bmatrix}
        1
    \end{bmatrix},
\end{equation*}
simply yields the partition-of-unity functions.


\subsubsection{Discontinuous finite elements}\label{subsubsec: disc}

For discontinuous finite element spaces, as employed in discontinuous, piecewise linear (P1-disc.) pressure spaces, there are no interface nodes that are shared between adjacent subdomains. 
Therefore, the discrete interface is empty in this case, and the construction of GDSW, GDSW$^\star$, and RGDSW cannot be applied directly. 
As motivated in \Cref{subsubsec: vector}, we want to use constant functions for the pressure.
Therefore, on each subdomain $\Omega_i$, we employ the single discrete characteristic function $\chi_{\Omega_i}$ (i.e., a function that is constant on $\Omega_i$ and zero elsewhere) as the coarse basis function, yielding a total of $N$ coarse basis functions. 
If this occurs, the coarse space will be abbreviated with ``Vol.'', since the characteristic function is prescribed on the entire volume and not just on the interface.


\subsection{Block Schwarz preconditioner} \label{sec: BlockSchwarz}

For the presented block preconditioners in \Cref{section: Block Preconditioner}, we can use the two-level overlapping Schwarz domain decomposition method from the previous section for the approximation of the blocks $F^{-1}$ and $S^{-1}$. 
The inverse of the Schur complement is not approximated directly; 
for example, in the case of $S_{\text{PCD}}^{-1}=- A_p^{-1} F_p M_p^{-1}$, we may approximate $A_p^{-1}$ and $M_p^{-1}$ with a Schwarz method. 
Nevertheless, to avoid new notation, we assume in the following that $S^{-1}$ shall be approximated.
We define the overlapping subdomain problems separately for the velocity and pressure. 
The velocity and pressure space $\mathbf{V}^h$ and $Q^h$, respectively, are decomposed into local spaces 
\begin{equation*}
\mathbf{V}_i^h = \mathbf{V}^h \cap (H^{1}_0(\Omega_i'))^{3} \quad\text{and}\quad Q_i^h = Q^h \cap H_0^{1}(\Omega_i'),
\end{equation*}
for $i=1,\dots,N$, on the overlapping subdomains $\Omega_i'$. 
Note that, here, we restrict ourselves to conforming finite element spaces; 
the precise definition for Q2--P1-disc. is different for the discontinuous pressure space, since $H^1$ cannot be used and since functions of the local finite element spaces with a zero boundary are not only zero at the boundary nodes of $\Omega_i'$ but also on the element layer next to the boundary of $\Omega_i'$. 
The operators restricting to the overlapping subdomains belonging to the velocity and pressure degrees of freedom are defined as
\begin{equation*}
	R_{u,i} : \mathbf{V}^h \longrightarrow \mathbf{V}_i^h 
	\quad
	\text{and} 
	\quad
	R_{p,i}: Q^h \longrightarrow Q_i^h,
\end{equation*}
respectively. 
Based thereon, the local overlapping problems are defined as
\begin{equation*}
	F_i = R_{u,i} F R^{T}_{u,i} 
	\quad \text{and} \quad 
	S_i = R_{p,i} S R_{p,i}^T,
\end{equation*}
where $S$, or an approximation of $S$, depend on the underlying strategy. 
The coarse matrix is constructed similarly to before via a Galerkin product as 
\begin{equation*}
	F_0 = \Phi_u^{T} F \Phi_u
	\quad
	\text{and} 
	\quad
	S_0= \Phi_p^{T} S \Phi_p. 
\end{equation*}
Let $\widetilde{R}_{u,i}^T$ and $\widetilde{R}_{p,i}^T$ be the scaled prolongation operators (cf. \Cref{sec:schwarz}).
To approximate the inverses of $F$ and $S$, we use the additive two-level overlapping Schwarz preconditioners
\begin{equation*}
	\hat{F}^{-1} = \Phi_u F_0^{-1} \Phi_u^{T} + \sum_{i=1}^{N} \widetilde{R}_{u,i}^T F_i^{-1} R_{u,i} 
	\quad \text{ and } \quad
	\hat{S}^{-1} = \Phi_p S^{-1}_0 \Phi^{T}_p + \sum_{i=1}^{N} \widetilde{R}_{p,i}^T S^{-1}_i R_{p,i}.
\end{equation*}


\subsection{Monolithic Schwarz preconditioner}\label{Section: OSP: Monolithic}

For the saddle-point system with system matrix $\mathcal{F}$ defined in \cref{Eq: Saddle point problem Navier Stokes}, we compare the block Schwarz preconditioners with monolithic Schwarz preconditioners. This technique has been introduced in~\cite{Klawonn_1998_MixedFE, Klawonn_2000_Mono_Block} with Lagrangian coarse spaces and extended to GDSW-type coarse spaces in~\cite{Heinlein_2019_MonoFluidFlow, heinlein_reduced_2019}. In monolithic Schwarz preconditioners, we construct a single Schwarz domain decomposition preconditioner for the full system $\mathcal{F}$ at once. The monolithic additive two-level Schwarz preconditioner reads
\begin{equation}\label{eq: 2LvlOS - Saddelpoint}
	\mathcal{M}_{\text{OS-2}}^{-1} 
	= 
	  \phi \mathcal{F}_0^{-1} \phi^{T} + \sum_{i=1}^{N} \TildeRiT \mathcal{F}_i^{-1} \Ri.
\end{equation}
Different from~\cref{eq:two_level schwarz}, we denote all matrices with calligraphic letters to indicate that they are block matrices with the same block structure as $\mathcal{F}$. The local subdomain matrices for \cref{eq: 2LvlOS - Saddelpoint} are given by $\mathcal{F}_i = \Ri \mathcal{F} \Ri^T$ with the restriction operators
\begin{equation*}
\Ri = \begin{bmatrix}
R_{u,i} & 0 \\
0 & R_{p,i} 
\end{bmatrix},
\end{equation*}
where $R_{u,i}$ and $R_{p,i}$ correspond to the restriction operators for velocity $\mathbf{u}$ and pressure $p$ degrees of freedom for $\Omega_i'$; 
$\TildeRiT$ is the scaled prolongation operator; cf. \Cref{sec:schwarz}. 
Then, also the coarse matrix $\mathcal{F}_0 = \phi^{T} \mathcal{F} \phi$ has a block structure given by the matrix
\begin{equation} \label{eq:phi_mono}
\phi = 
\begin{bmatrix}
\phi _{\mathbf{u},\mathbf{u_0}} & \phi_{\mathbf{u},p_0} \\
\phi_{p,\mathbf{u_0}} & \phi_{p,p_0}
\end{bmatrix},
\end{equation}
which has the coarse basis functions as its columns. 
$\mathbf{u_0}$ and $p_0$ indicate the velocity and pressure functions in the coarse space. 
The interface values of the coarse functions are constructed as for the block matrices; 
the remaining interface values are always set to zero. 
Hence, we obtain the structure
\begin{equation*}
	\phi_\Gamma 
	=
	\begin{bmatrix}
		\phi _{\Gamma,\mathbf{u},\mathbf{u_0}} & 0 \\
		0 & \phi _{\Gamma,p,p_0}
	\end{bmatrix}.
\end{equation*}
For a block Schwarz preconditioner, the same $\phi _{\Gamma,\mathbf{u},\mathbf{u_0}}$ and $\phi _{\Gamma,p,p_0}$ are used. 
Only the subsequent extension operator from the interface to the subdomains differs. 
The extension from the interface to the subdomains is computed monolithically, analogously to \cref{eq: energy minimizin extension} but with submatrices of $\mathcal{F}$. As a result, the coarse basis $\phi$ also contains components reflecting the coupling between $\mathbf{u}$ and $p$; cf.~\cref{eq:phi_mono}. 

Since the local saddle-point subdomain matrices $\mathcal{F}_i$ are extracted from $\mathcal{F}$, they posses homogeneous Dirichlet boundary conditions for both velocity and pressure in the case of conforming finite elements. 
Additionally, we enforce a zero mean pressure value on the local overlapping subdomains as in \textit{Version~2} in \cite[Sect.~4.1]{Klawonn_1998_MixedFE}, which is beneficial in some cases, even though the pressure is uniquely determined by the problem setting; see also \cite[Sect.~5.1]{Heinlein_2019_MonoFluidFlow}. 
To this end, we introduce a pressure projection in the first level of the two-level overlapping Schwarz preconditioner 
\begin{equation}\label{eq: Pressure Projection}
\mathcal{M}_{\text{OS-2}}^{-1}
= \phi \mathcal{F}_0^{-1} \phi^{T} 
+ \sum_{i=1}^{N} \TildeRiT  \overline{\mathcal{P}}_i \mathcal{F}_i^{-1} \Ri
\end{equation}
with local projection operators $\overline{\mathcal{P}}_i$ corresponding to the overlapping subdomain $\Omega_i$. 
Here, the projection operators are of the form
\begin{equation*}
\overline{\mathcal{P}}_i
= \begin{bmatrix}
    I_{u,i} & 0 \\
    0 & \overline{P}_{p,i}
    \end{bmatrix}, \quad \overline{P}_{p,i} = I_{p,i}  - a_i (a_i^{T} a_i )^{-1} a_i^{T}.
\end{equation*}
The vector $a_i$ is defined as $ a_i = \Ri a$ via integrals of the pressure basis functions $\psi_j$:
\begin{equation*}
    a = \begin{bmatrix}
    0 \\
    a_p
    \end{bmatrix},\quad
    a_p = 
    \begin{bmatrix}
        \int_\Omega \psi_{1}  \dx 
        & \hdots & 
        \int_\Omega \psi_{N_p}  \dx 
    \end{bmatrix}^T.
\end{equation*}
Results in \cite{Klawonn_1998_MixedFE,Heinlein_2019_MonoFluidFlow} show that this approach, similar to introducing a Lagrange multiplier, can improve the convergence of the iterative linear solver.


\section{Implementation}\label{Sec: Implementation}

The software framework of this paper is based on the \texttt{Trilinos} library, in particular, we employ its domain decomposition package \texttt{FROSch} (Fast and Robust Overlapping Schwarz) \cite{Heinlein:2020:FRO} and its block-preconditioning package \texttt{Teko}~\cite{cyr_teko_2016,Cyr_teko_stabilized} to construct the preconditioners under investigation. 
Moreover, for the implementation of the model problems and the corresponding finite element discretizations, we use the \texttt{FEDDLib} (Finite Element and Domain Decomposition Library)~\cite{feddlib}, which strongly relies on several \texttt{Trilinos} packages and assembles input directly in a format to make use of the full capabilities of \texttt{FROSch} and \texttt{Teko}.

We have carried out some minor modifications to \texttt{Teko} in \texttt{Trilinos}\footnote{ Our implementation and results in this paper are based on \texttt{Trilinos} in the version with commit ID 23ccc58 (master branch); see the GitHub repository \cite{trilinos}.} for the use of PCD and LSC in the \texttt{FEDDLib}. 
Furthermore, the implementation of the pressure projection is not yet contained in the official \texttt{FROSch} repository.


\subsection{Trilinos and FEDDLib}

The \texttt{Trilinos} library~\cite{trilinos} is a collection of interoperable software packages for high-performance scientific computing. 
Since the list of packages is extensive, we will focus on those of relevance to our paper. 
The main building block for the preconditioners investigated in this paper is the parallel linear algebra package \texttt{Tpetra}, which provides parallel vector and matrix classes as well as functions that allow communication of those objects. 
It entails a software stack of packages that operate on \texttt{Tpetra} objects, such as direct and iterative solvers and preconditioners. 
\texttt{Tpetra} automatically incorporates access to the \texttt{Kokkos} performance portability framework~\cite{edwards_kokkos_2014,trott_kokkos_2021} and the corresponding kernel library \texttt{Kokkos Kernels}~\cite{rajamanickam_kokkos_2021}. 
As a result, functions of \texttt{Tpetra} objects allow for node-level parallelization on CPUs and GPUs via \texttt{Kokkos}' node type template parameter; 
here, we will only focus on distributed memory parallelization and, hence, we will not consider the effects of node-level parallelization via \texttt{Kokkos}. 
The older \texttt{Epetra} linear algebra package is currently being deprecated together with the whole \texttt{Epetra}-only stack of packages. 
Therefore, even though our software framework supports the \texttt{Epetra} stack, we employ only the \texttt{Tpetra} stack. 

The model problems are implemented and assembled using the \texttt{FEDDLib} \cite{feddlib}. 
Currently, it is split into three main packages, \texttt{amr}, \texttt{core}, and \texttt{problems}. 
The \texttt{amr} package is for adaptive mesh refinement, which is not employed in this work. 
The \texttt{core} package includes the general implementation of the finite element assembly, IO tools, the interface to the \texttt{Tpetra} parallel linear algebra, mesh classes, and further smaller utility functions. 
The \texttt{problem} package provides implementations of different physics problems supported by \texttt{FEDDLib}, including the computational fluid dynamics problems considered in this work. 
Further problems can be added easily using the interface to the Mathematica code generation package AceGen~\cite{korelc2016automation}; 
the interface is discussed in more detail in~\cite{balzani_computational_nodate}. 
Furthermore, \texttt{FEDDLib}'s \texttt{problem} package includes interfaces to the nonlinear solver package \texttt{NOX}, which will also be used for solving the nonlinear problems in this work, and the unified solver interface \texttt{Stratimikos}. 

\texttt{Stratimikos} enables the use of iterative solvers implemented in \texttt{Belos} and of interfaces for direct solvers, which are accessible via the \texttt{Amesos2} package. 
In order to accelerate the convergence of iterative solvers, \texttt{Trilinos} provides various preconditioners, including: 
\begin{itemize}[leftmargin=2.0em,topsep=4pt,itemsep=2pt]
    \item \texttt{Ifpack2}: one-level Schwarz preconditioners, supporting incomplete factorizations on the subdomains;
    \item \texttt{Teko}: block preconditioners;
    \item \texttt{FROSch}: multi-level Schwarz preconditioners;
    \item \texttt{MueLu}: algebraic multigrid, as a solver or preconditioner.
\end{itemize}
As mentioned before, we focus on \texttt{Teko} and \texttt{FROSch} preconditioners; 
see \Cref{section: Block Preconditioner,sec:schwarz}, respectively, for details on the preconditioners employed in this work. 
Details on \texttt{Teko} are given in~\Cref{Subsection:Teko} and on \texttt{FROSch} in \Cref{sec:frosch}.


\subsection{Fast and Robust Overlapping Schwarz (FROSch)} \label{sec:frosch}

\texttt{FROSch} (Fast and Robust Overlapping Schwarz)~\cite{Heinlein:2020:FRO,heinlein_parallel_2016a} is a parallel domain decomposition preconditioning package within the \texttt{Trilinos} software library~\cite{heroux_overview_2005,trilinos:2025:arxiv}. It builds on the Schwarz framework~\cite{toselli:2005:ddm} to construct preconditioners by combining elementary Schwarz operators, for instance, implementing the two-level Schwarz preconditioner by combining the first and second level operators; cf.~\cref{eq:one_level schwarz,eq:two_level schwarz}.

A key feature of \texttt{FROSch} is its use of extension-based coarse spaces that avoid explicit geometric information for the problems of this paper (cf.~\Cref{subsec: Coarse Spaces}) and, in particular, the need for coarse triangulations; see \Cref{sec:schwarz} for a detailed discussion of the methodology. 
The coarse spaces are based on the GDSW preconditioner~\cite{dohrmann_family_2008,dohrmann_domain_2008} and have been adapted for block systems~\cite{Heinlein_2019_MonoFluidFlow}. 
Moreover, \texttt{FROSch} also implements the reduced-dimension RGDSW variants~\cite{heinlein_improving_2018}. Multilevel extensions~\cite{heinlein_multilevel_2023} further enhance scalability, enabling \texttt{FROSch} to scale to over 220\,000 cores~\cite{heinlein_parallel_2023}.

\texttt{FROSch} integrates with \texttt{Trilinos}' unified solver interface, \texttt{Stratimikos}, and supports both parallel linear algebra frameworks, \texttt{Epetra} and \texttt{Tpetra}, via the lightweight \texttt{Xpetra} interface; due to the recent deprecation of Epetra, its support is merely a legacy feature. 
The overlapping domain decomposition is constructed algebraically using the sparsity pattern of the parallel distributed input matrix. With \texttt{Tpetra}, \texttt{FROSch} leverages the performance portability of \texttt{Kokkos} and \texttt{Kokkos Kernels} for efficient parallelization on CPUs and GPUs~\cite{yamazaki_experimental_2023}. 

\texttt{FROSch} preconditioners are compatible with iterative Krylov solvers from the \texttt{Belos} package~\cite{bavier_amesos2_2012}, either directly or through \texttt{Stratimikos}. 
For more details, see~\cite{heinlein_parallel_2016a,Heinlein_2019_MonoFluidFlow,Heinlein:2020:FRO}.


\subsection{Block preconditioners in Teko}
\label{Subsection:Teko}

\texttt{Teko} is a \texttt{Trilinos} package for block\hyp{}preconditioning multi-physics problems. 
It can be used to build generic preconditioning strategies for arbitrary problems, the simplest of which is a diagonal preconditioner. 
Specifically for the Navier--Stokes problem exist a number of block preconditioners within \texttt{Teko}. 
We will consider the block preconditioners presented in \Cref{section: Block Preconditioner} and different modifications that can be applied. 
The SIMPLE, SIMPLEC, LSC, and PCD preconditioners are part of the \texttt{NS} (Navier--Stokes) namespace within \texttt{Teko}. 
For a detailed description of \texttt{Teko} and the Navier--Stokes preconditioners within \texttt{Teko}, we refer to \cite{cyr_teko_2016,Cyr_teko_stabilized}.

The SIMPLE and SIMPLEC preconditioner can be constructed algebraically based on the underlying system to be solved, without the assembly of additional matrices. 
For LSC and PCD, to make components of the Schur complement that require additional construction available to \texttt{Teko}, callback functions are defined, e.g., for $A_p$ and $F_p$ of the PCD preconditioner. 
These operators are constructed outside of \texttt{Teko}, in our case in the \texttt{FEDDLib}. 
Only the callback function for the velocity mass matrix is optional; 
if it is not provided, a diagonal matrix based on $F$ is substituted as an approximation.

For the approximations $\hat{F}^{-1}$ and $\hat{S}^{-1}$ of the inverses of the fluid matrix and Schur complement, respectively, the specification of inverse types is necessary. 
If no inverse type is specified, the default is to use a direct solver from the \texttt{Amesos2} package. 
Ideally, a less expensive approximation strategy is specified. 
For our simulations, we will use different two-level overlapping Schwarz preconditioners from the \texttt{Trilinos} package \texttt{FROSch}; see \Cref{Table: approximations of system matrix used in numerical results} for an overview of the used combinations of preconditioning strategies.


\subsubsection{PCD block preconditioner}\label{Teko: PCD}

The PCD preconditioner in \texttt{Teko} is derived from the abstract class \texttt{LU2x2Strategy}. 
To use it, the parameter \texttt{Type} is set to \textit{Block LU2x2} and the \texttt{Strategy Name} to \textit{NS PCD Strategy}; see \Cref{Code2}. 
For the LSC and SIMPLE preconditioners, the \texttt{Type} directly refers to the block preconditioner, e.g., it is set to \textit{NS LSC} instead of \textit{LU2x2Strategy}; 
see \Cref{Code1,Code3} for the LSC and SIMPLE parameter files.

\begin{lstlisting}[float, floatplacement=tb, language=xml, caption={Pressure Convection--Diffusion (PCD): Parameter XML file for \texttt{Teko} (\texttt{Trilinos}).}, label=Code2]
<ParameterList name="PCD">
	<Parameter name="Type" type="string" value="Block LU2x2"/>
	<!-- Using LDU decomposition for preconditioning structure (true) or upper block-triangular (false) -->
	<Parameter name="Use LDU" type="bool" value="false"/>
	<Parameter name="Strategy Name" type="string" value="NS PCD Strategy"/>
	<ParameterList name="Strategy Settings">
		<Parameter name="PCD Operator" type="bool" value="true"/>
		<Parameter name="Pressure Laplace Operator" type="bool" value="true"/>
		<Parameter name="Pressure Mass Operator" type="bool" value="true"/>
		<!--FROSch preconditioner for approximation of inverses -->
		<Parameter name="Inverse F Type" type="string" value="FROSchVelocity"/>
		<Parameter name="Inverse Laplace Type" type="string" value="FROSchPressure"/>
		<!-- Approximation of inverse of pressure mass matrix: FROSch preconditioner (FROSchPressure) or diagonal approximation (Diagonal, AbsRowSum) -->
		<Parameter name="Inverse Mass Type" type="string" value="FROSchPressure"/>
		<Parameter name="Flip Schur Complement Ordering" type="bool" value="true"/> 
	</ParameterList>
</ParameterList>
\end{lstlisting}

\begin{table}[!tb]
\caption{\textbf{Comparison of block-triangular and LDU approach for PCD preconditioner.} Stationary simulation of backward-facing step (see \Cref{sect:problem:BFS}) with P2--P1 discretization, $H/h = 9$, and 1\,125 cores.
	PCD block preconditioner with boundary condition type (BC--2) and RGDSW--RGDSW coarse space. 
	Required number of nonlinear steps in parentheses. 
	See \Cref{Sec: Results} for parameters of the nonlinear and linear solver.
	Total time consists of setup and solve time.}
	\label{table: Block Comp PCD detailed LDU vs block triangular}
\begin{tabular*}{\linewidth}{@{}LRR@{} }
\toprule
									& Full LDU 				& Block-triangular   \\ \cmidrule{2-3}
 \#\,avg. iterations				 &	 \textbf{78} (5)					& 103 (5)						 	\\  \midrule
  Setup time 				 &	14.0\,\unit{\s}					&	14.4\,\unit{\s}						 	\\  
  Solve time		   			&	 77.1\,\unit{\s}				& 56.9\,\unit{\s}						 	\\ \cmidrule{2-3}
  Total 				 	&	91.1\,\unit{\s}					&	\textbf{71.3\,\unit{\s}}		    \\ \bottomrule
\end{tabular*}
\end{table}

For PCD, the order of the factors within the Schur complement approximation can be flipped by specifying the parameter \texttt{Flip Schur Complement Ordering}:
\begin{itemize}[leftmargin=2.0em,topsep=4pt,itemsep=2pt]
\item $S_{\text{PCD}}^{-1} = -M_p^{-1} F_p A_p^{-1}$, specified with \textit{false} (``default''),
\item $S_{\text{PCD}}^{-1} = -A_p^{-1} F_p M_p^{-1}$, specified with \textit{true} (``flipped''; cf. \cref{eq:PCD:Schur complement approximation}).
\end{itemize}
The ``default'' Schur complement ordering corresponds to an earlier version of the Schur complement approximation based on a commutator that uses a gradient operator on the pressure instead of a divergence operator on the velocity; 
see in \cite{Elman_BC_Navier_Stokes} (2.8) for the ``default'' and (2.2) for the ``flipped'' version; 
see also \cite[Remark~9.3]{Elman_fluid_book}. 

The Schur complement is constructed and defined as a linear operator. 
To provide all required operators for the PCD preconditioner, callback functions are defined for $M_p$, $A_p$, and $F_p$. 
\texttt{Teko} then creates approximations for $M_p^{-1}$ and $A_p^{-1}$ according to the parameter file.
The inverse of~$A_p$ is approximated with an overlapping Schwarz method via \texttt{FROSch}. 
This is indicated by \texttt{Inverse Laplace Type} set to \texttt{FROSchPressure}, where this value refers to another\footnote{See, for example, in the parameter file \path{feddlib/problems/examples/unsteadyNavierStokes/parametersTeko.xml} in \texttt{FEDDLib}'s GitHub repository \cite{feddlib} (commit ID 1246cdd), where a corresponding entry exists in the sublist \texttt{Inverse Factory Library}.} parameter list that contains the settings for \texttt{FROSch}.
The inverse of~$M_p$ will either be approximated via \texttt{FROSch} or by using the same approach as in LSC and SIMPLE for $M_u$ and $F$, respectively; that is, we can use a diagonal approximation of $M_p^{-1}$:
\begin{equation*}
H_{M_p}^D = \operatorname{diag}(M_p)^{-1}, \qquad 
H_{M_p}^{\Sigma} = \delta^{ij}\Big(\sum_{k=1}^{N_p} |(M_p)_{i,k}|\Big )^{-1},
\end{equation*}
referred to as \textit{Diagonal} and \textit{AbsRowSum} scaling, respectively. 
Note that the pressure convection--diffusion operator needs to be updated in each Newton or fixed-point iteration, while $A_p$ and $M_p$ do not change.

The PCD preconditioner introduced in \Cref{subsection: PCD} is a block-triangular preconditioner \cref{eq:block-tri precond}. 
In \texttt{Teko}, the PCD preconditioner can be applied as a full LDU block preconditioner (based on \cref{eq:Navier Stokes LDU}) or as a block-triangular preconditioner. 
The results in \Cref{table: Block Comp PCD detailed LDU vs block triangular} show, compared to the block-triangular preconditioner, how the LDU solve can improve the iteration count, but the total time is larger as its application is more costly. 
Consequently, we will use the classic block-triangular approach: 
\begin{enumerate}[leftmargin=2.5em,topsep=4pt,itemsep=2pt]
\item $p = \hat{S}_{\text{PCD}}^{-1} f_p$,
\item $\mathbf{u} = \hat{F}^{-1} (\mathbf{f}_u - B^{T} p)$.
\end{enumerate}
Here, $F^{-1}$ and $S_{\text{PCD}}^{-1}$ have been replaced with their approximations based on a Schwarz method, $\hat{F}^{-1}$ and $\hat{S}_{\text{PCD}}^{-1}$; 
see \Cref{Table: approximations of system matrix used in numerical results} for an overview of the approximations used in preconditioners. 

As covered in \Cref{Sec: PCD Boundary}, we need to additionally account for the boundary conditions that would naturally arise from the original Schur complement~$S$. 
The boundary information is integrated into $A_p$ and $F_p$ by the \texttt{FEDDLib}. 
When we set the Dirichlet boundary conditions in the different operators, we zero out the row and keep the diagonal entry; cf.~\cite[S325]{cyr_teko_2016}. 
The use of Robin boundary conditions additionally requires the assembly of a matrix on the surface elements of the respective boundary.


\subsubsection{LSC block preconditioner}\label{sect:teko:lsc}

To use the LSC preconditioner in \texttt{Teko}, the parameter \texttt{Type} is set to \textit{NS LSC}; see \Cref{Code1}. 
The LSC preconditioner, just like PCD, is an upper block-triangular preconditioner that is applied in two steps:
\begin{enumerate}[leftmargin=2.5em,topsep=4pt,itemsep=2pt]
\item $p = \hat{S}_{\text{LSC}}^{-1} f_p$,
\item $\mathbf{u} = \hat{F}^{-1} (\mathbf{f}_u - B^{T} p)$.
\end{enumerate}
The Schur complement $S_{\text{LSC}}$ is defined in \cref{Eq: S_LSC}. 
Similarly to before, $F^{-1}$ and $S_{\text{LSC}}^{-1}$ have been replaced with their approximations $\hat{F}^{-1}$ and $\hat{S}_{\text{LSC}}^{-1}$ based on a Schwarz method. 
We have different options to treat the inverse of $M_u$: 
$M_u$ can be passed along to \texttt{Teko} via a callback function. 
Inverting $M_u$ is expensive, but it can be well approximated with $H_{M_u}^D$ or $H_{M_u}^{\Sigma}$ (see \Cref{subsection: LSC} for the definition). 
Via the parameter file, $H_{M_u}^D$ can be selected with \textit{Diagonal} and $H_{M_u}^{\Sigma}$ with \textit{AbsRowSum}.
If $M_u$ is not provided via a callback function, by default it is replaced with $\operatorname{diag}(F)$. 
If $M_u$ is replaced with the identity matrix, 
it leads to the BFBt preconditioner in \cref{Eq: S_BFBT}, presented in \cite{Elman_OG_BFBT}.

\begin{lstlisting}[float, floatplacement=tb, language=xml, caption={Least-Squares Commutator (LSC): Pressure-Laplace Parameter XML file for \texttt{Teko} (\texttt{Trilinos}).}, label=Code1]
<ParameterList name="LSC-Pressure-Laplace">
	<Parameter name="Type" type="string" value="NS LSC"/>
	<Parameter name="Strategy Name" type="string" value="Pressure Laplace"/> 
	<ParameterList name="Strategy Settings">
		<!-- FROSch preconditioner for approximation of inverses -->
		<Parameter name="Inverse Velocity Type" type="string" value="FROSchVelocity"/>
		<Parameter name="Inverse Pressure Type" type="string" value="FROSchPressure"/>
		<!-- Diagonal approximation type of the inverse of the velocity mass matrix -->
		<Parameter name="Scaling Type" type="string" value="AbsRowSum"/>
		<!-- Use of mass matrix (true) or F -->
		<Parameter name="Use Mass Scaling" type="bool" value="true"/>
		<!-- Assuming stable (true) or stabilized  discretization (false) -->
		<Parameter name="Assume Stable Discretization" type="bool" value="true" />
 	</ParameterList>
</ParameterList>
\end{lstlisting}

Similarly to the PCD preconditioner, we also have the option to replace $B M_u^{-1} B^{T}$ for continuous pressure elements with the pressure-Laplacian $A_p$, the resulting Schur complement approximation we denote by $S_{\text{LSC}_{A_p}}$; see \cref{Eq: S_LSC_Ap}. 
In that case, we need to account for boundary conditions; 
we use the same method as in \Cref{Teko: PCD} to implement the boundary conditions in $A_p$ and always prescribe a Dirichlet condition on the outlet, which corresponds to the boundary condition strategy (BC--2) of PCD.

The preconditioner is part of the \texttt{NS} (Navier--Stokes) subpackage of \texttt{Teko} as \texttt{InvLSC}-\texttt{Strategy} or \texttt{PresLaplaceLSCStrategy}, where the \textit{Strategy Name} specifies the strategy (see \Cref{Code1}).

To compute the spectral radius of a matrix~$M$ in the stabilized variant, \texttt{Teko} uses a block Krylov--Schur method with block size 5 and tolerance 0.05 for the convergence criterion $\|M x-\lambda x\|_2 / \|\lambda\|_2$, where $x$ is an approximation of an eigenvector and $\lambda$ of the corresponding eigenvalue. 
This is carried out in each Newton iteration if the preconditioner is rebuilt.


\subsubsection{SIMPLE block preconditioner}

The \texttt{SIMPLE} preconditioner is also part of the \texttt{NS} sub package of \texttt{Teko}. 
It is selected by setting the parameter \texttt{Type} to \textit{NS SIMPLE}; see \Cref{Code3}. 
The type used to approximate the inverse of $F$ determines whether it is called SIMPLE or SIMPLEC: 
The type \textit{Diagonal} corresponds to SIMPLE and \textit{AbsRowSum} to SIMPLEC. 
The inverses of $F$ and the Schur complement are replaced with approximations $\hat{F}^{-1}$ and $\hat{S}_{\text{SIMPLE}}^{-1}$ based on a Schwarz method. 
The under-relaxation parameter is always set to $\alpha=0.9$; cf. \cite[p.~350]{Cyr_teko_stabilized}. 
Unlike in the case of PCD and LSC, the application of SIMPLE and SIMPLEC requires three steps:
\begin{enumerate}[leftmargin=2.5em,topsep=4pt,itemsep=2pt]
\item $\mathbf{u}^{*} = \hat{F}^{-1} \mathbf{f}_u$,
\item $p = \alpha \hat{S}_{\text{SIMPLE}}^{-1} ( f_p - B \mathbf{u}^{*})$,
\item $\mathbf{u} = \mathbf{u}^{*} - \frac{1}{\alpha} H_F B^{T} p$. 
\end{enumerate}

\begin{lstlisting}[float, floatplacement=tb, language=xml, caption={SIMPLEC: Parameter XML file for \texttt{Teko} (\texttt{Trilinos}).}, label=Code3]
<ParameterList name="SIMPLE">
	<Parameter name="Type" type="string" value="NS SIMPLE"/>
	<!--FROSch preconditioner for approximation of inverses -->
	<Parameter name="Inverse Velocity Type" type="string" value="FROSchVelocity"/>
	<Parameter name="Inverse Pressure Type" type="string" value="FROSchPressure"/>
	<!-- Definition of diagonal approximation type for H_F: Diagonal (SIMPLE), AbsRowSum (SIMPLEC), Lumped, BlkDiag-->
	<Parameter name="Explicit Velocity Inverse Type" type="string" value="AbsRowSum"/>
	<!-- Under-relaxation parameter -->
	<Parameter name="Alpha" type="double" value=".9"/>
</ParameterList>
\end{lstlisting}


\subsection{Comparison methodology}

In order to compare monolithic and block preconditioning for the Navier--Stokes equations, we will evaluate the time to set up the preconditioners, the time taken by the Krylov method (henceforth called ``solve'' time or cost), and the average iteration count per Newton step. 
The solve time includes the application of the preconditioner. 
The PCD and LSC preconditioners are block-triangular preconditioners and, thus, consist of two steps. 
On the other hand, the SIMPLE preconditioner contains the application of a backward and forward substitution and consists in total of three application steps. 
Independent of this, the block preconditioners contain approximations of inverses based on a Schwarz method. 
The monolithic preconditioner contains only the approximation of the inverse of $\mathcal{F}$ based on a Schwarz method. 
For the monolithic preconditioner, the setup time consists of the approximation of the inverse of $\mathcal{F}$, and for the block 
preconditioner it consists of the approximation of the inverses of $F$ and the Schur complement (which includes approximations of, e.g., $A_p^{-1}$ using a Schwarz method). 
\Cref{Table: Teko and Monolithic} illustrates the different operators that need to be approximated for the different preconditioning strategies. 
\Cref{Table: approximations of system matrix used in numerical results} gives an overview of the combinations of the different preconditioning strategies.

\begin{table}[!tb]
\caption{\textbf{Overview of approximations of inverses used for different preconditioning strategies.} 
$\operatorname{\mathbf{A_S}}\left[\cdot\right]$ denotes the application of an additive Schwarz preconditioner. 
In the case of a monolithic preconditioner, the entire system $\mathcal{F}^{-1}$ is approximated with a Schwarz method. 
In all other cases, the matrices $F^{-1}$ and $S^{-1}$ are approximated and then a block preconditioner like \cref{eq:block-tri precond} is applied.}\label{Table: approximations of system matrix used in numerical results}
\renewcommand{\arraystretch}{1.6}
\begin{tabular*}{\linewidth}{@{} LLC @{} }
\toprule
                        & \multicolumn{1}{c}{$\hat{S}^{-1}$}
                        & \multicolumn{1}{c}{$\hat{F}^{-1}$}
\\
\midrule
PCD                     & $-\AdditiveSchwarz\left[A_p^{-1}\right] 
                          F_p 
                          \AdditiveSchwarz\left[M_p^{-1}\right]$
                        & $\AdditiveSchwarz\left[F^{-1}\right]$
\\
SIMPLEC                 & $\phantom{-}\AdditiveSchwarz\Big[\big(-C-B H_F^{\Sigma} B^{T}\big)^{-1}\Big]$
                        & $\AdditiveSchwarz\left[F^{-1}\right]$
\\
LSC                     & $-\AdditiveSchwarz\Big[\big(B H_{M_u}^{\Sigma}B^T\big)^{-1}\Big]
                          \big(B H_{M_u}^{\Sigma} F H_{M_u}^{\Sigma} B^{T}\big)
                          \AdditiveSchwarz\Big[\big(B H_{M_u}^{\Sigma}B^T\big)^{-1}\Big]$
                        & $\AdditiveSchwarz\left[F^{-1}\right]$
\\
LSC$_{A_p}$             & $-\AdditiveSchwarz\left[A_p^{-1}\right]
                          \big(B H_{M_u}^{\Sigma} F H_{M_u}^{\Sigma} B^{T}\big)
                          \AdditiveSchwarz\left[A_p^{-1}\right]$
                        & $\AdditiveSchwarz\left[F^{-1}\right]$
\\
LSC$_{\text{stab},A_p}$ & $-\AdditiveSchwarz\left[A_p^{-1}\right] 
                          \big(B H_{M_u}^{\Sigma} F H_{M_u}^{\Sigma} B^{T}\big)
                          \AdditiveSchwarz\left[A_p^{-1}\right]
                          - \alpha D^{-1}$
                        & $\AdditiveSchwarz\left[F^{-1}\right]$
\\
\midrule
Monolithic              & \multicolumn{2}{c}{$\AdditiveSchwarz\left[\mathcal{F}^{-1}\right]$}
\\
\bottomrule
\end{tabular*}
\end{table}


\section{Numerical results}\label{Sec: Results}

To examine the performance of the different preconditioners, we consider two three-dimensional problem settings: 
the well-known backward-facing step and a realistic geometry of an artery from \cite{Artery_Balzani_2012}. 
We show results for stationary and transient simulations with volume force $\mathbf{f}=0$. 
In \Cref{subsect:monolithic preconditioner}, the monolithic preconditioner is tested for variations of the coarse spaces and further modifications. 
Subsequently, in \Cref{subsect:block preconditioners}, results of block preconditioners are shown, again for different coarse spaces. 
The preconditioners are always applied as right preconditioners.
Finally, in \Cref{subsect:results:comparison}, for specific choices of coarse spaces and modifications, the performance of monolithic and block preconditioners is compared.

The parallel results were obtained on the Fritz supercomputer at Friedrich-Alexander-Universität Erlangen-Nürnberg. Fritz has 992 compute nodes, each with two Intel Xeon Platinum 8360Y \textit{Ice Lake} processors and 256~GB of DDR4 RAM.

For the monolithic as well as the block preconditioners, different coarse spaces can be applied to the velocity and pressure components. The coarse spaces GDSW, GDSW\expStar{}, and RGDSW (see \Cref{subsec: Coarse Spaces}) are used. 
This offers the possibility to combine different coarse space combinations to improve convergence. 
Similarly to the notation of a discretization, for example, P2--P1, where the first term stands for the velocity and the second for the pressure discretization, we refer to the coarse spaces as, for example, GDSW\expStar{}--RGDSW. 
To shorten the notation, sometimes the abbreviations G for GDSW, R for RGDSW, and G\expStar{} for GDSW\expStar{} are used in figures. 
Different strategies for the reuse of parts of the preconditioner exist; 
they will be introduced based on the monolithic preconditioner but are equally used for block preconditioners; see \Cref{Subsubsection: Setup}. 
Note that we employ a scaled additive Schwarz method in the first level; see \Cref{sec:schwarz}.

Newton's method terminates when the relative residual $\|\Residual(X^{k+1})\|_2/\|\Residual(X^0)\|_2$ or the Newton update $\|X^{k+1}-X^k\|_2$ reach a tolerance of $10^{-8}$. 
The initial guess for Newton's method is $X^0=0$. 
For time-dependent problems this refers only to the first time step. 
Newton's method is used with GMRES along with an adaptive forcing term $\eta_k$ as described in \Cref{Sec: Inexact Newton + Forcing Term}. 
The initial forcing term is set to $\eta_0 = 10^{-3}$, the minimum forcing term to $\eta_{\min}=10^{-8}$, the maximum forcing term to $\eta_{\max}=10^{-3}$, and the constants $\alpha=1.5$ and $\gamma=0.9$; cf.~\cref{eq: forcing term bounds}. 
In each nonlinear iteration, the GMRES iteration is started with zero as the initial guess.

For the Newton method in \texttt{NOX}, we employ a line search with backtracking as the globalization strategy; see \cite[Section 5.4]{diss_hochmuth_2020} for more details and \cite{InexactNewton} for a full derivation of the method and globalization techniques. 
The starting step length is set to~1. 
The other required parameters are left at their respective default values; see the \texttt{NOX} documentation \cite{trilinos-website}. 

The initial mesh partition defining the nonoverlapping subdomains is either structured and constructed manually, or it is constructed by METIS \cite{karypis:1999:metis}, which results in an unstructured partition of the mesh. 
The subdomains are of similar diameter, which we denote by~$H$; i.e., $H\approx\operatorname{diam}(\Omega_i)$.
We choose an algebraically determined overlap of $\hat{\delta}=1$ for the construction of the local overlapping stiffness matrices; cf. \Cref{sec:schwarz}.

For the transient simulations we apply the BDF-2 time stepping scheme; see \cref{Eq: Sys with BDF2} and \cite{BDF2}.

Due to comprehensive legends and space restrictions, we occasionally show the legend only in one subfigure if it is valid for all subfigures of the respective figure.


\subsection{Test cases}

We use two geometries to test the preconditioners: a backward-facing step geometry with a structured mesh that is decomposed into cube-shaped subdomains and a realistic artery with an unstructured mesh and an unstructured domain decomposition.


\subsubsection{Test case 1: Backward-Facing Step (BFS)}
\label{sect:problem:BFS}

\begin{figure}[!tb] 
    \centering
    \begin{minipage}[b]{0.49\textwidth}%
        \centering
        \includegraphics[width=0.9\textwidth]{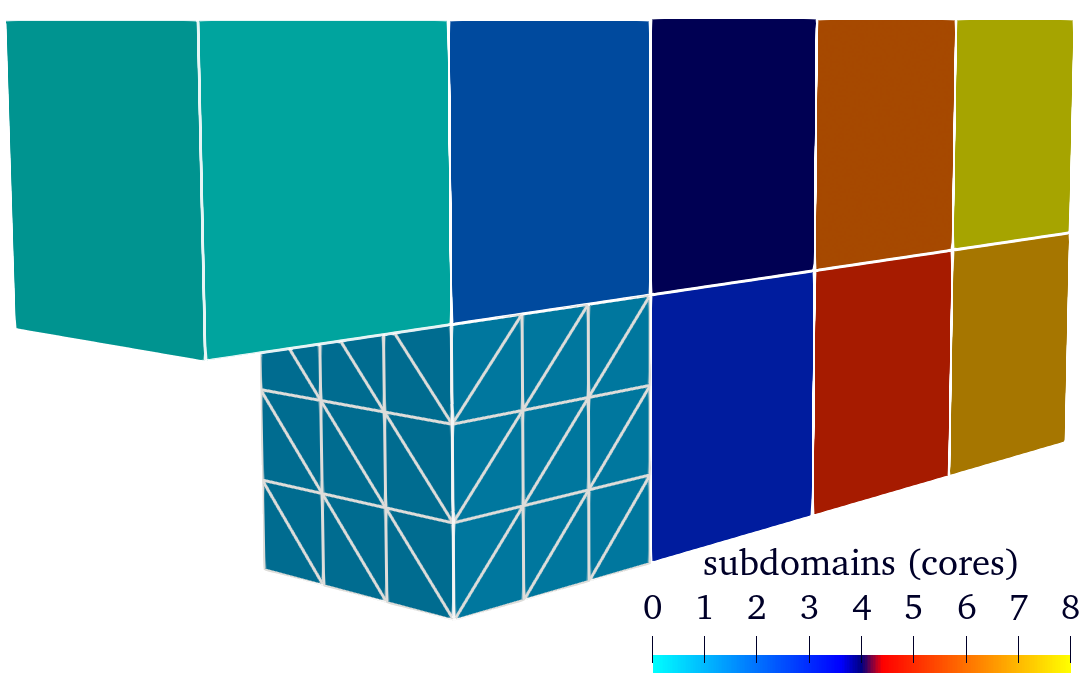}
    \end{minipage}%
    \begin{minipage}[b]{0.49\textwidth}%
        \centering
       \includegraphics[width=0.9\textwidth]{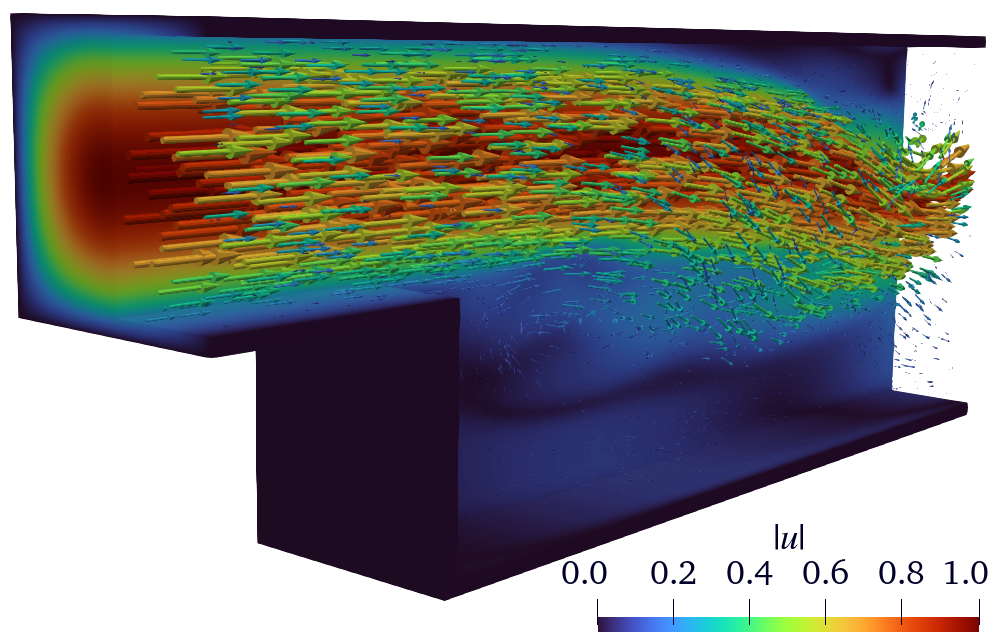}
    \end{minipage}%
    \caption{
    \textbf{Test case 1: Backward-facing step (BFS) problem.}    
    \textbf{Left:} Backward-facing step geometry with structured mesh partition into 9 subcubes, each with side length \qty{1}{\cm}. 
    The subcubes correspond to the subdomains, and in the case displayed, we have $H/h=3$. 
    For a higher number of processor cores, the subcubes are divided up further. 
    For example, for 243 processor cores, each subcube holds $3\times 3 \times 3$ subdomains. 
    \textbf{Right:} Solution to transient BFS problem at $t=\qty{10.0}{\s}$, using time step size $\Delta t = \qty{0.05}{\s}$. 
    Reynolds number is 3\,200 with kinematic viscosity $\nu=\qty[per-mode = symbol]{6.25e-4}{\cm^2\per\s}$. 
    P2--P1 discretization with $H/h=9$, computed on 243~cores.
    }\label{figure1: Boundary conditions BFS and Inflow Profiles}
\end{figure}

\begin{figure}[!tb]
    \centering
    \includegraphics[width=0.7\textwidth,trim={0.25cm 0.15cm 0.15cm 0.0cm},clip]{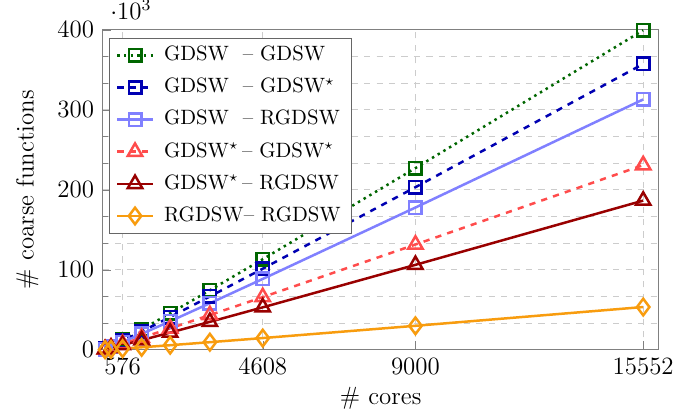}%
    \caption{
        \textbf{Number of coarse functions for different number of processor cores (which equals the number of subdomains) for the backward-facing step (BFS) problem.}
        The first name corresponds to the interface partition of unity for the velocity component and the second for the pressure component.}\label{figure: Coarse functions BFS}
\end{figure}

The three-dimensional backward-facing step (BFS) test case is used with a structured mesh and a structured mesh partition; see \Cref{figure1: Boundary conditions BFS and Inflow Profiles} (left), and see (right) for the solution of a sample case. 
%
We test the weak scaling properties of the monolithic and block preconditioners for stationary and transient flow for this geometry. 
To this end, the characteristic subdomain size $H/h$ is kept constant, while the number of subdomains increases with the number of processor cores.
For the different finite element discretizations and subdomain resolutions, the number of degrees of freedom varies; see \Cref{Table:BFS_Hh_Dofs}. 
The corresponding coarse space dimensions are shown in \Cref{figure: Coarse functions BFS} for a selection of coarse space combinations. 

\begin{table}[!tb]
\caption{\textbf{Subdomain resolution and number of degrees of freedom of the backward-facing step geometry for a structured mesh partition.} Subdomain resolution $H/h$ and number of degrees of freedom for different discretizations and a varying number of subdomains~$N$ (processor cores).}\label{Table:BFS_Hh_Dofs}
\begin{tabular*}{\linewidth}{@{} LRRRRR@{} }
\toprule
			& 	 & 	\multicolumn{4}{c}{Number of degrees of freedom} \\
			\cmidrule{3-6} 
Discretization & $H/h$ & $N=243$           & $N=1\,125$       & $N=4\,608$       & $N=15\,552$  \\ \toprule
Q1--Q1 / P1--P1 & 17 & $4.9\cdot 10^{6}$ & $22\cdot 10^{6}$ & $92\cdot 10^{6}$ & $308\cdot 10^{6}$ \\
Q2--Q1 / P2--P1 &  9 & $4.5\cdot 10^{6}$ & $21\cdot 10^{6}$ & $85\cdot 10^{6}$ & $285\cdot 10^{6}$ \\
Q2--P1-disc.    &  9 & $5.1\cdot 10^{6}$ & $23\cdot 10^{6}$ & $95\cdot 10^{6}$ & $320\cdot 10^{6}$ \\
\bottomrule
\end{tabular*}
\end{table}

We prescribe a parabolic-like inflow profile on the inlet boundary of the domain $\partial \Omega_{\text{in}}$ with a maximum inflow velocity of $v_{\max}$ and an inflow height $L=\qty{1}{\cm}$:
\begin{equation*}
\mathbf{u}(0,y,z) = \left[ \big ( 16 v_{\max} y z (L-y) (L-z)\big) \big / L^4 , 0 , 0 \right],\quad (y,z)\in[0,L]^2.
\end{equation*}
To define the Reynolds number of the BFS problem, we use the outlet height of \qty{2}{\cm} as the characteristic length and the maximum inflow velocity $v_{\max}$ as $v$ in \cref{eq:Reynolds number:general}:
\begin{equation}\label{Eq: Re BFS}
    \Reynum = \frac{2 v_{\max}}{\nu}.
\end{equation}
For $v_{\max}=\qty[per-mode=symbol]{1}{\cm\per\s}$ and a kinematic viscosity of $\nu=\qty[per-mode=symbol]{0.01}{\cm^2\per\s}$, the Reynolds number is 200. 
If not stated otherwise, this is the default configuration for the BFS problem.


\subsubsection{Test case 2: realistic artery} \label{subsubsec: test case realistic artery}

\begin{figure}[!tb]
    \centering
    \begin{minipage}[b]{0.55\textwidth}%
        \includegraphics[width=\textwidth]{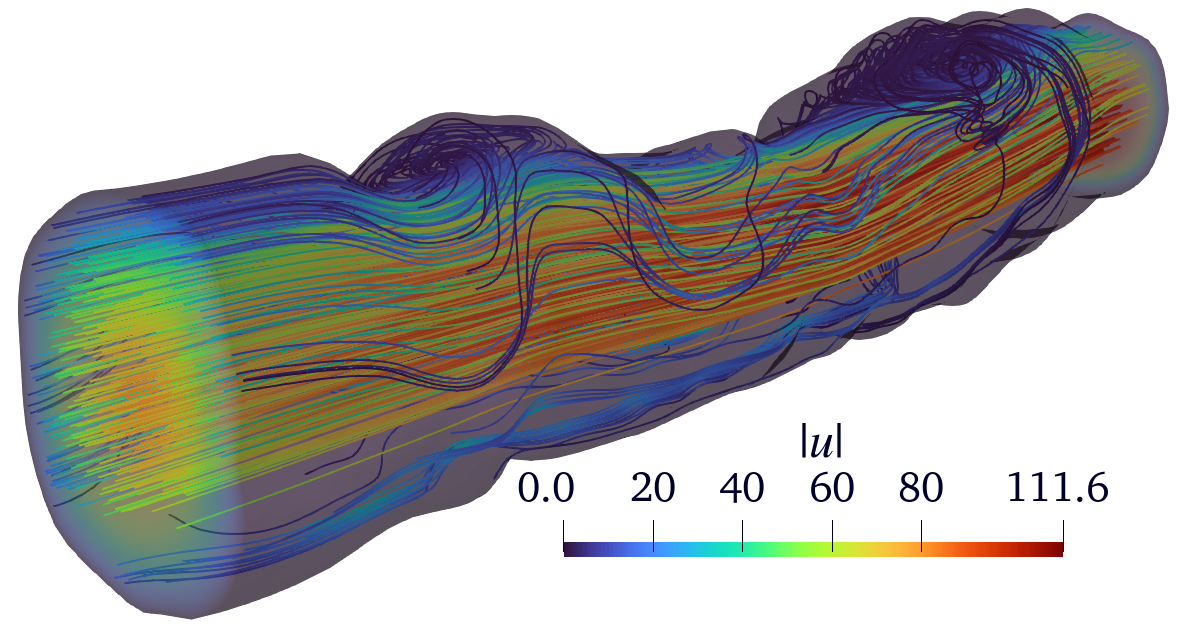}
    \end{minipage}%
    \hspace*{0.025\textwidth}
    \begin{minipage}[b]{0.4\textwidth}%
        \includegraphics[width=\textwidth]{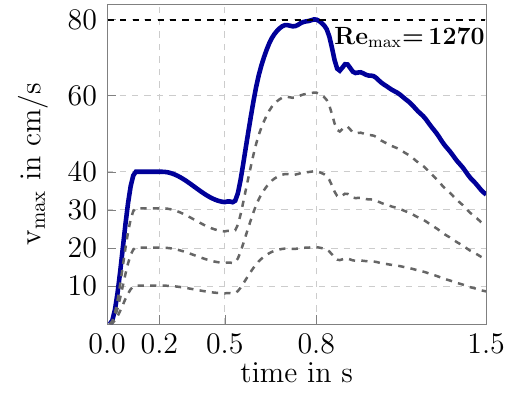}
    \end{minipage}%
    \caption{
	\textbf{Test case 2: Realistic artery problem.}    
	Geometry taken from \cite{Artery_Balzani_2012}.
    \textbf{Left:} Magnitude of solution $\mathbf{u}_h$ at $t=\qty{0.8}{\s}$ with initial maximum inflow velocity $v_{\max}=\qty[per-mode=symbol]{40}{\cm\per\s}$. 
		\textbf{Right:} Inflow profile over time of the realistic artery problem. 
        Ramp phase until $t=\qty{0.1}{\s}$ and subsequent initial maximum inflow velocity $v_{\max}=\qty[per-mode=symbol]{40}{\cm\per\s}$. 
		Maximum Reynolds number calculated based on maximum velocity and on approximate maximum diameter of artery; see \cref{Def: Re artery max}. 
		The heartbeat and the corresponding flow rate were constructed in \cite{balzani:2015:nmf} based on inflow pressure data from \cite{Hemolab}.}\label{figure2: realisic artery}
\end{figure}

For the second test case, we consider a tetrahedralized geometry of the realistic artery in \Cref{figure2: realisic artery} (left) from \cite{Artery_Balzani_2012} with 170\,000 P1 nodes and 915\,336 elements. 
The full fluid flow problem with a P2--P1 discretization has 4.1~million degrees of freedom. 
The focus of this test case is on transient flow simulations where the flow profile contains one heart beat; see \Cref{figure2: realisic artery} (right).

The interior (lumen) diameter of the artery varies between \qty{0.2}{\cm} and \qty{0.34}{\cm}. 
The full geometry from \cite{Artery_Balzani_2012} also contains layers of plaque, media, and adventitia (not shown here). 
The (outer) diameter of this artery is relatively even and approximately between \qty{0.44}{\cm} and \qty{0.47}{\cm}. 
This is in the range of mean diameters of renal arteries as stated in \cite[Table~2]{renal_artery_diam}; 
see also \cite[Sect.~3.1]{renal_artery_flow}, where a catheter with an interior diameter of \qty{0.2}{\cm} was used to emulate a renal artery. 
Motivated by this correlation, we prescribe a flow rate and a range for the Reynolds number that is realistic for this kind of artery. 
The inflow profile at the inlet $\partial \Omega_{\text{in}}$ of the domain is constructed by solving a Poisson problem on the inlet. 
Since the inlet is fairly round, this approximately yields a parabolic flow profile. 
The proposed blood flow rates in \cite[Sect.~3.3]{renal_artery_flow} range between \qty[per-mode=symbol]{0.1}{\cm^3\per\s} and \qty[per-mode=symbol]{4.0}{\cm^3\per\s} and are constant over time. 
For the (not perfectly round) inlet of our geometry, these flow rate values translate to a maximum inflow velocity ranging from 3.0 to \qty[per-mode=symbol]{120.0}{\cm\per\s}. 

The kinematic viscosity is set to $\nu=\qty[per-mode=symbol]{0.03}{\cm^2\per\s}$. 
Since the proposed artery has no fixed diameter, and the velocity field throughout the artery varies, we will estimate the maximum Reynolds number by
\begin{equation}\label{Def: Re artery max}
\Reynum_{\text{art}}(t) = \frac{D_{\max} \cdot u_{\max}}{\nu},
\quad \Reynum_{\text{art},\max} = \max_t \Reynum_{\text{art}}(t),
\end{equation}
where $D_{\max}$ is the maximum diameter, and $u_{\max}$ is the maximum velocity magnitude $\vert \mathbf{u} \vert$ attained in a finite element node. 
Consequently, with the maximum measured diameter of \qty{0.34}{\cm} and $u_{\max} = \qty[per-mode=symbol]{112}{\cm\per\s}$ for a maximum initial inflow velocity of $v_{\max}=\qty[per-mode=symbol]{40}{\cm\per\s}$ (after an initial ramp phase until \qty{0.1}{\s}), we have $\Reynum_{\text{art},\max} = 1\,270$; see \Cref{figure2: realisic artery} (left) for a solution plot.

In \Cref{Table: Velocity and Re in artery}, based on a selection of initial flow velocities, the corresponding maximum Reynolds numbers and time steps are given. This corresponds to measurements taken from the renal artery of 10~patients in \cite[Table~2]{renal_artery_Re}. The peak Reynolds number of the patients ranged between 870 and 1\,320, with an average value of 1\,145.

Using the realistic artery and flow problem, we can further analyze the solvers' abilities to adapt to a large range of Reynolds numbers and to complex flows. 
Furthermore, strong scaling is tested.

\begin{table}[!tb]
\caption{
\textbf{Detailed configuration for the realistic artery problem (Test case 2).}
Configurations for the maximum initial inflow velocity, the time step size, and approximate maximum Reynolds number. The maximum Reynolds number for this test is defined as $\Reynum_{\text{art},\max}$ in \cref{Def: Re artery max}.}
\label{Table: Velocity and Re in artery}
\begin{tabular*}{\linewidth}{@{} LRRRR@{} }
\toprule
$v_{\max}$ 													 &  \qty[per-mode=symbol]{10}{\cm\per\s}		& \qty[per-mode=symbol]{20}{\cm\per\s}		& \qty[per-mode=symbol]{30}{\cm\per\s}	   		& \qty[per-mode=symbol]{40}{\cm\per\s} \\ \midrule
$\Delta t$ 																&  0.005\,\unit{\s} 			& 0.0025\,\unit{\s} 	&  	0.001875\,\unit{\s} & 0.00125\,\unit{\s} \\
$\Reynum_{\text{art},\max}$ 		&   345				&	665					&	970			   &	1\,270				\\
\bottomrule
\end{tabular*}
\end{table}


\subsection{Monolithic preconditioner and coarse-problem recycling strategies}
\label{subsect:monolithic preconditioner}

The monolithic overlapping two-level Schwarz preconditioner is used for the saddle-point problem~\cref{Eq: Saddle point problem Navier Stokes} as described in \Cref{Section: OSP: Monolithic}. 
In \cite{heinlein_reduced_2019} the GDSW and RGDSW coarse spaces were evaluated and compared for three-dimensional Navier--Stokes problems. 
We will take the authors' findings into account about strategies for reusing parts of the preconditioner (cf. \cite[Sec.~5.4]{heinlein_reduced_2019}). 
In particular, we employ the strategy of reusing the coarse basis functions, which produced the lowest total time -- consisting of solve and setup time -- in \cite{heinlein_reduced_2019}. 
The problems of the first level and the coarse problem of the Schwarz preconditioner are solved sequentially; 
for example, the $N$ subdomain problems of the first level are solved on $N$ processor cores, and subsequently a subset of these cores are used to solve the coarse problem in parallel. 
We do not employ a hybrid Schwarz approach but use a fully additive one; see \Cref{sec:schwarz}.
Further improvements with respect to the compute time can be achieved by using additional cores to solve the coarse problem; cf. \cite[Sec.~4.4]{heinlein_reduced_2019}. 
The number of processor cores assigned to the coarse solve depends on the coarse space dimension (shown in \Cref{figure: Coarse functions BFS} for the backward-facing step problem) and ranges between 2 and 16; see \Cref{Table: coarse procs} for details in case of the BFS problem. 

\begin{table}[!tb]
\caption{\textbf{Number of cores assigned to solve the coarse problem for the results shown in \Cref{figure: Average time per Newton step monolithic stationary all discs} for the stationary BFS problem.}
	Specifically, this corresponds to the P2--P1 discretization with the GDSW$^\star$--RGDSW coarse space combination. 
	For RGDSW--RGDSW fewer cores are used due to the lower number of coarse functions; compare with \Cref{figure: Coarse functions BFS}.}\label{Table: coarse procs}
\begin{tabular*}{\linewidth}{@{} LRRRRRRRR@{} }
	\toprule						
	\#\,cores			 & 243 	& 576  &		1\,125		&		1\,944 		 		  	&		3\,087  			&	 4\,608 &  9\,000    & 15\,552   \\ \midrule
	\#\,cores for coarse solve		&	4		&		4	&			6			&			6						&			8				&  12 		& 16 		& 16 \\ \bottomrule
\end{tabular*}
\end{table}


\subsubsection{Recycling strategies}\label{Subsubsection: Setup}

\Cref{figure: Coarse functions BFS} shows the dependence of the number of coarse functions on the used coarse spaces. Generally, the higher the number of coarse functions is, the higher the setup costs are. 
The reuse strategies described in \cite{heinlein_reduced_2019} and included in the software package \texttt{FROSch} are 
 the reuse of the symbolic factorization of $\mathcal{F}_0$ (SF),
 of the coarse basis functions $\phi$ (CB), and
of the coarse matrix $\mathcal{F}_0$ (CM).

The more extensive the reuse strategy is, the more it affects the linear iteration count. 
The reuse of the symbolic factorization (SF) is always applied. 
In general, the preconditioner is rebuilt or updated in each Newton iteration.
If the coarse basis is reused (CB), only the coarse matrix is recomputed in each Newton iteration. 

Additionally to the previously described reuse strategies, the setup cost can be reduced further by only rebuilding or updating the preconditioner in the first $k$ Newton iterations (of each time step or of the stationary simulation), and then reusing it. 
This is crucial for the monolithic preconditioner, where the setup cost is higher compared to block-preconditioning strategies.

\begin{table}[!tb]
\caption{\textbf{Comparison of setup and solve times of the monolithic preconditioner if recomputing it the first $k$ Newton iterations.}
Reuse of symbolic factorization (SF) and coarse basis (CB). 
Recompute preconditioner (except for (SF) and (CB)) up until $k$ Newton iterations, reuse entire preconditioner based on $k$th iteration afterwards.
Stationary BFS problem with P2--P1 discretization with $H/h=9$. Average GMRES iteration count per Newton step. 
Total number of Newton iterations required to reach convergence is~5. 
RGDSW--RGDSW coarse space combination for the velocity (first component) and pressure (second component). 
}
\label{Table: MonoStationary_SetupReuse_R_R}
\begin{tabular*}{\linewidth}{@{} LLRRRR@{} }
	\toprule						
			&$k$															&		5						&		4		 		  	&		3  			&	 2 \\ \midrule
	\multirow{4}{*}{\shortstack{RGDSW--RGDSW}}	
												&\#\,avg. iter. 	 	&  52.8			&   52.8 					&   54.0	 				& 78.4	\\ 
												&Setup cost 		 	&  29.4\,\unit{\s} 		&  24.3\,\unit{\s}				& 19.4\,\unit{\s}					& 14.8\,\unit{\s} \\
												&Solve cost 		 	&  35.9\,\unit{\s} 		&  34.9\,\unit{\s}				&  35.9\,\unit{\s}  	 			&  56.1\,\unit{\s} \\																											&Total					    &	65.3\,\unit{\s}			&	59.2\,\unit{\s}					&	\textbf{55.3\,\unit{\s}} &	70.9\,\unit{\s}	\\ \midrule
\end{tabular*}
\end{table}

Based on the results in \Cref{Table: MonoStationary_SetupReuse_R_R} for the stationary BFS case, we can save computing time by reusing the complete preconditioner after three Newton steps. 
For $k=3$ we achieve the best total time, since the linear iteration count is only slightly increased, while the setup time is small. 
The reuse of the coarse basis is not always the best choice, for example, for the GDSW--GDSW coarse space combination, since it increases the solve time to the point where it outweighs the time savings in the setup part. 
As it is preferable to use the RGDSW coarse space due to its smaller dimension, for the stationary tests we will employ the CB reuse strategy and reuse the preconditioner in the Newton method after iteration $k=3$. 
In transient simulations it is sufficient to update the preconditioner only in the first ($k$=1) Newton step of each time step, probably since the required average number of Newton iterations per time step is only approximately three.
Furthermore, the CB reuse strategy is applied in transient simulations over all time steps; that is, the coarse basis is only computed once for the entire simulation, but due to the changing system matrix, the coarse matrix and, thus, the preconditioner are updated once in each time step. 
A justification for the reuse of the coarse basis is also given later in \Cref{figure: Iterations CFL comparison}, which shows a direct comparison of reusing and not reusing the coarse basis for building the preconditioner. 
It reveals only a slight difference with respect to the average iteration count. 

The stated recycling strategies are also applied to the block preconditioners. 
Additionally, for the LSC and PCD block preconditioners, we reuse the approximations of $A_p^{-1}$ and $M_p^{-1}$ that are used for the construction of the Schur complement approximation, since they are independent of the solution and only need to be built once.

Regarding the construction of the coarse matrix for the monolithic two-level overlapping Schwarz preconditioner, we can additionally choose between using the unaltered coarse basis \cref{eq:phi_mono} (\emph{full $\phi$}) or excluding the off-diagonal blocks:
\begin{equation*}
\phi \leftarrow
\begin{bmatrix}
\phi _{\mathbf{u},\mathbf{u_0}} & 0 \\
0 & \phi_{p,p_0}
\end{bmatrix}.
\end{equation*}
This can have an effect on the scaling behavior of the preconditioner and will be further discussed in the following sections.


\subsubsection{Stationary Navier--Stokes problem}\label{Subsubsection:Mono_Stationary}

Different finite element discretizations are considered for the stationary BFS problem. 
First, we consider the P1--P1 discretization with the Bochev--Dohrmann stabilization and the inf-sup stable P2--P1 discretization for the comparison of different coarse space combinations and to show results of the new GDSW$^\star$ coarse space. 
For the inf-sup stable discretizations, formulation \cref{eq: Pressure Projection} with the pressure projection is used. 
For results omitting the pressure projection in case of stable discretizations, see \Cref{figure: Iterations monolithic stationary comparison without pP} in the appendix.

In \Cref{figure: Iterations monolithic stationary} we observe that P1--P1 elements yield similar results for all coarse space combinations, independent of using the full $\phi$ or excluding the off-diagonal blocks to construct the coarse matrix. 
Weak scaling can be observed for up to 15\,552 cores;
see also \Cref{figure: Iterations monolithic stationary all discs,figure: Average time per Newton step monolithic stationary all discs} for results of RGDSW--RGDSW up to 15\,552 cores. 
For the stable P2--P1 element, in contrast, we observe a substantial influence in \Cref{figure: Iterations monolithic stationary comparison} for some combinations of coarse spaces and using the full or decoupled $\phi$.

\begin{figure}[!tb]
    \centering
    \begin{minipage}[t]{0.49\textwidth}%
        \includegraphics[width=\textwidth,trim={0.3cm 0.15cm 0.15cm 0.15cm},clip]{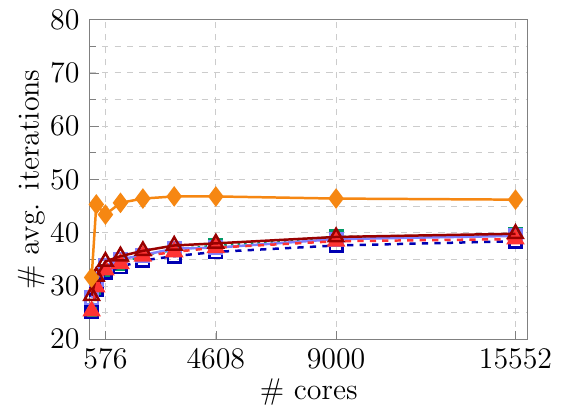}
    \end{minipage}%
    \hspace*{1pt}
    \begin{minipage}[t]{0.49\textwidth}%
        \includegraphics[width=\textwidth,trim={0.3cm 0.15cm 0.15cm 0.15cm},clip]{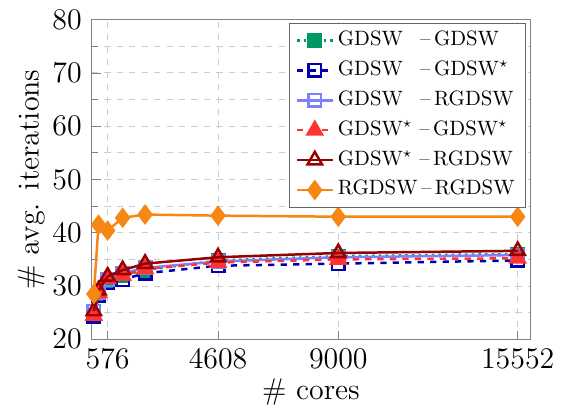}
    \end{minipage}%
    \caption{\textbf{Comparison of different coarse space combinations for the velocity (first component) and pressure (second component) for a P1--P1 discretization and the monolithic preconditioner.}
    Average GMRES iteration count per Newton step. 
    Total number of Newton iterations required to reach convergence is~5. 
    Weak scaling test. 
    Stationary BFS problem with $H/h=17$. 
    $\Reynum=200$. 
    Pressure projection is not used. 
    \textbf{Left:} Excluding off-diagonal blocks in $\phi$ to build the coarse matrix.
    \textbf{Right:} Using full $\phi$ to build the coarse matrix.
    }\label{figure: Iterations monolithic stationary}
\end{figure}

\begin{figure}[!tb]
    \begin{minipage}[t]{0.495\textwidth}%
        \includegraphics[width=\textwidth,trim={0.3cm 0.15cm 0.15cm 0.15cm},clip]{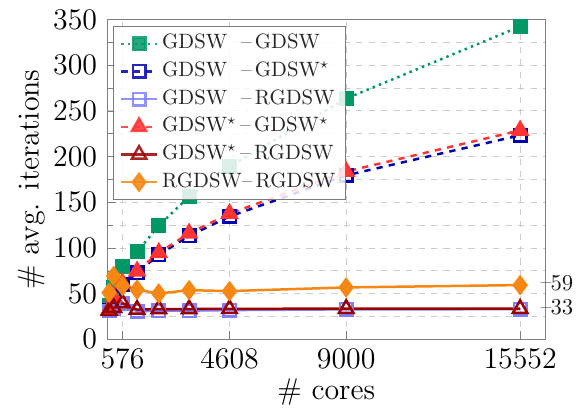}
    \end{minipage}%
    \hspace*{\fill}
    \begin{minipage}[t]{0.495\textwidth}%
        \includegraphics[width=\textwidth,trim={0.3cm 0.15cm 0.15cm 0.15cm},clip]{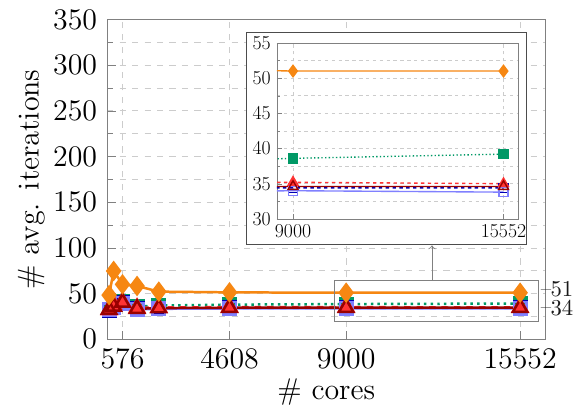}
    \end{minipage}%
    \caption{\textbf{Comparison of different coarse space combinations for the velocity (first component) and pressure (second component) for a P2--P1 discretization and the monolithic preconditioner.}
    Average GMRES iteration count per Newton step. 
    Total number of Newton iterations required to reach convergence is~5. 
    Weak scaling test. 
    Stationary BFS problem with $H/h=9$. 
    $\Reynum=200$. 
    Pressure projection included in the first level; 
    see \Cref{figure: Iterations monolithic stationary comparison without pP} in the appendix for results if the pressure projection is not used. 
    See \Cref{figure: Iterations monolithic stationary ldc} in the appendix for results of a lid-driven cavity problem.
    \textbf{Left:} Excluding off-diagonal blocks in $\phi$ to build the coarse matrix.
    \textbf{Right:} Using full $\phi$ to build the coarse matrix.}
    \label{figure: Iterations monolithic stationary comparison}
\end{figure}

We focus first on the results in \Cref{figure: Iterations monolithic stationary comparison} (left), where the off-diagonal blocks of $\phi$ are excluded. 
It stands out that, if RGDSW for the pressure is combined with a larger coarse space like GDSW or GDSW$^\star$ for the velocity, we achieve the lowest average GMRES iteration count. 
If a smaller coarse space for the velocity is used than for the pressure, the algorithm does not reach the GMRES convergence criterion (results are not shown here). 
The comparison to an inf-sup condition of a finite element discretization seems natural. 
On the other hand, the ``equal-order'' approach RGDSW--RGDSW is scalable with a low number of iterations, albeit not the lowest one. 
The cause for the correlation of the iteration count and the relative size of the velocity and pressure coarse spaces may also be independent of an inf-sup condition; further research is required. 

For the results using the full $\phi$ in \Cref{figure: Iterations monolithic stationary comparison} (right), we observe similar good scaling results but for all coarse space combinations. 
If the full $\phi$ is used to construct the coarse matrix, it seems imperative that also the pressure projection is applied in the first level; 
if the pressure projection is omitted, we lose scalability for all combinations in \Cref{figure: Iterations monolithic stationary comparison without pP} (right). 
On the other hand, the analogue to \Cref{figure: Iterations monolithic stationary comparison} (left) is \Cref{figure: Iterations monolithic stationary comparison without pP} (left), which shows the same general scaling behavior; that is, the pressure projection does not seem to affect the case in which $\phi$ is decoupled.
Independent of the use of the decoupled or full $\phi$, GDSW\expStar{}--RGDSW is the best performing configuration for the P2--P1 finite element discretization if the pressure projection is used. 
It remains to be thoroughly investigated, however, to determine how the various settings and coarse spaces impact each other. 
\Cref{figure: Iterations monolithic stationary ldc} in \Cref{appendix:stationary:monolithic} qualitatively shows the same results as \Cref{figure: Iterations monolithic stationary comparison} for a lid-driven cavity, indicating that the results are not specific to the backward-facing step problem.

\Cref{Table:MonoStationary_RuntimeCompP2} shows details corresponding to \Cref{figure: Iterations monolithic stationary}~(right) (P1--P1 discretization) and \Cref{figure: Iterations monolithic stationary comparison}~(right) (P2--P1 discretization) for a selection of coarse space combinations. 
Even though the average iteration count for the RGDSW--RGDSW combination is highest for P1--P1, it offers the lowest total cost due to its low setup cost. 
Notably, the setup cost for GDSW is higher than for GDSW$^\star$. The lowest setup cost offers the RGDSW coarse space. 
In \Cref{Table:MonoStationary_RuntimeCompP2} we observe how GDSW$^\star$--RGDSW performs best with respect to the total time (setup+solve) for the P2--P1 discretization. 

\begin{table}[!tb]
\caption{\textbf{Comparison of setup and solve time of the three best coarse space combinations for P1--P1 and P2--P1 discretization on 4\,608 cores.} Stationary BFS problem with configurations according to \Cref{Table:BFS_Hh_Dofs}. Total number of Newton iterations required to reach convergence is~5. Use of Bochev--Dohrmann stabilization for P1--P1 elements. Use of full $\phi$. Use of pressure projection for P2--P1. $\Reynum=200$. 
}\label{Table:MonoStationary_RuntimeCompP2}
\begin{tabular*}{\linewidth}{@{} LLRRR@{} }
\toprule
				&																		 & 			GDSW--RGDSW & GDSW$^\star$--RGDSW 	& RGDSW--RGDSW  \\ \toprule
				&\#\,avg. iterations 												  &					34.6  	   		& 35.4										&  43.2 \\ \midrule
\multirow{3}{*}{\shortstack{P1--P1}}	&Setup time &					26.7\,\unit{\s}	   		& 17.5\,\unit{\s}									&  12.8\,\unit{\s} \\
																	&	Solve time &					21.3\,\unit{\s}	   		& 18.4\,\unit{\s}									&  22.3\,\unit{\s} \\
																	&Total			 &						48.0\,\unit{\s}	   		& 	35.9\,\unit{\s}				    			& \textbf{35.1\,\unit{\s}} \\ \toprule

				&\#\,avg. iterations						 &					    34.4	   			&  34.6										&   51.4\\ \midrule
\multirow{3}{*}{\shortstack{P2--P1}}	&Setup time &					 47.6\,\unit{\s}	   	&  23.7\,\unit{\s}									&   20.2\,\unit{\s} \\
																	&	Solve time &					 27.3\,\unit{\s}	   	&  22.5\,\unit{\s}									&   32.8\,\unit{\s} \\
																	&Total			  &					74.9\,\unit{\s}  		& \textbf{46.2\,\unit{\s}}						&  53.0\,\unit{\s} \\
\bottomrule
\end{tabular*}
\end{table}

For a comparison of all the different discretizations, P1--P1, Q1--Q1, P2--P1, Q2--Q1, and Q2--P1-disc., we will focus on GDSW$^\star$--RGDSW and RGDSW--RGDSW (see \Cref{figure: Iterations monolithic stationary all discs}), based on the previous findings which combination gave the best results. 

\begin{figure}[!tb]
    \begin{minipage}[t]{0.495\textwidth}%
        \includegraphics[width=\textwidth,trim={0.3cm 0.15cm 0.15cm 0.15cm},clip]{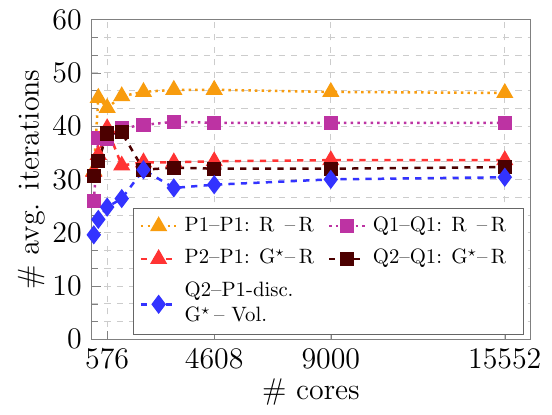}
    \end{minipage}%
    \hspace*{1pt}
    \begin{minipage}[t]{0.495\textwidth}%
        \includegraphics[width=\textwidth,trim={0.3cm 0.15cm 0.15cm 0.15cm},clip]{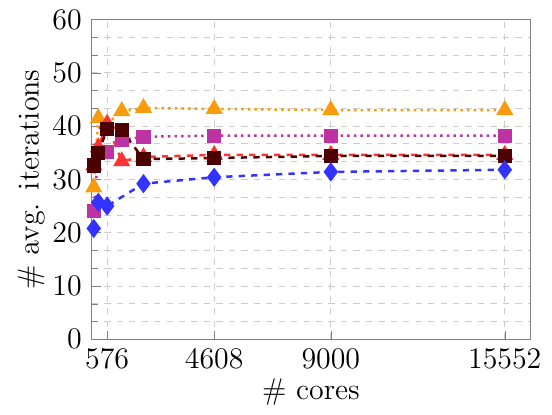}
    \end{minipage}%
	\caption{
	\textbf{Comparison of different finite element discretizations for the monolithic preconditioner.}
    Combinations of different coarse spaces for the velocity (first component) and pressure (second component). 
    G: GDSW. G\expStar{}: GDSW\expStar{}. R: RGDSW. 
    Average GMRES iteration count per Newton step. Total number of Newton iterations required to reach convergence is~5. Weak scaling test. Stationary BFS problem with configurations according to \Cref{Table:BFS_Hh_Dofs}. $\Reynum=200$. Use of Bochev--Dohrmann stabilization for P1--P1 and Q1--Q1 elements. Pressure projection included in the first level for P2--P1 and Q2--Q1 discretization. 
     \textbf{Left:} Excluding off-diagonal blocks in $\phi$ to build the coarse matrix.
    \textbf{Right:} Using full $\phi$ to build the coarse matrix. }\label{figure: Iterations monolithic stationary all discs}
\end{figure}

Due to the discontinuous pressure in Q2--P1-disc., there does not exist an interface for the pressure component. 
As a result, all coarse spaces define the same space on the pressure: the span of constant functions on each subdomain; see \Cref{subsubsec: disc}. 
As mentioned at the beginning of this section, if the pressure coarse space is too large in comparison to the velocity coarse space, the algorithm does not converge. 
We observe the same for Q2--P1-disc. if RGDSW is used for the velocity. 
The explanation may be the same as before, since, for example in case of 243~subdomains, the RGDSW pressure coarse space for a continuous P1 pressure space has 124 coarse functions (vertices of the domain decomposition), but any GDSW-type coarse space for P1-discontinuous pressure functions has 243 coarse functions. 
Thus, RGDSW--RGDSW for Q2--P1 has a smaller pressure coarse space than for Q2--P1-disc.
Consequently, the Q2--P1-disc. test results are based on the GDSW$^\star$ coarse space for the velocity.

The different configurations all yield similar good weak scaling results for up to 15\,552 processor cores. 
The stabilized P1--P1 and Q1--Q1 discretizations produce slightly higher average GMRES iteration counts than their stable counterparts. 
Excluding the off-diagonal blocks in $\phi$ or using it fully has only a minor impact on the results (compare \Cref{figure: Iterations monolithic stationary all discs} (left) and (right)). 
Corresponding to \Cref{figure: Iterations monolithic stationary all discs} (left), \Cref{figure: Average time per Newton step monolithic stationary all discs} shows the total time, consisting of setup and solve time. 
Again, all the different configurations yield similar good weak scaling results with respect to the total time for up to 15\,552 processor cores. 
Included are timings for P2--P1 and RGDSW--RGDSW, a combination that was not included in \Cref{figure: Iterations monolithic stationary all discs} but shown previously in \Cref{figure: Iterations monolithic stationary comparison}; 
the results show that for a smaller number of subdomains, the number of iterations and the total time using GDSW\expStar{} for the velocity is smaller than if RGDSW is used. 
However, the timings show that this changes for a larger number of subdomains.

\begin{figure}[!tb]
    \centering
    \includegraphics[width=0.95\textwidth,trim={0.3cm 0.15cm 0.15cm 0.15cm},clip]{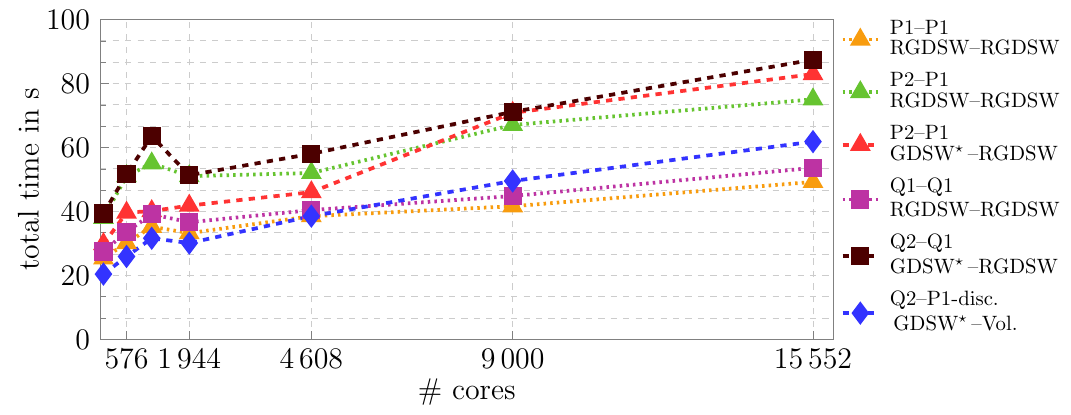}%
    \caption{\textbf{Comparison of total time (setup + solve time) for different finite element discretizations for the monolithic preconditioner.}
    Combinations of different coarse spaces for the velocity (first component) and pressure (second component). 
    Average GMRES  iteration count per Newton step. Total number of Newton iterations required to reach convergence is~5.
    Weak scaling test; see \Cref{Table: coarse procs} for the number of cores used to solve the coarse problem based on GDSW\expStar{}--RGDSW.
    Stationary BFS problem with configurations according to \Cref{Table:BFS_Hh_Dofs}. 
    $\Reynum=200$. 
    Use of Bochev--Dohrmann stabilization for P1--P1 and Q1--Q1 elements. 
    Pressure projection included in the first level for P2--P1 and Q2--Q1 discretization. 
    Use of full $\phi$. 
    }\label{figure: Average time per Newton step monolithic stationary all discs}
\end{figure}


\subsubsection{Transient Navier--Stokes problem}

For the transient backward-facing step problem, we simulate up to $t=\qty{10.0}{\s}$. 
We choose a time step size of $\Delta t=\qty{0.05}{\s}$, giving 200 time steps in total (after which an almost stationary solution is reached). 
We omit the pressure projection in the transient case, as it only slightly influences linear convergence but increases the compute time. 
We use the coarse space combinations that performed best in the stationary case. 
Furthermore, we focus on the construction of the coarse matrix where the off-diagonal blocks of $\phi$ are omitted. 

\begin{figure}[!tb]
    \centering
    \begin{minipage}[t]{0.495\textwidth}%
        \includegraphics[width=\textwidth,trim={0.3cm 0.15cm 0.15cm 0.15cm},clip]{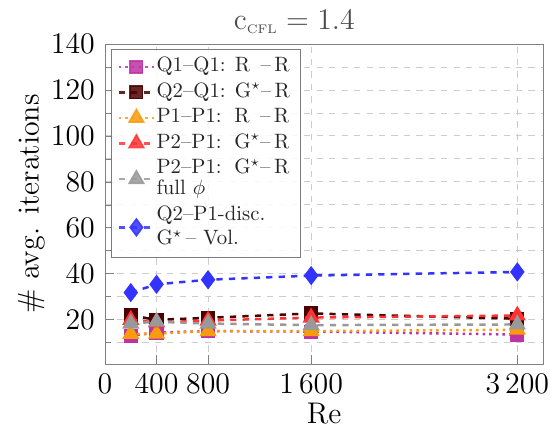}
    \end{minipage}%
    \hspace*{1pt}
    \begin{minipage}[t]{0.495\textwidth}%
        \includegraphics[width=\textwidth,trim={0.3cm 0.15cm 0.15cm 0.15cm},clip]{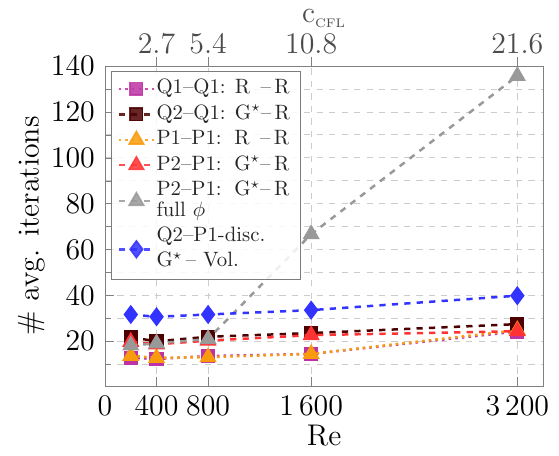}
    \end{minipage}%
    \caption{
	\textbf{Comparison of different finite element discretizations for the monolithic preconditioner for an increasing Reynolds number.} 
    For the definition of the Reynolds number, see \cref{Eq: Re BFS}. 
    Average GMRES iteration count per Newton step. 
    Transient BFS problem with configurations according to \Cref{Table:BFS_Hh_Dofs}. 
    Time step size is constant, $\Delta t = \qty{0.05}{\s}$, simulation until $t=\qty{10}{\s}$, which amounts to 200 time steps in total. 
    Use of Bochev--Dohrmann stabilization for P1--P1 and Q1--Q1 elements. 
    243 processor cores. 
    Combinations of different coarse spaces for the velocity (first component) and pressure (second component). 
    G: GDSW. G\expStar{}: GDSW\expStar{}. R: RGDSW. 
    No use of pressure projection. 
    Exclude off-diagonal blocks in $\phi$ except for one case (see legend). 
    For detailed results, see \Cref{Table:Mono RE Results table}. 
    \textbf{Left:} Increasing Reynolds number with constant $v_{\max}=\qty[per-mode=symbol]{1}{\cm\per\s}$. 
	The viscosity $\nu$ is decreased from \qty[per-mode=symbol]{0.01}{\cm^2\per\s} for $\Reynum=200$ to $\nu = \qty[per-mode=symbol]{6.25e-4}{\cm^2\per\s}$ for $\Reynum=3\,200$. 
	The CFL number is~1.4. 
	Average number of Newton iterations is between 3.0 and~3.1.
	\textbf{Right:} Increasing Reynolds number with constant $\nu=\qty[per-mode=symbol]{0.01}{\cm^2\per\s}$. 
	The velocity $v_{\max}$ is increased from \qty[per-mode=symbol]{1}{\cm\per\s} for $\Reynum=200$ to \qty[per-mode=symbol]{16}{\cm\per\s} for $\Reynum=3\,200$. 
	The CFL number increases with higher velocity. 
	Average number of Newton iterations is between 3.0 and~4.1.}\label{figure: Iterations monolithic transient all discs}
\end{figure}

In \Cref{figure: Iterations monolithic transient all discs} we investigate how the average iteration count of the linear solver (over all Newton and time steps) depends on the Reynolds number; see \Cref{Table:Mono RE Results table} for details. 
See \cref{Eq: Re BFS} for the Reynolds number definition of the BFS problem. 
The Reynolds number can be increased in two ways: either by decreasing the viscosity~$\nu$ (see \Cref{figure: Iterations monolithic transient all discs} (left)) or by increasing the maximum velocity (see \Cref{figure: Iterations monolithic transient all discs} (right)). 
The solver and preconditioner combination is very robust with respect to the Reynolds number. 

Both approaches of adjusting the Reynolds number yield similar results for the monolithic preconditioner if the off-diagonal blocks of $\phi$ are omitted; 
this is not the case if the full $\phi$ is used or for the block preconditioners PCD and SIMPLEC, as we will see in \Cref{subsect:results:comparison}.
When the CFL number increases (\Cref{figure: Iterations monolithic transient all discs} (right)), the approach using the full $\phi$ is not robust. 
However, rebuilding the coarse basis instead of reusing it will improve these results such that the number of iterations is comparable to the other coarse space variants. 
Nonetheless, since reusing the coarse basis is computationally more efficient and omitting the off-diagonal blocks proved to be more robust, we will continue to use that strategy. 
It may be, however, that in different circumstances rebuilding the coarse basis $\phi$ in each Newton step and not discarding the off-diagonal blocks is advantageous; this needs to be further investigated. 

The insight that the CFL number is the relevant metric in the considered case of varying the Reynolds number in different ways is also reflected by considering the dimensionless form of the Navier--Stokes equations. 
There, $\mathbf{u}$ is transformed to $\mathbf{u}^{*}=\mathbf{u}/u_{\text{ref}}$ via a reference velocity (for example, the maximum velocity at the inlet). 
The relevant influence on the solver by means of changing the Reynolds number via the velocity or the viscosity is then found in the transformed time step $\Delta s=u_{\text{ref}} \Delta t/L$, where $L$ is a characteristic length scale (see \cite[p.~334 and p.~412]{Elman_fluid_book} for details). 
This is comparable to a CFL number of the type $u_{\text{ref}} \Delta t/\Delta x$, since $L$ and $\Delta x$ are fixed.

We remark that using a high CFL number may be inadvisable, even for stable time discretizations, depending on the desired accuracy of the numerical solution. 
It is, however, not a focus of this work to address this but to analyze the robustness of the solvers with respect to different parameters.


\subsection{Block preconditioners}\label{subsect:block preconditioners}

We conduct similar tests for the block preconditioners as for the monolithic preconditioner. 
A description of the used recycling strategies for block preconditioners is given in \Cref{Subsubsection: Setup}.
First, the stationary BFS problem is used to analyze the different block preconditioners with respect to weak scaling. 
Then, weak scaling and the robustness with respect to an increasing Reynolds number is analyzed for the transient BFS problem.

In the stationary case, we will approximate the inverse of $M_p$ arising in the PCD preconditioner with the diagonal approximation that is denoted \textit{AbsRowSum}; see \Cref{Teko: PCD}. 
This was advantageous to approximating the inverse of $M_p$ with a Schwarz method. 
For the transient case though, approximating $M_p$ with a two-level overlapping Schwarz method proved to be best. 
The inverse of $M_u$ arising in the LSC-type preconditioners is also approximated via the \textit{AbsRowSum} diagonal scaling; see \Cref{sect:teko:lsc}.


\subsubsection{Stationary Navier--Stokes problem}

Similar to \Cref{Subsubsection:Mono_Stationary} we use Test~Case~1 with the backward-facing step in three dimensions. 
First we compare in \Cref{table: Block Comp detailed case stationary} the available block preconditioners from \texttt{Teko} and a diagonal block preconditioner implemented in the \texttt{FEDDLib}.
The diagonal block preconditioner approximates the Schur complement with $-\frac{1}{\nu} M_p$; 
the inverse of the approximate Schur complement and the inverse of $F^{-1}$ are then approximated with the Schwarz method.  
The BFBt method is defined by \cref{Eq: S_BFBT}, as it is a variant of the LSC block preconditioner. 

\begin{table}[!tb]
\caption{\textbf{Comparison of block preconditioners for different Reynolds numbers for stationary BFS problem.} 
P2--P1 discretization with $H/h=9$, 1\,125 cores. 
(BC--2) is applied in PCD and LSC$_{A_p}$. 
RGDSW--RGDSW coarse spaces for the velocity (first component) and pressure (second component). 
Average GMRES iteration count per Newton step and number of Newton steps in  parentheses.
}\label{table: Block Comp detailed case stationary}
\begin{tabular*}{\linewidth}{@{} LRRRRRR@{} }
\toprule
Preconditioner	& Diagonal		& BFBt 	 & LSC	 	  &		LSC$_{A_p}$  & PCD & SIMPLEC	 \\
\midrule
$\Reynum=\phantom{0}20$ 				& 	127 (4)	 	&  180 (4)		& 	169 (4)		&	 	 	154 (4)						& 	 	82 (4)		&	 278 (4)	 	\\
$\Reynum=200$ 				&	266 (5) 		 	&  292 (5)		& 	 257 (5)		& 			207 (5)						& 	 	103 (5)		&	 316	 (5)	\\ \bottomrule
\end{tabular*}
\end{table}

The SIMPLE preconditioner is known to not perform well in the stationary case; see, for example, \cite[Section 4.2]{Cyr_teko_stabilized}, \cite[Section 3]{Elman_LSC_2006}, or \cite[Section 8.3] {prec_nav_stokes}. 
Only the PCD block preconditioner delivers moderate results. 
The results for LSC may improve by using an additional scaling, like the ones introduced in \cite{Elman_BC_Navier_Stokes,Prec_Twophase_NS_2019}; 
however, we did not see an improvement with the scaling in \cite{Prec_Twophase_NS_2019} for our setting, and an implementation of the scaling matrix from \cite{Elman_BC_Navier_Stokes} is not available to us -- see \cite{Elman_BC_Navier_Stokes} for a derivation of a scaling matrix in a two-dimensional setting. 
Consequently, and particularly later when we compare monolithic with block preconditioners, we will focus mainly on the PCD preconditioner and at times show results based on the SIMPLEC preconditioner for reference.

\begin{table}[!tb]
\caption{\textbf{Comparison of boundary conditions in PCD block preconditioner for stationary BFS problem.} 
For different boundary conditions applicable to PCD, see \Cref{Table: PCD BC-Table}. 
P2--P1 discretization with $H/h=9$, 243 cores, and $\Reynum=200$.
RGDSW--RGDSW coarse space combination for the velocity and Schur complement components. 
Average GMRES iteration count per Newton step.
}\label{table: Block PCD BC stationary}
\begin{tabular*}{\linewidth}{@{} LRRR@{} }
\toprule
Boundary condition strategy		& (BC--1)		& (BC--2)  & (BC--3)		 \\
\midrule
\#\,avg. itererations			&	 $> 500$		 	&  	87.7			& 	 	93.7	\\ \bottomrule
\end{tabular*}
\end{table}

In \Cref{sec: BlockSchwarz} we showed how the two-level overlapping Schwarz preconditioners are applied to block systems. 
The arising inverses $\hat{F}^{-1}$ and $\hat{S}^{-1}$ of the block preconditioners in \Cref{section: Block Preconditioner} are approximated with different Schwarz preconditioners. 

For the stationary case, using the boundary condition strategy (BC--2) for PCD proved to be slightly better than the default (BC--3); see \Cref{table: Block PCD BC stationary}. 
Consequently (BC--2) will be used for the PCD preconditioner for the stationary BFS problem, but we will revert to the default (BC--3) for transient simulations. 

We select the RGDSW--RGDSW coarse space combination as the default strategy: 
\Cref{table: Block Comp PCD detailed case stationary} in the appendix shows that different coarse space combinations deliver similar results if the Reynolds numbers is varied by changing the viscosity. 
Furthermore, the setup and solve times are comparable; 
mixing different coarse spaces for PCD, for example GDSW\expStar{}--RGDSW, does not improve the average iteration count as in the monolithic case. 
Note that RGDSW has the lowest setup cost, since it has the smallest coarse space dimension. 
Especially for a higher number of processor cores, this will become visible; cf. \Cref{figure: Iterations mono vs block stationary} (right). 

Due to the decoupling of the sub problems via the block-preconditioning approach, the choice of the coarse space for the two-level overlapping Schwarz method is then secondary as indicated by the results in \Cref{figure: Iterations Block transient} and \Cref{table: Block Comp PCD detailed case stationary}.
These results also show that the PCD block preconditioner is robust once the boundary conditions are correctly specified (cf. \Cref{Sec: PCD Boundary}). 

\begin{figure}[!tb]
    \begin{minipage}[b]{0.495\textwidth}%
        \centering
        \includegraphics[width=\textwidth,trim={0.3cm 0.15cm 0.15cm 0.15cm},clip]{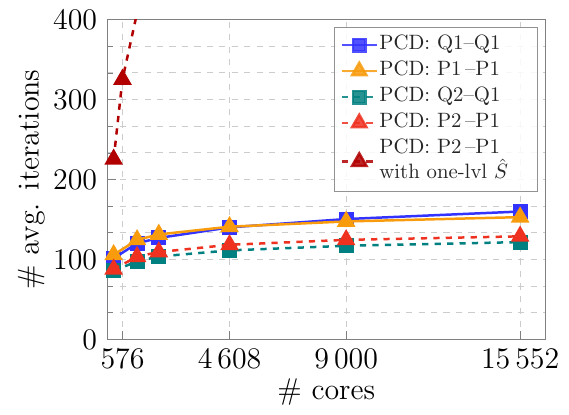}
    \end{minipage}%
    \hspace*{\fill}
    \begin{minipage}[b]{0.495\textwidth}%
        \centering
        \includegraphics[width=\textwidth,trim={0.3cm 0.15cm 0.15cm 0.15cm},clip]{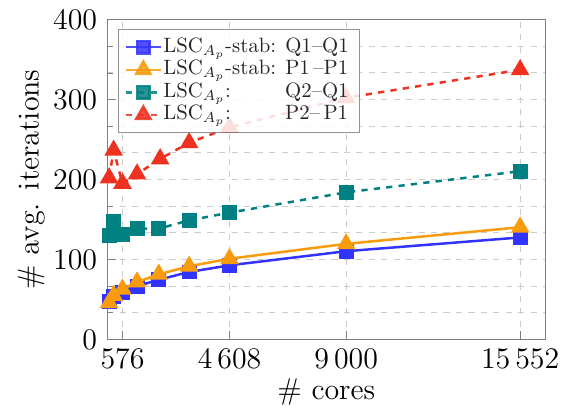}
    \end{minipage}%
    \caption{
    \textbf{Comparison of different finite element discretizations for block-triangular preconditioners.}
    Use of RGDSW--RGDSW coarse space for the velocity (first component) and for Schur complement components (second component) for both PCD and LSC (except for one result using PCD; see legend). 
    (BC--2) boundary condition strategy used for PCD and LSC. 
    Average GMRES iteration count per Newton step. 
    Total number of Newton iterations required to reach convergence is~5.
    Weak scaling test. 
    Stationary BFS problem with configurations according to \Cref{Table:BFS_Hh_Dofs}. 
    $\Reynum=200$. 
    Use of Bochev--Dohrmann stabilization for P1--P1 and Q1--Q1 elements. 
    No use of pressure projection. 
    \textbf{Left:} PCD block preconditioner. 
    One-level Schwarz approximation for $\hat{S}^{-1}$ denoted as one-lvl $\hat{S}$. 
    \textbf{Right:} LSC$_{A_p}$ block preconditioner; see \Cref{subsection: LSC}.
    }\label{figure: Iterations block stationary pcd}
\end{figure}

In \Cref{figure: Iterations block stationary pcd} weak scaling results up to 15\,552 processor cores are shown for the block-triangular preconditioners PCD and LSC. 
Most results are based on the application of a two-level Schwarz method to approximate $\hat{F}^{-1}$. 
For reference, we also show a combination that uses a one-level Schwarz approximation for $\hat{S}^{-1}$. 
We test the inf-sup stable P2--P1 and Q2--Q1 and the stabilized P1--P1 and Q1--Q1 elements but omit the Q2--P1-disc. element, since the implementation of the discrete pressure Laplace operator $A_p$ cannot be used for discontinuous elements. 
In that case, a different type of operator is required for PCD and LSC$_{A_p}$; 
its implementation is not available to us, and its definition is not straight forward; see \cite[p.~368--370]{Elman_fluid_book}.

For PCD we observe weak scalability for all discretizations in \Cref{figure: Iterations block stationary pcd}~(left). 
Using only a one-level overlapping Schwarz preconditioner for the Schur complement proves to be inadequate; the average iteration count exceeds 400. 

The results using LSC are less clear; see \Cref{figure: Iterations block stationary pcd}~(right) and \Cref{figure: Iterations LSC stationary}.
LSC$_{\text{stab},A_p}$, LSC$_{\text{stab}}$, and PCD show comparably good scaling results for the stabilized Q1--Q1 and P1--P1 elements; see \Cref{figure: Iterations block stationary pcd,figure: Iterations LSC stationary}.
However, unstabilized LSC and LSC$_{A_p}$ for stable discretizations do not scale well; 
for a P2--P1 discretization, both exceed 300 average GMRES iterations.
Using Q2--Q1 elements, the convergence for LSC$_{A_p}$ is significantly better than for P2--P1 but still significantly worse than for equal-order elements. 
Both discretizations show the $h$ sensitivity of the LSC preconditioner; compare, for example, with \cite[Table 4]{Prec_Twophase_NS_2019}.
Only the stabilized LSC variants are close to reaching an asymptotically constant number of GMRES iterations.
See also \cite{prec_nav_stokes} (e.g., Table~11) for some related results of the LSC preconditioner and the stationary BFS problem or \cite{Prec_Twophase_NS_2019} for a comparison of PCD and LSC preconditioners for a lid-driven cavity problem using Q2--Q1 elements.



\subsubsection{Transient Navier--Stokes problem}

As the PCD preconditioner gave consistently good results in the previous section, we focus our analysis on it for the transient backward-facing step problem but also show results for SIMPLEC as a reference. 
We iterate for 50 time steps with the step size $\Delta t=\qty{0.02}{\s}$. 

As stated in section \Cref{Sec: PCD Boundary}, different boundary conditions applied to the Schur complement components in the PCD preconditionier can affect the performance. 
Unlike for the stationary case in \Cref{table: Block PCD BC stationary}, where strategy (BC--2) worked best and (BC--1) failed to converge in a decent number of iterations, for the transient case, all three strategies perform equally well; see \Cref{Block PCD Comp detailed case transient} in the appendix.
Thus, we will revert to the default setting of (BC--3) for transient simulations.

\begin{figure}[!tb]
    \centering
    \begin{minipage}[t]{0.495\textwidth}%
        \centering
        \includegraphics[height=4.75cm,trim={0.3cm 0.15cm 0.15cm 0.15cm},clip]{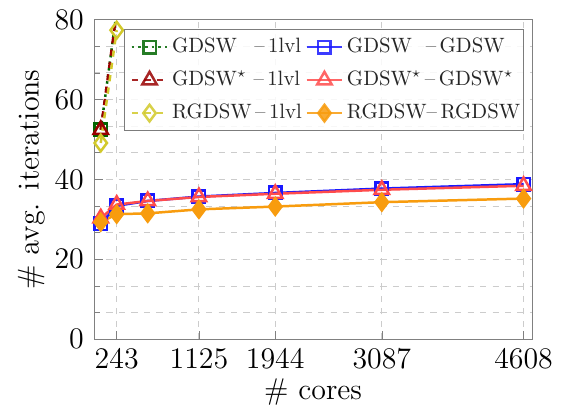}
    \end{minipage}%
    \begin{minipage}[t]{0.495\textwidth}%
        \centering
        \includegraphics[height=4.75cm,trim={0.3cm 0.15cm 0.15cm 0.15cm},clip]{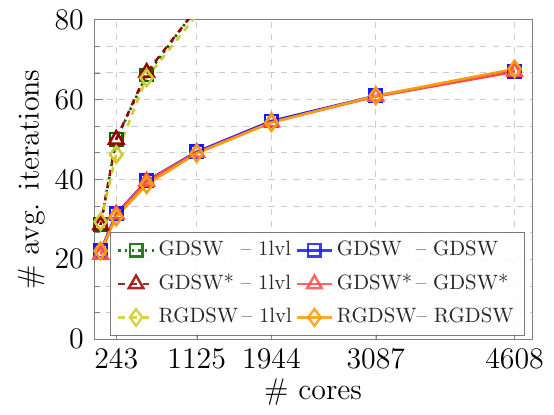}
    \end{minipage}%
    \caption{
\textbf{Comparison of different coarse space combinations for the velocity (first component) and Schur complement components (second component) for P2--P1 discretization for block preconditioners.} 
    Weak scaling test for different combinations of one- and two-level approaches. 
    Average GMRES iteration count per Newton step over all time steps. Average number of Newton iterations per time step required to reach convergence is~3.
    Transient BFS problem with $H/h=9$. 
    Time step $\Delta t=\qty{0.02}{\s}$, simulation until $t=\qty{1.0}{\s}$, which amounts to 50 time steps. 
    $\Reynum=200$.
    The CFL number varies due to the increasing mesh resolution and constant time step size. 
    The maximum CFL number is 1.44 for 4\,608 cores.
   \textbf{Left:} Pressure convection--diffusion (PCD) block preconditioner with (BC--3) boundary condition strategy for the Schur complement setup. 
   \textbf{Right:} SIMPLEC block preconditioner.
     }
     \label{figure: Iterations Block transient}
\end{figure}

We analyze the weak scaling properties of the PCD and SIMPLEC block preconditioners with \Cref{figure: Iterations Block transient}. 
As for the stationary problem, using a one-level method for the pressure component is infeasible due to very high iteration counts. 
The other shown variants of PCD and SIMPLEC qualitatively perform similarly as in the stationary case; 
PCD is weakly scalable and does significantly better than SIMPLEC.

\subsection{Comparison of specific monolithic and block preconditioners}
\label{subsect:results:comparison}

For the comparison of monolithic and block preconditioners, we will only address the configurations that gave the most robust results in the previous sections. 
Furthermore, we focus mostly on PCD as a block preconditioner but show a number of SIMPLEC results for reference. 
We contrast the weak scaling behavior for the backward-facing step and the strong scaling behavior for the realistic artery. 
We also study the preconditioners' robustness by increasing the Reynolds number. 
This is especially relevant for physically realistic scenarios. 
In \cite{DEPARIS_unsteady_NS_PCD} similar aspects relevant to hemodynamics simulations were analyzed; 
among other tests, the authors compared the SIMPLE and PCD block preconditioners for varying viscosities and Reynolds numbers.


\subsubsection{Test case 1: Backward-Facing Step}

\paragraph{Stationary simulation}

\begin{figure}[!tb]
    \centering
    \begin{minipage}[b]{0.495\textwidth}%
        \centering
        \includegraphics[height=4.6cm,trim={0.3cm 0.15cm 0.15cm 0.15cm},clip]{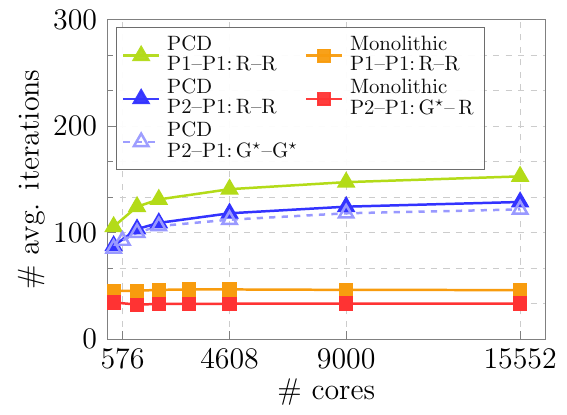}
    \end{minipage}%
    \begin{minipage}[b]{0.495\textwidth}%
        \centering
        \includegraphics[height=4.6cm,trim={0.3cm 0.15cm 0.15cm 0.15cm},clip]{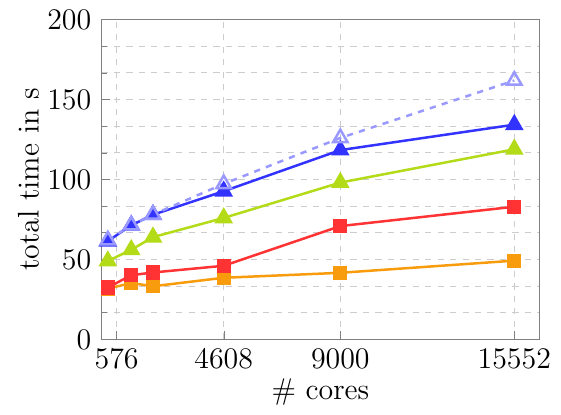}
    \end{minipage}%
    \caption{\textbf{Weak scaling comparison of block and monolithic preconditioners.} 
    Use of GDSW\expStar{} (G\expStar{}) and RGDSW (R) coarse spaces for the velocity (first component) and pressure (second component). 
    Stationary BFS problem with configurations according to \Cref{Table:BFS_Hh_Dofs}. 
    $\Reynum=200$. 
    Use of Bochev--Dohrmann stabilization for P1--P1 elements. 
    Pressure projection included in the first level for monolithic preconditioner with P2--P1 discretization. 
    Exclude off-diagonal blocks in $\phi$ for monolithic preconditioner. 
    Boundary condition strategy for PCD is (BC--3). 
    \textbf{Left:} Average GMRES iteration count per Newton step. Total number of Newton iterations required to reach convergence is~5. 
    \textbf{Right:} Total time in seconds consisting of setup and solve time.
     }
	\label{figure: Iterations mono vs block stationary}
\end{figure}

In \Cref{figure: Iterations mono vs block stationary} the average iteration count and total time for the PCD and monolithic preconditioner is shown. 
Since the average iteration count of the monolithic preconditioner is significantly lower, the total time is much lower as well. 
Despite a slightly lower number of iterations, the RGDSW--RGDSW combination for PCD performs better in case of the P2--P1 discretization and a large number of processor cores compared to the GDSW\expStar--GDSW\expStar coarse space combination, which can be attributed to the smaller cost associated with the coarse problem.

It stands out that stabilized P1--P1 elements lead to smaller compute times compared to P2--P1, despite a larger iteration count and a slightly larger number of degrees of freedom (cf. \Cref{Table:BFS_Hh_Dofs}). 
This is probably due to a more favorable sparsity pattern. 
However, since a P2--P1 discretization should give a better finite element approximation, one cannot infer from the results that P1--P1 elements are preferrable.


\paragraph{Transient simulation}

\begin{figure}[!tb]
    \begin{minipage}[b]{0.48\textwidth}%
        \centering
        \includegraphics[height=4.55cm,trim={0.3cm 0.15cm 0.15cm 0.15cm},clip]{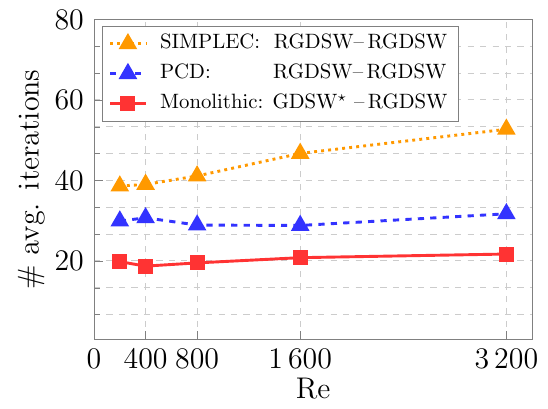}
    \end{minipage}%
    \hspace*{\fill}
    \begin{minipage}[b]{0.51\textwidth}%
        \centering
        \includegraphics[height=4.88cm,trim={0.3cm 0.15cm 0.15cm 0.0cm},clip]{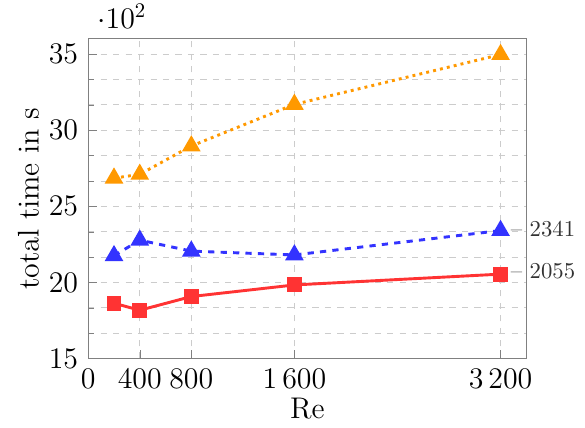}%
    \end{minipage}%
    \caption{\textbf{Comparison of block and monolithic preconditioners for increasing Reynolds number.}
     Use of GDSW\expStar{} and RGDSW coarse spaces for the velocity (first component) and pressure (second component). 
    Transient BFS problem with P2--P1 discretization and $H/h=9$. 
    Increasing Reynolds number due to decreasing viscosity. 
    No pressure projection is used. 
    Exclude off-diagonal blocks in $\phi$ for monolithic preconditioner. 
    Computation of 200 time steps with time step size $\Delta t = \qty{0.05}{\s}$, simulation until $t=\qty{10}{\s}$. 
    CFL number is constant and equal to 1.35. 
    Computations on 243 cores. 
    Boundary condition strategy in PCD is (BC--3). 
    See \Cref{BlockTransientComp_PCD_LSC viscosity-Reynolds} for detailed numbers corresponding to the figure. 
    \textbf{Left:} Average GMRES iteration count per Newton step per time step. Average number of Newton iterations per time step required to reach convergence is~3.
    \textbf{Right:} Total time in seconds consisting of setup and solve time.
    }\label{figure: Iterations mono block transient Re}
\end{figure}

The transient backward-facing step simulation consists of 200 time steps with $\Delta t=\qty{0.05}{\s}$ until $t=\qty{10}{\s}$ is reached. 
We restrict ourselves to a P2--P1 discretization and the monolithic preconditioner with the GDSW\expStar{}--RGDSW or RGDSW--RGDSW coarse space combination; 
PCD and SIMPLEC always use RGDSW--RGDSW. 
For the comparison of monolithic and block preconditioners, we will consider three aspects: 
Firstly, we will compare the robustness with respect to an increasing Reynolds number. 
The Reynolds number is varied by adjusting the viscosity. 
As we have seen before, changing the velocity instead is more challenging, as this also changes the CFL number.
In \Cref{figure: Iterations mono block transient Re} (left), the average linear iteration count is shown for an increasing Reynolds number. 
Detailed iteration counts and timings are included in the appendix in \Cref{BlockTransientComp_PCD_LSC viscosity-Reynolds}. 
The PCD and monolithic preconditioners are robust with respect to the increase of the Reynolds number, but also SIMPLEC shows only a small increase in the number of iterations. 
Despite using a larger coarse space for the velocity component, the monolithic preconditioner requires the least amount of CPU time; see \Cref{BlockTransientComp_PCD_LSC viscosity-Reynolds} and \Cref{figure: Iterations mono block transient Re} (right).

\begin{figure}[!tb]
    \centering
    \begin{minipage}[b]{0.495\textwidth}%
        \centering
        \includegraphics[width=\textwidth,trim={0.3cm 0.15cm 0.15cm 0.15cm},clip]{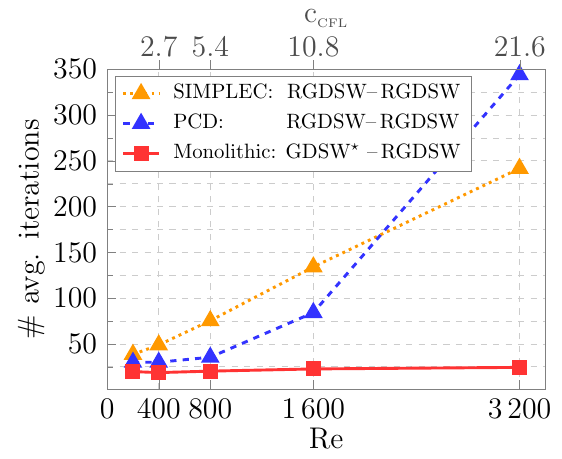}
    \end{minipage}%
    \hspace*{1pt}
    \begin{minipage}[b]{0.495\textwidth}%
        \centering
        \includegraphics[width=\textwidth,trim={0.3cm 0.15cm 0.15cm 0.15cm},clip]{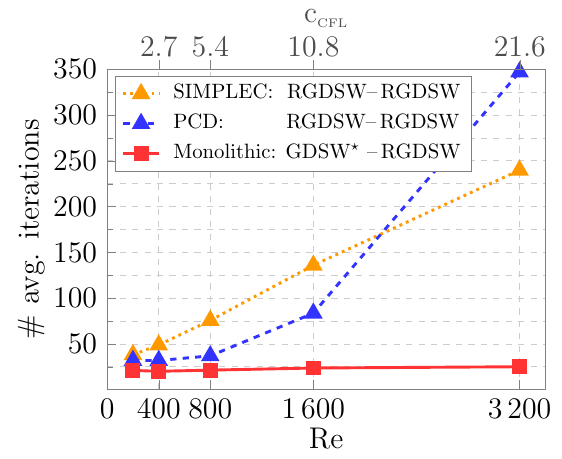}
    \end{minipage}%
    \caption{\textbf{Comparison of block and monolithic preconditioners for increasing Reynolds and CFL number.} 
     Use of GDSW\expStar{} and RGDSW coarse spaces for the velocity (first component) and pressure (second component). 
    Transient BFS problem with P2--P1 discretization and $H/h=9$. 
    Increasing Reynolds number due to increasing velocity. 
    No pressure projection is used. 
    Exclude off-diagonal blocks in $\phi$ for monolithic preconditioner. 
    Constant time step size $\Delta t = \qty{0.05}{\s}$, simulation until $t=\qty{10}{\s}$. 
    CFL condition increases with increasing velocity. 
    Boundary condition strategy of PCD is (BC--3). Average GMRES iteration count per Newton step per time step. 
    Computations of all methods with reuse of coarse basis (CB) (\textbf{Left}) and without reuse (\textbf{Right}); see \Cref{Subsubsection: Setup}.
    Note that the two plots look identical, since the reuse of the coarse basis (CB) has a very minor effect on the iteration count in this test setting. 
    }
    \label{figure: Iterations CFL comparison}
\end{figure}

Secondly, we will consider robustness with respect to an increasing CFL number. 
In \Cref{figure: Iterations mono block transient Re} the CFL number remained constant. 
In \Cref{figure: Iterations CFL comparison} we raise the Reynolds number by increasing the maximum inflow velocity, which increases the CFL number, since the time step size remains constant. 
For the monolithic preconditioner, the average iteration count per Newton step over all time steps is still almost constant over the range 1.35 to 21.6 of CFL numbers. 
PCD and SIMPLEC, however, show a large increase of the average iteration count with the increasing CFL number. 
A sensitivity of PCD and SIMPLEC with respect to the CFL number was also observed in \cite[Fig.~9]{Cyr_teko_stabilized}.
Not reusing the coarse basis, as it was stated in \Cref{Subsubsection: Setup}, has almost no impact on the results; compare \Cref{figure: Iterations CFL comparison} (left) with (right). 
This justifies the reuse of the coarse basis throughout the simulation process to save setup time.

Thirdly, we test the weak scaling behavior; see \Cref{Table:MonoStationary_RuntimeCompP2_Re1600}. 
Due to the increasing number of processor cores, the mesh resolution increases and in turn also the CFL number. 
For 4\,608 processor cores, the CFL number is still moderate with $\CFL=3.6$. 
Based on the results in \Cref{figure: Iterations CFL comparison}, this should not impact the average iteration count by much. 
All preconditioners, except for SIMPLEC, perform well. 
The results for 4\,608 cores and SIMPLEC were not computed due to the high computational cost.

\begin{table}[!tb]
\caption{\textbf{Comparison of weak scaling between block and monolithic preconditioners.} 
Transient BFS problem with P2--P1 discretization and $H/h=9$. 
$\Reynum=3\,200$ with $\nu=\qty[per-mode = symbol]{6.25e-4}{\cm^2\per\s}$ and inflow velocity $v_{\max}=\qty[per-mode=symbol]{1}{\cm\per\s}$. 
Average number of Newton iterations per time step required to reach convergence is~3.
Due to the increasing number of subdomains, the mesh resolution increases and in turn also the CFL number $\CFL$. 
No pressure projection is used. 
PCD boundary condition strategy is (BC--3). 
SIMPLEC and PCD use the RGDSW--RGDSW coarse space combination.
Exclude off-diagonal blocks in $\phi$ in monolithic preconditioner. 
Detailed results for 243 cores and varying Reynolds numbers are given in \Cref{BlockTransientComp_PCD_LSC viscosity-Reynolds}. The results for 4\,608 cores and SIMPLEC were not computed due to the high computational cost. 
In the case of 4\,608 processor cores -- to achieve the best results -- 8 cores are used to solve the coarse problems for PCD, 12 for the monolithic preconditioner and RGDSW--RGDSW, and 16 for the monolithic preconditioner and GDSW\expStar{}--RGDSW.}
\label{Table:MonoStationary_RuntimeCompP2_Re1600}
\begin{tabular*}{\linewidth}{@{} LLRRRRR@{} }
\toprule
			&	& & \multicolumn{2}{c}{Block preconditioner} 			&  \multicolumn{2}{c}{Monolithic preconditioner} 	\\ \cmidrule{4-5} \cmidrule{6-7}
\#\,cores			&		$\CFL$ 										&											 & \footnotesize SIMPLEC 			& \footnotesize PCD 		& \footnotesize RGDSW--RGDSW & \footnotesize GDSW\expStar{}--RGDSW \\ \midrule 
	\multirow{4}{*}{\phantom{0\,}243} & \multirow{4}{*}{1.4}
																					&\#\,avg.\,iter. & 	52.7			&  31.7 			&		35.0			&	 21.7\\ \cmidrule{3-7}
																	&				&Setup time 	&  595\,\unit{\s}		&  580\,\unit{\s} 		&		785\,\unit{\s}	 	&	787\,\unit{\s}\\
																	&				&	Solve time 	&  2\,901\,\unit{\s}	&  1\,761\,\unit{\s}		&		2\,008\,\unit{\s} 	&	1\,268\,\unit{\s}	\\
																	&				&Total			 	 &  3\,496\,\unit{\s}	&  2\,341\,\unit{\s} 		&		2\,793\,\unit{\s} 	& \textbf{2\,055\,\unit{\s}}\\ \midrule
\multirow{4}{*}{1\,125}	& \multirow{4}{*}{2.3}
																					&\#\,avg.\,iter. &  96.8		&  41.5   	       &		38.4			&  20.5 \\ \cmidrule{3-7}
																	&				&Setup time 	&  709\,\unit{\s}		&  684\,\unit{\s}		    &		904\,\unit{\s}		&	935\,\unit{\s}	\\
																	&				&	Solve time 	&  6\,882\,\unit{\s}	&  2\,924\,\unit{\s} 		&		2\,788\,\unit{\s}	&	1\,512\,\unit{\s}	\\
																	&				&Total				&  7\,591\,\unit{\s}	&  3\,608\,\unit{\s}		&		3\,692\,\unit{\s}	& \textbf{2\,447\,\unit{\s}}	\\ \midrule
	\multirow{4}{*}{4\,608} & \multirow{4}{*}{3.6}				
																					&\#\,avg.\,iter.	&   ---				&  43.7	       &		43.0			&  19.0\\ \cmidrule{3-7}
																	&				&Setup time 	& 		---			&  745\,\unit{\s}		    &		994\,\unit{\s}		&  1\,179\,\unit{\s}\\
																	&				&	Solve time 	& 		---			&  3\,567\,\unit{\s} 		&		 3\,642\,\unit{\s}	&  1\,664\,\unit{\s}\\
																	&				&Total				&  	---			&  4\,312\,\unit{\s}		&		4\,636\,\unit{\s}	& \textbf{2\,843\,\unit{\s}}	\\
\bottomrule
\end{tabular*}
\end{table}

Due to the low number of iterations achieved with the monolithic preconditioner and GDSW\expStar{}--RGDSW, this combination delivers by far the lowest solve time. 
Consequently, with respect to total time, it performed best.


\subsubsection{Test case 2: realistic artery -- transient problem}

\begin{figure}[!tb]
    \centering
    \begin{tikzpicture}
        \node[inner sep=0,below left] at (0,0) {\includegraphics[height=4.48cm,trim={0.3cm 0.15cm 0.15cm 0.15cm},clip]{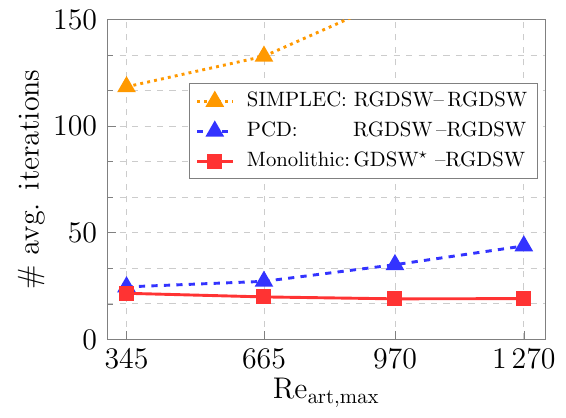}};
        \node[inner sep=0,below right] at (0.2cm,-0.07cm) {\includegraphics[height=4.3cm,trim={0.3cm 0.15cm 0.15cm 0.15cm},clip]{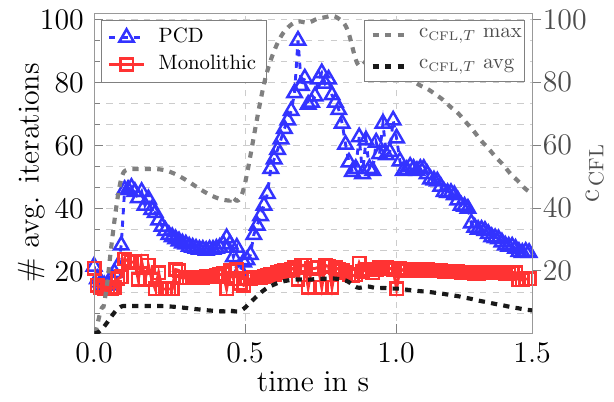}};
    \end{tikzpicture}%
  \caption{%
    \textbf{Comparison of block and monolithic preconditioners for the transient realistic artery problem (\Cref{figure2: realisic artery}) by varying the Reynolds and CFL number.}
    Transient problem with P2--P1 discretization and 4.1~million degrees of freedom. Average number of Newton iterations per time step required to reach convergence is~3.
    Computations on 416 processor cores. 
    No pressure projection is used. 
    Use of GDSW\expStar{} and RGDSW coarse spaces for the velocity (first component) and pressure (second component). 
    Excluding off-diagonal blocks of $\phi$ in monolithic preconditioner. 
    Boundary condition strategy of PCD is (BC--3). 
    \textbf{Left:} 
    Increasing Reynolds number due to increasing initial velocity $v_{\max}$ from 10 to \qty[per-mode=symbol]{40}{\cm\per\s}; see \Cref{Table: Velocity and Re in artery}. 
    \textbf{Right:} 
    Configuration with $ \Reynum_{\text{art},\max}$ = 1\,270; see \Cref{Table: Velocity and Re in artery}. 
    Average GMRES iteration count per Newton step for each time step. 
    Value is plotted every tenth time step. 
    The maximum elementwise CFL number $\CFL(t)$ is shown (light gray); 
    its maximum over all time steps is $\CFL=100$. 
    The average of elementwise CFL numbers $\CFLavg(t)$ is depicted as well (dark gray); its average over time is 11. 
    Due to the high average GMRES iteration count, the SIMPLEC preconditioner is excluded from the figure.
    }\label{figure: Iterations RE transient}
\end{figure}

To test the preconditioners for a range of inflow velocities, we adjust the time step size accordingly to obtain a fairly constant CFL number; see \Cref{Table: Velocity and Re in artery}. 
In \Cref{figure: Iterations RE transient} (left) we can observe the response of the preconditioning strategies to the maximum Reynolds number defined in \cref{Def: Re artery max}.
The PCD and monolithic preconditioner deliver robust results. 
If we take a closer look at the average number of iterations per Newton step for the different time steps, the monolithic approach is very robust to the velocity changes, since the iteration count over time is almost constant in \Cref{figure: Iterations RE transient} (right). 
The graph of the PCD preconditioner is similar to the graph of the CFL number; that is, it shows that its performance deteriorates with a larger CFL number.

\begin{figure}[!tb]
    \centering
    \begin{minipage}[b]{0.495\textwidth}%
        \centering
        \includegraphics[height=4.77cm,trim={0.3cm 0.15cm 0.15cm 0.15cm},clip]{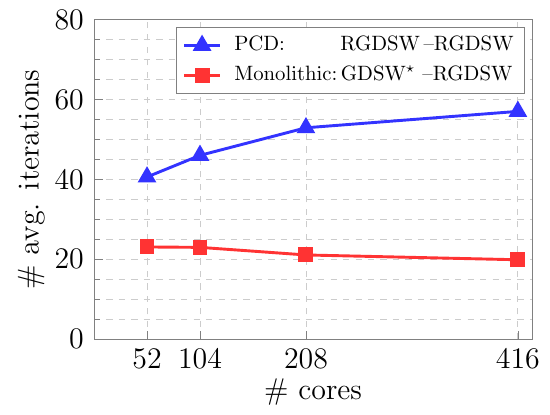}
    \end{minipage}%
    \hspace*{1pt}
    \begin{minipage}[b]{0.495\textwidth}%
        \centering
        \includegraphics[height=4.8cm,trim={0.3cm 0.15cm 0.15cm 0.15cm},clip]{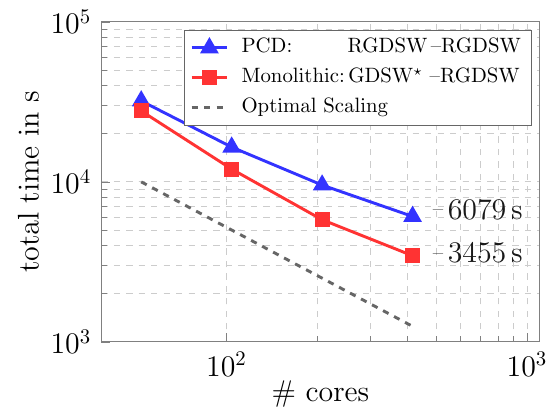}
    \end{minipage}%
    \caption{\textbf{Comparison of strong scaling of block and monolithic preconditioners for transient realistic artery problem (\Cref{figure2: realisic artery}).} 
    Use of GDSW\expStar{} and RGDSW coarse spaces for the velocity (first component) and pressure (second component). 
   Transient problem with P2--P1 discretization and 4.1~million degrees of freedom. Average number of Newton iterations per time step required to reach convergence is~3.
    Computations on 52 to 416 processor cores; see \Cref{Table: Strong Scaling Procs Coarse Solve} for the number of cores used to solve the coarse problem. 
    Time step size $\Delta t=\qty{0.002}{\s}$, initial maximum velocity of \qty[per-mode=symbol]{40}{\cm\per\s}, resulting in $\Reynum_{\text{art},\max}= 1\,270$. 
    The maximum CFL number is $\CFL=160$. 
    For the average elementwise CFL number, we have $\CFLavg(t)=28$. 
    The average of $\CFLavg(t)$ over all time steps is 18. 
    No pressure projection is used. 
    Excluding off-diagonal blocks in $\phi$. 
    Boundary condition strategy of PCD is (BC--3).
    \textbf{Left:} Comparison of average GMRES iteration count per Newton step per time step for different numbers of processor cores and different preconditioning techniques. 
    \textbf{Right}: Total time consisting of setup and solve time.
    }
\label{figure: Strong Scaling artery}
\end{figure}

Finally, we show strong scaling results for the artery in \Cref{figure: Strong Scaling artery}; the mesh has 4.1~million degrees of freedom. 
Between 52 and 416 processor cores are used, of which 2 to 8 are used to solve the coarse problem; see \Cref{Table: Strong Scaling Procs Coarse Solve}. 
The unstructured mesh partition is constructed with METIS. 
We use an initial velocity of \qty[per-mode=symbol]{40}{\cm\per\s} but increase the time step with respect to \Cref{Table: Velocity and Re in artery} to $\Delta t = \qty{0.002}{\s}$ to reach results for 52~processor cores within a reasonable time frame. 
We test the monolithic preconditioner with GDSW\expStar{}--RGDSW and the PCD block preconditioner with RGDSW--RGDSW. 
The results show strong scaling of both the PCD and the monolithic preconditioner. 
However, the monolithic preconditioner is more robust and requires significantly fewer iterations to converge, which results in a substantially lower computing time.

Based on the results in this work and the considered inflow--outflow problems, choosing the monolithic two-level overlapping Schwarz preconditioner with the GDSW\expStar{}--RGDSW coarse space combination provides a consistently robust performance with a low compute time. 
The PCD block preconditioner can achieve similarly reliable results, provided the boundary conditions of the preconditioner are chosen appropriately and the CFL number is not too high; cf., e.g., \Cref{figure: Iterations RE transient}.

\begin{table}[!tb]
\caption{\textbf{Number of processor cores assigned to the coarse problem for the results in \Cref{figure: Strong Scaling artery}.} 
    Number of coarse functions in  parentheses. 
    Coarse spaces: The monolithic preconditioner uses GDSW\expStar{}--RGDSW and the PCD block preconditioner RGDSW--RGDSW.}
\label{Table: Strong Scaling Procs Coarse Solve}
\begin{tabular*}{\linewidth}{@{} LRRRR@{} }
\toprule
\#\,cores 		& 52 			& 104 			& 208 			& 416 \\ \midrule
PCD 			& 2\phantom{0\,} (404) 	& 3 (1\,084)	& 4	(2\,771)	& 	8\phantom{0} (6\,057)\\
Monolithic 	&	2	(1\,013)	& 4 (2\,518)	& 	6	(5\,912)	& 12 (12\,975)\\
\bottomrule
\end{tabular*}
\end{table}


\section*{Acknowledgments}

Financial funding from the Deutsche Forschungsgemeinschaft (DFG) through the Priority Program 2311 ``Robust coupling of continuum-biomechanical in silico models to establish active biological system models for later use in clinical applications - Co-design of modeling, numerics and usability'', project ID 465228106, is greatly appreciated.

The authors gratefully acknowledge the scientific support and HPC resources provided by the Erlangen National High Performance Computing Center (NHR@FAU) of the Friedrich-Alexander-Universität Erlangen-Nürnberg (FAU) under the NHR project k105be. NHR funding is provided by federal and Bavarian state authorities. NHR@FAU hardware is partially funded by the German Research Foundation (DFG) -- 440719683.


\appendix


\section{Stationary Navier--Stokes problem}
\label{appendix:stationary:monolithic}
\begin{figure}[h]
    \begin{minipage}[t]{0.495\textwidth}%
        \includegraphics[width=\textwidth,trim={0.3cm 0.15cm 0.15cm 0.15cm},clip]{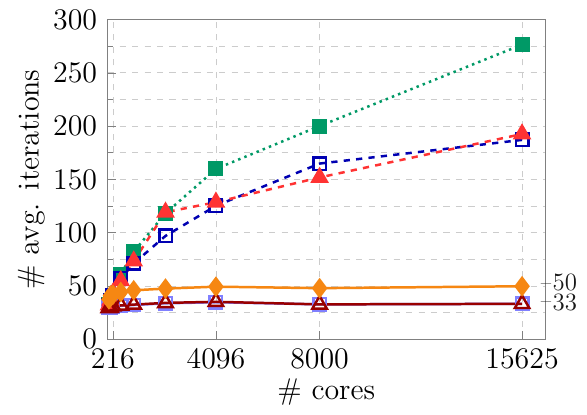}
    \end{minipage}%
    \hspace*{\fill}
    \begin{minipage}[t]{0.495\textwidth}%
        \includegraphics[width=\textwidth,trim={0.3cm 0.15cm 0.15cm 0.15cm},clip]{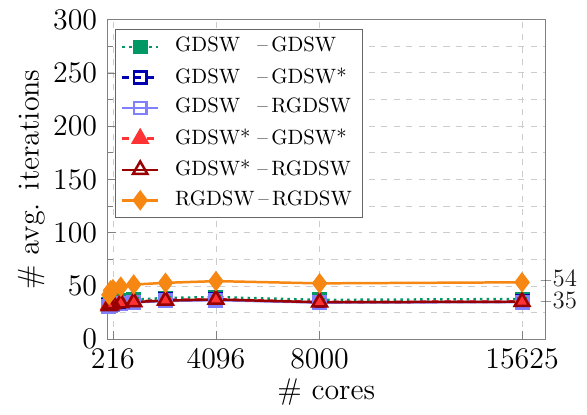}
    \end{minipage}%
    \caption{\textbf{Comparison of different coarse space combinations for the velocity (first component) and pressure (second component) for P2--P1 discretization for the monolithic preconditioner.}
    Average GMRES iteration count per Newton step. 
    Weak scaling test. 
    Stationary lid-driven cavity problem with $H/h=9$; see \Cref{appendix:stationary:monolithic} for a description. 
    $\Reynum=200$. 
    Using pressure projection in first level of monolithic preconditioner. 
    See \Cref{figure: Iterations monolithic stationary comparison} for corresponding BFS problem.
    \textbf{Left:} Excluding off-diagonal blocks in $\phi$ to build the coarse matrix.
    \textbf{Right:} Using full $\phi$ to build the coarse matrix.
    }
  \label{figure: Iterations monolithic stationary ldc}
\end{figure}

In \Cref{figure: Iterations monolithic stationary ldc} results for a stationary, three-dimensional lid-driven cavity problem are shown to ensure that the results seen before in \Cref{figure: Iterations monolithic stationary comparison} are not particularly problem specific: 
The results for the monolithic preconditioner, different combinations of coarse spaces, and the use of the pressure projection \cref{eq: Pressure Projection} qualitatively show the same behavior as before for the stationary backward-facing step problem. 

The setup uses the classic formulation, where the entire top part of a unit cube is moved at constant velocity in one direction; see, for example, \cite[Sect.~2.1]{Heinlein_2019_MonoFluidFlow}. 
We use a structured mesh of the unit cube $[0,1]^3$, decomposed into cubic subdomains. 
The subdomain resolution is $H/h=9$, giving for a P2--P1 discretization 4.1 million degrees of freedom if $6^3=216$ subdomains are used and 286.7 million degrees of freedom if $25^3=15\,625$ subdomains are used. 
The kinematic viscosity is $\nu=0.005$, giving $\Reynum=200$. 
The lid includes its boundary points (classic lid-driven cavity formulation). 
The Dirichlet condition for the velocity at the lid is $(1,0,0)$ in each node and zero elsewhere.

It was not necessary to ensure uniqueness of the pressure by fixing a pressure node or introducing a projection into GMRES; 
note that prescribing a zero mean integral via a Lagrange multiplier can be computationally inefficient in a parallel setting, as the associated vector in the system matrix is dense and couples globally across all subdomains. 

\vspace*{\fill}

\vspace*{\fill}


\vspace*{\fill}

\begin{figure}[!tb]
    \begin{minipage}[t]{0.51\textwidth}%
        \includegraphics[height=4.6cm,trim={0.3cm 0.15cm 0.15cm 0.15cm},clip]{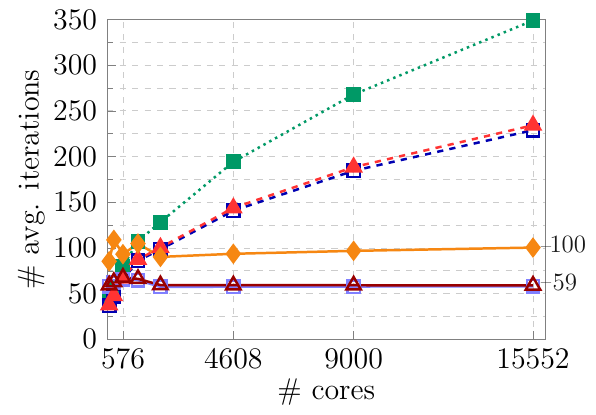}
    \end{minipage}%
    \hspace*{\fill}
    \begin{minipage}[t]{0.48\textwidth}%
        \includegraphics[height=4.6cm,trim={0.3cm 0.15cm 0.15cm 0.15cm},clip]{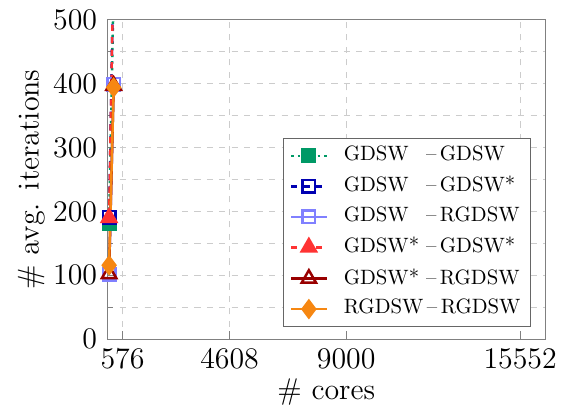}
    \end{minipage}%
    \caption{ \textbf{Comparison of different coarse space combinations for the velocity (first component) and pressure (second component) for a P2--P1 discretization and the monolithic preconditioner.}
    Average GMRES iteration count per Newton step. 
    Weak scaling test. 
    Stationary BFS problem with $H/h=9$. 
    $\Reynum=200$. 
    Omitting pressure projection; see \Cref{figure: Iterations monolithic stationary comparison} for the corresponding results if the pressure projection is used. 
    \textbf{Left:} Excluding off-diagonal blocks in $\phi$ to build the coarse matrix.
    \textbf{Right:} Using full $\phi$ to build the coarse matrix.
    }
 \label{figure: Iterations monolithic stationary comparison without pP}
\end{figure}

\vspace*{\fill}

\begin{table}[!tb]
\caption{\textbf{Comparison of different coarse spaces for PCD block preconditioner with different Reynolds numbers for P2--P1 disretization.} Combinations of different coarse spaces for the velocity (first component) and Schur complement components (second component). Stationary BFS problem with $H/h = 9$. (BC--2) is applied in PCD. 1\,125 cores. Reynolds number is increased by decreasing the viscosity. 
For the definition of the Reynolds number, see \cref{Eq: Re BFS}. 
Average GMRES iteration count per Newton step and number of Newton steps in  parentheses. 
Total time consists of setup and solve time.
}\label{table: Block Comp PCD detailed case stationary}
\begin{tabular*}{\linewidth}{@{}LLRRRR@{} }
\toprule
\multirow{2}{*}{\shortstack[l]{Coarse\\[1pt] space}}
& Velocity &GDSW 		& GDSW\expStar{} 																																		& GDSW\expStar{}  & RGDSW		 \\ 
& Pressure & GDSW		& GDSW\expStar{} 	
																																				&RGDSW  & RGDSW		 \\ 																							
\midrule
\multirow{4}{*}{$\Reynum=\phantom{0}20$ }	& \#\,avg. iterations		
																					& 		76 (4)					&		79 (4)			&		76 (4)			&		82 (4)								 	\\  \cmidrule{3-6}
												& Setup time				& 		15.0\,\unit{\s}					&		13.7\,\unit{\s}			&		15.0\,\unit{\s}			&		12.9\,\unit{\s}							 	\\  
												& Solve time				& 		35.3\,\unit{\s}					&		35.1\,\unit{\s}			&		33.8\,\unit{\s}			&		35.2\,\unit{\s}								 	\\  
												& Total 						& 		50.3\,\unit{\s}  				&		48.8\,\unit{\s}			&		48.8\,\unit{\s}			&		48.1\,\unit{\s}							 	\\
\midrule
\multirow{4}{*}{$\Reynum=200$ }	& \#\,avg. iterations	
																					& 		98 (5)					&	100 (5)				&		102 (5) 			&	103 (5)									 	\\  \cmidrule{3-6}
												& Setup time				& 		15.5\,\unit{\s}					&	14.2	\,\unit{\s}			&		13.9\,\unit{\s}				&	12.9\,\unit{\s}									 	\\  
												& Solve time				& 		57.7\,\unit{\s}					&	56.7\,\unit{\s}				&		57.2\,\unit{\s}				&	56.9\,\unit{\s}									 	\\  
												& Total 						& 	73.2\,\unit{\s}						&	70.9\,\unit{\s}				&		70.1\,\unit{\s}				&	69.8\,\unit{\s}									 	\\  
\midrule
\multirow{4}{*}{$\Reynum=400$ }	& \#\,avg. iterations	
																				
																					& 		163 (6)					&	167 (6)				&		174 (6)				&	173 (6)									 	\\  \cmidrule{3-6}
												& Setup time				& 		15.9\,\unit{\s}					&	13.8\,\unit{\s}				&		13.7\,\unit{\s}				&	13.0\,\unit{\s}									 	\\  
												& Solve time				& 		122.3\,\unit{\s}				&	120.2\,\unit{\s}			&		125.0\,\unit{\s}			&	119.8\,\unit{\s}									 	\\  
												& Total 						& 		138.2\,\unit{\s}				&	134.0\,\unit{\s}			&		138.7\,\unit{\s}			&	132.8\,\unit{\s}							 	\\  
\bottomrule
\end{tabular*}
\end{table}

\begin{figure}[!tb]
    \begin{minipage}[c]{0.51\textwidth}%
        \includegraphics[height=4.6cm,trim={0.3cm 0.15cm 0.15cm 0.15cm},clip]{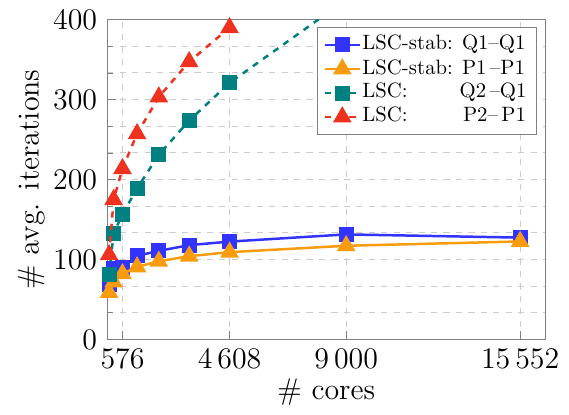}
    \end{minipage}%
        \begin{minipage}[c]{0.47\textwidth}%
      \caption{ \textbf{Comparison of different finite element discretizations for LSC block-triangular preconditioner; see \Cref{subsection: LSC}.}
    Use of RGDSW--RGDSW coarse space for the velocity (first component) and for Schur complement components (second component).
    Average GMRES iteration count per Newton step. 
    Total number of Newton iterations required to reach convergence is~5.
    Weak scaling test. 
    Stationary BFS problem with configurations according to \Cref{Table:BFS_Hh_Dofs}. 
    $\Reynum=200$. 
    Use of Bochev--Dohrmann stabilization for P1--P1 and Q1--Q1 elements. }
    \label{figure: Iterations LSC stationary}
        \end{minipage}%
\end{figure}

\vspace*{\fill}


\clearpage
\newpage
\section{Transient Navier--Stokes problem}

\vspace*{\fill}

\begin{table}[h]
\caption{
\textbf{Comparison of different finite element discretizations for the monolithic preconditioner for an increasing Reynolds number.} 
Average GMRES iteration count per Newton step per over all time steps. 
Transient BFS problem with configurations according to \Cref{Table:BFS_Hh_Dofs}. 
Time step size is constant, $\Delta t = \qty{0.05}{\s}$, simulation until $t=\qty{10}{\s}$, which amounts to 200 time steps in total. 
Use of Bochev--Dohrmann stabilization for P1--P1 and Q1--Q1 elements. 
243 processor cores. 
Combinations of different coarse spaces for the velocity (first component) and pressure (second component). 
No use of pressure projection. 
Exclude off-diagonal blocks in $\phi$. 
Results corresponding to \Cref{figure: Iterations monolithic transient all discs}.}\label{Table:Mono RE Results table}
\begin{tabular*}{\linewidth}{@{}  LLR >{\raggedleft\arraybackslash}p{1.25cm} >{\raggedleft\arraybackslash}p{1.35cm} >{\raggedleft\arraybackslash}p{1.35cm} >{\raggedleft\arraybackslash}p{1.35cm} >{\raggedleft\arraybackslash}p{1.85cm} @{}}
\toprule
\multicolumn{2}{@{}l}{Discretization} 
                                  &            & P1--P1 & P2--P1         & Q1--Q1 & Q2--Q1         & Q2--P1-disc. \\
\midrule
\multirow{2}{*}{\shortstack[l]{Coarse\\[1pt] space}}
                      & \multicolumn{2}{r}{Velocity} & RGDSW  & GDSW\expStar{} & RGDSW  & GDSW\expStar{} & GDSW\expStar{} \\
                      & \multicolumn{2}{r}{Pressure} & RGDSW  & RGDSW          & RGDSW  & RGDSW          & GDSW\expStar{} \\
\midrule
\Reynum                    & $\nu$     & $v_{\max}$ & \multicolumn{5}{c}{\#\,avg. iter. per Newton step over all time steps} \\ \midrule
\phantom{0\,}200      & 0.01      &  1         & 13.7   & 19.8           & 12.8   & 21.8           & 31.7 \\
\phantom{0\,}400      & 0.005     &  1         & 13.8   & 18.7           & 14.3   & 19.8           & 35.3 \\
\phantom{0\,}800      & 0.0025    &  1         & 14.9   & 19.5           & 14.9   & 20.7           & 37.3 \\
          1\,600      & 0.00125   &  1         & 14.9   & 20.8           & 14.7   & 22.6           & 39.1 \\
          3\,200      & 0.000625  &  1         & 15.5   & 21.7           & 13.4   & 20.3           & 40.7 \\
\midrule
\phantom{0\,}200      & 0.01      &  1         & 13.7   & 19.8           & 12.8   & 21.8           & 31.7 \\
\phantom{0\,}400      & 0.01      &  2         & 12.8   & 18.7           & 12.3   & 20.0           & 30.7 \\
\phantom{0\,}800      & 0.01      &  4         & 13.1   & 20.3           & 13.6   & 22.0           & 31.7 \\
          1\,600      & 0.01      &  8         & 14.5   & 22.7           & 14.5   & 23.6           & 33.6 \\
          3\,200      & 0.01      & 16         & 24.9   & 24.4           & 24.3   & 27.5           & 39.9 \\
\bottomrule
\end{tabular*}
\end{table}

\vspace*{\fill}

\begin{table}[h]
\caption{\textbf{Comparison of different boundary condition strategies for the PCD block preconditioner and a transient simulation.} 
Boundary conditions in PCD are defined according to \Cref{Sec: PCD Boundary}. 
Transient BFS problem with P2--P1 discretization. 
Average iteration count per Newton step over all time steps, number of Newtons steps in parentheses. 
GDSW--GDSW coarse space for the velocity (first component) and pressure (second component). 
Number of processor cores is~243. 
Time step $\Delta t=\qty{0.02}{\s}$, simulation until $t=\qty{1.0}{\s}$, which amounts to 50~time steps. 
(*) Results corresponding to \Cref{figure: Iterations Block transient}~(left).}\label{Block PCD Comp detailed case transient}
\begin{tabular*}{\linewidth}{@{} LRRR@{} }
\toprule
Boundary condition strategy & (BC--1) & (BC--2) & (BC--3) (*) \\
\midrule
\#\,avg. iterations & 33.5 (3.1) & 33.4 (3.1) & 33.4 (3.1) \\
\bottomrule
\end{tabular*}
\end{table}

\vspace*{\fill}


\clearpage
\newpage
\section{Comparison of block and monolithic preconditioners for varying Reynolds numbers}

\vspace*{\fill}

\begin{table}[h]
\caption{\textbf{Comparison of block and monolithic preconditioners for an increasing Reynolds number.} 
    Average GMRES iteration count per Newton step per time step. 
    Total time in seconds consisting of setup and solve time. 
    Transient BFS problem with P2--P1 discretization and $H/h=9$. 
    Increasing Reynolds number due to decreasing viscosity. 
    Computation of 200 time steps with time step size $\Delta t = \qty{0.05}{\s}$, simulation until $t=\qty{10}{\s}$. 
    CFL number is constant and equal to 1.35. 
    Computations on 243 cores. 
    Boundary condition strategy in PCD is (BC--3).
    Use of GDSW\expStar{} and RGDSW coarse spaces.
    Results corresponding to \Cref{figure: Iterations mono block transient Re}. 
    }
\label{BlockTransientComp_PCD_LSC viscosity-Reynolds}
\begin{tabular*}{\linewidth}{@{} LLRRRRR@{} }
\toprule
			& 		\Reynum		& 	 	  200 	      & 400 				  &  800 				& 1\,600				     & 3\,200 	 \\ \midrule
								&  $\nu$	& 	$1.0{\cdot}10^{-2}$& $5.0{\cdot}10^{-3}$& $2.5{\cdot}10^{-3}$  & 	$1.3{\cdot}10^{-3}$	& 	$6.25{\cdot}10^{-4}$	  \\ \midrule
\multirow{4}{*}{\shortstack[l]{SIMPLEC \\[2pt] \footnotesize RGDSW--\\ \footnotesize RGDSW}}	
								& \#\,avg.\,iter.			&  38.6\,(3.0)	&  	39.0\,(3.0)	&  41.1\,(3.0)		& 46.7\,(3.0)		&	52.7\,(3.0)		\\	\cmidrule{3-7}	
								&	Setup time		& 	  596\,\unit{\s}		&  	594\,\unit{\s}			&  598\,\unit{\s}			& 599\,\unit{\s}	 		&  595\,\unit{\s}			\\ 
							    &Solve	time		& 	  2\,089\,\unit{\s}		&  	2\,116\,\unit{\s}		&  2\,298\,\unit{\s}		& 2\,570\,\unit{\s}		&  2\,901\,\unit{\s}				\\
								&	Total			& 	  2\,685\,\unit{\s}		&  	2\,710\,\unit{\s}		&  2\,896\,\unit{\s}		& 3\,169\,\unit{\s}		&  3\,496\,\unit{\s}	\\	 
\midrule
\multirow{4}{*}{\shortstack[l]{PCD \\[2pt] \footnotesize RGDSW--\\ \footnotesize RGDSW}}	
								& \#\,avg.\,iter.			&  29.9\,(2.8)	&  	30.7\,(3.0)		&  28.9\,(3.0)		&  28.8\,(3.0)	&	 31.7\,(3.0)		\\	\cmidrule{3-7}	
								&	Setup time		& 	  584\,\unit{\s}		&  	578\,\unit{\s}			&  577\,\unit{\s}				& 581\,\unit{\s}	 		&  580\,\unit{\s}			\\ 
							    &Solve	time		& 	  1\,592\,\unit{\s}		&  	1\,700\,\unit{\s}		&  1\,629\,\unit{\s}			& 1\,600\,\unit{\s}	 	&  1\,761\,\unit{\s}				\\
								&	Total			& 	  2\,176\,\unit{\s}		&  	2\,278\,\unit{\s}		&  2\,206\,\unit{\s}			& 2\,181\,\unit{\s}	 	&  2\,341\,\unit{\s}	\\	 	
\midrule
\multirow{4}{*}{\shortstack[l]{Monolithic \\[2pt] \footnotesize RGDSW--\\ \footnotesize RGDSW}}	
								& \#\,avg.\,iter.			&  26.9\,(2.8)		&  	29.2\,(3.0)		&  30.5\,(3.0)	&  32.7\,(3.0)	&	 35.0\,(3.0)		\\	\cmidrule{3-7}
								&	Setup time		& 	  795\,\unit{\s}			&  	789\,\unit{\s}			&  789\,\unit{\s}			& 785\,\unit{\s}			&  785\,\unit{\s}			\\ 	
							    & Solve	time		& 	  1\,458\,\unit{\s}			&  	1\,671\,\unit{\s}		&  1\,747\,\unit{\s}		& 1\,873\,\unit{\s}		&  2\,008\,\unit{\s}				\\
								&	Total			& 	  2\,253\,\unit{\s}			&  	2\,460\,\unit{\s}		&  2\,536\,\unit{\s}		& 2\,658\,\unit{\s}	 	&  2\,793\,\unit{\s}	\\	 	\midrule

\multirow{4}{*}{\shortstack[l]{Monolithic \\[2pt] \footnotesize GDSW$^\star$--\\ \footnotesize RGDSW}}	
								& \#\,avg.\,iter.			&  19.8\,(2.8)		&  	18.7\,(2.8)		&  19.5\,(2.9)	&  20.8\,(3.0)	&	 21.7\,(3.1)		\\	\cmidrule{3-7}	
								&	Setup time		& 	  790\,\unit{\s}			&  	783\,\unit{\s}			&  786\,\unit{\s}			& 782\,\unit{\s}	 	&  787\,\unit{\s}			\\ 
							    & Solve	time	& 	 1\,074\,\unit{\s}			&  	1\,035\,\unit{\s}		&  1\,122\,\unit{\s}		& 1\,203\,\unit{\s}	 	&  1\,268\,\unit{\s}				\\
								&	Total			& \textbf{1\,864\,\unit{\s}}&  	\textbf{1\,818\,\unit{\s}}		&  \textbf{1\,908\,\unit{\s}}		& \textbf{1\,984\,\unit{\s}}	 	&  \textbf{2\,055\,\unit{\s}}	\\	 
\bottomrule
\end{tabular*}
\end{table}

\vspace*{\fill}


\clearpage
\newpage

\bibliographystyle{siamplain}
\bibliography{bibfile}

\end{document}